\documentclass[a4paper,fleqn]{cas-sc}

\usepackage[numbers]{natbib}
\def\tsc#1{\csdef{#1}{\textsc{\lowercase{#1}}\xspace}}
\tsc{WGM}
\tsc{QE}
\usepackage{subcaption}
\usepackage{multirow}%
\usepackage{amsmath,amssymb,amsfonts,mathtools}%
\usepackage{amsthm}%
\usepackage{mathrsfs}%
\usepackage[title]{appendix}%
\usepackage{xcolor}%
\usepackage{textcomp}%
\usepackage{transparent}

\usepackage{docincludes} %exudyn definitions
\usepackage{tabularx}

\usepackage{acronym}
\usepackage{cleveref}
\usepackage{import}

\usepackage{longtable}

\usepackage[mathcal]{euscript} % produces an upright \mathcal letters

\theoremstyle{plain}

\theoremstyle{definition}

\theoremstyle{remark}

\usepackage[ruled,vlined]{algorithm2e}

\usepackage{tikz}
\usetikzlibrary{arrows.meta, positioning, fit, backgrounds}

\usepackage{placeins}

\usepackage{bookmark}

\usepackage{graphicx}
\usepackage{subcaption}

\usepackage[utf8]{inputenc} % if not already present
\DeclareUnicodeCharacter{2212}{$-$}
\DeclareUnicodeCharacter{2009}{\,}

\def\ba#1\ea{\begin{align}#1\end{align}}
\newcommand{\bbm}{\begin{bmatrix}}
\newcommand{\ebm}{\end{bmatrix}}
\newcommand{\ben}{\begin{enumerate}}
\newcommand{\een}{\end{enumerate}}

\renewcommand{\eqs}[2]{Eqs.~(\ref{#1}-\ref{#2})}

\newcommand{\tab}[1]{Tab.~\ref{#1}}

\newcommand{\sect}[1]{Sect.~\ref{#1}}
\newcommand{\sects}[2]{Sects.~\ref{#1}-\ref{#2}}

\newcommand{\appx}[1]{Appx.~\ref{#1}}
\renewcommand{\tp}{^\top}

\newcommand{\vecOp}{\textrm{vec}}
\renewcommand{\vec}[1]{\vecOp\!\left({#1}\right)}

\newcommand{\vecToSkew}[1]{\widetilde{#1}}

\newcommand{\RealNumbers}{\mathbb{R}}

\newcommand{\NaturalNumbers}{\mathbb{N}}
\newcommand{\ComplexNumbers}{\mathbb{C}}

\newcommand{\identity}{\mathbf{I}} 
\newcommand{\zero}{\mathbf{0}} 
\newcommand{\bigO}{\mathcal{O}} 

\newcommand{\SOThree}{SO(3)}
\newcommand{\SThree}{S^3}
\newcommand{\SEThree}{SE(3)}

\newcommand{\LieAlg}{\mathfrak{g}}

\newcommand{\tildeOperator}{\widetilde{(\bullet)}}

\renewcommand{\exp}{\textrm{exp}}

\newcommand{\expSOThree}{\textrm{exp}_{\SOThree}}

\newcommand{\expSThree}{\textrm{exp}_{\SThree}}

\newcommand{\Ad}{\operatorname{Ad}}

\newcommand{\elLieAlgSca}{w}
\newcommand{\elLieAlgVec}{\mathbf{\elLieAlgSca}}
\newcommand{\elLieAlg}{\widetilde{\elLieAlgVec}}

\newcommand{\rb}{\mathbf{r}}

\newcommand{\ab}{\mathbf{a}}

\newcommand{\Mb}{\mathbf{M}}

\newcommand{\gb}{\mathbf{g}}

\newcommand{\vb}{\mathbf{v}}
\newcommand{\vbSkew}{\widetilde{\mathbf{v}}}
\newcommand{\vbDot}{\dot{\vb}}
\newcommand{\lagMul}{\lambda}
\newcommand{\lagMulb}{\boldsymbol{\lagMul}}
\newcommand{\lagMulKin}{\eta}
\newcommand{\lagMulKinb}{\boldsymbol{\lagMulKin}}

\newcommand{\auxVar}{\mu}
\newcommand{\auxVarb}{\boldsymbol{\auxVar}}

\newcommand{\phib}{\mathbf{\Phi}}
\newcommand{\phibDot}{\dot{\phib}}

\newcommand{\RotMatScalar}{D}

\newcommand{\RotMat}{\mathbf{\RotMatScalar}}

\newcommand{\ine}{_{\textrm{ine}}}
\renewcommand{\int}{_{\textrm{int}}}
\newcommand{\ext}{_{\textrm{ext}}}
\newcommand{\jnt}{_{\textrm{jnt}}}
\newcommand{\csp}{_{\textrm{csp}}}

\newcommand{\aux}{_{\textrm{aux}}}
\newcommand{\fp}{_{\textrm{fp}}}
\newcommand{\eqp}{_{\textrm{e}}}

\newcommand{\ordMax}{n_{\textrm{max}}}

\newcommand{\CRes}{\boldsymbol{\mathcal{F}}}
\newcommand{\CResS}{\overline{\CRes}}
\newcommand{\MRes}{\rb}

\newcommand{\qVec}{\mathbb{q}}
\newcommand{\kinRecMat}{\mathbf{Q}}
\newcommand{\kinRecMatInv}{\mathbf{P}}
\newcommand{\qVecExuDot}{\boldsymbol{v}}
\newcommand{\qVecExuDDot}{\dot{\qVecExuDot}}
\newcommand{\SSA}{\boldsymbol{\mathcal{A}}}
\newcommand{\SSB}{\boldsymbol{\mathcal{B}}}

\newcommand{\SpecSubSpace}{\mathcal{E}}

\newcommand{\dimension}{n}
\newcommand{\dimG}{\dimension_G}
\newcommand{\dimJnt}{\dimension\jnt}
\newcommand{\dimAux}{\dimension\aux}
\newcommand{\dimCsp}{\dimension\csp}
\newcommand{\dimq}{\dimension_\qVec}
\newcommand{\dimSS}{\dimension_X}
\newcommand{\dimIVM}{\dimension_r}
\newcommand{\nB}{\dimension_{\textrm{b}}}
\newcommand{\ivm}{\mathcal W}
\newcommand{\dimRed}{\bar{n}_{X}}
\newcommand{\nJ}{\dimension_{\textrm{j}}} % number of joints

\newcommand{\cSpaceBody}{\mathcal{G}}

\newcommand{\EpSca}{e} 
\newcommand{\EpVec}{\mathbf{\EpSca}}

\newcommand{\pV}{\ttheta}
\newcommand{\pU}{\ttheta_U}
\newcommand{\pO}{\ttheta_\Omega}

\newcommand{\Pbridge}{\Xm}

\newcommand{\WmRed}{\Wm^{\mathrm{r}}}

\begin{document}
\let\WriteBookmarks\relax
\def\floatpagepagefraction{1}
\def\textpagefraction{.001}

% Short title
\shorttitle{}    

% Short author
\shortauthors{}  

% Main title of the paper 
\title [mode = title]{Residual-Based Spectral Submanifold Reduction of Constrained Multibody Systems on Matrix Lie Groups}

% Title footnote mark
% eg: \tnotemark[1]
%\tnotemark[1] 

% Title footnote 1.
% eg: \tnotetext[1]{Title footnote text}
%\tnotetext[1]{} 

% First author
%
% Options: Use if required
% eg: \author[1,3]{Author Name}[type=editor,
%       style=chinese,
%       auid=000,
%       bioid=1,
%       prefix=Sir,
%       orcid=0000-0000-0000-0000,
%       facebook=<facebook id>,
%       twitter=<twitter id>,
%       linkedin=<linkedin id>,
%       gplus=<gplus id>]

\author[1]{Stefan Holzinger}[
orcid=0000-0002-8678-2992,
]

% Corresponding author indication
\cormark[1]

% Footnote of the first author
%\fnmark[]

% Email id of the first author
\ead{stefan.holzinger@jku.ac.at}

% URL of the first author
%\ead[url]{}

% Credit authorship
% eg: \credit{Conceptualization of this study, Methodology, Software}
\credit{contributed to conceptualization, methodology, software development, formal analysis, investigation, visualization, and writing-review and editing}

% Address/affiliation
\affiliation[1]{organization={
Institute of Machine Design and Fluid Power,
Johannes Kepler University},
            addressline={Altenberger Straße 69}, 
            city={Linz},
%          citysep={}, % Uncomment if no comma needed between city and postcode
            postcode={4041}, 
            state={Upper Austria},
            country={Austria}}

\author[2]{Tiago Martins}[
orcid=0000-0002-3200-2225,
]

%% Footnote of the second author
%\fnmark[]
%
%% Email id of the second author
%\ead{}
%
%% URL of the second author
%\ead[url]{}
%
% Credit authorship
\credit{contributed to conceptualization, methodology, formal analysis, and writing-review and editing}
%
%% Email id of the second author
%\ead{}

% Address/affiliation
\affiliation[2]{organization={Chair of Applied Mechanics, TUM School of Engineering and Design, Technical University of Munich},
            addressline={Boltzmannstr. 15, Garching b. M\"unchen}, 
            city={Munich},
%          citysep={}, % Uncomment if no comma needed between city and postcode
            postcode={85748}, 
            state={Bavaria},
            country={Germany}}

\author[2]{Daniel Rixen}[
orcid=0000-0002-2303-4292,
]

% Credit authorship
\credit{contributed to methodology and writing, review, and editing}

% Corresponding author text
\cortext[1]{Corresponding author}

% Footnote text
%\fntext[1]{}

% For a title note without a number/mark
%\nonumnote{}

% Here goes the abstract
\begin{abstract}
The parameterization method yields efficient and accurate reduced-order models for nonlinear dynamical systems, but existing procedures that compute them assume polynomial vector fields on finite-dimensional vector spaces. Constrained multibody systems instead obey differential-algebraic equations on matrix Lie groups, where rigid-body configurations evolve on manifolds such as $\RealNumbers^3 \times SO(3)$ and $SE(3)$.
Minimally parameterizing such manifolds yields a vector-space representation but creates highly nonlinear, non-polynomial equations that one must approximate or recast with auxiliary variables before reduced order modeling. 
Existing implementations require nonlinearities as polynomial multilinear maps, which become cumbersome to construct explicitly for general assembled multibody systems.
 
This work develops a parameterization framework for constrained multibody systems on matrix Lie groups.
A polynomial computational representation reformulates the Lie group kinematics without altering the underlying geometrically exact formulation. 
For rigid and floating frame of reference bodies, polynomial kinematics are obtained through rotation-matrix components or Euler parameters, while geometrically exact beam elements use an Euler--Rodrigues representation that eliminates the logarithmic map and expresses the exponential map polynomially, thereby retaining screw-motion interpolation.
The reduced-order models satisfy linear joint constraints exactly, and a rotation-preserving parameterization keeps the configuration on the Lie group to machine precision.
%The parameterization method is formulated directly on the residual of the assembled equations of motion, to avoid decomposing the nonlinear forces into multilinear maps.
The proposed approach formulates the parameterization method directly in terms of the assembled residual, which collects the equations of motion, kinematic equations, auxiliary dynamics, and constraint equations. 
It requires only residual evaluations and their Fréchet derivatives, avoiding the explicit construction of multilinear force maps.
Truncated Taylor arithmetic delivers all Taylor coefficients of a given order from residual evaluation, and a non-intrusive coupling to codes with real-valued residuals is demonstrated with the multibody code Exudyn.
Numerical examples demonstrate accurate prediction of backbone curves, forced response curves, and internal resonances at substantially reduced computational cost.
\end{abstract}

% Use if graphical abstract is present
%\begin{graphicalabstract}
%\includegraphics{}
%\end{graphicalabstract}

%% Research highlights
%\begin{highlights}
%\item 
%\item 
%\item 
%\end{highlights}

% Keywords
% Each keyword is seperated by \sep

\maketitle

%%%%%%%%%%%%%%%%%%%%%%%%%%%%%%%%%%%%%%%%%%%%%%%%%%%%%%%%%%%%%%%%%%%%%%%%%%%%%%%%
%%%%%%%%%%%%%%%%%%%%%%%%%%%%%%%%%%%%%%%%%%%%%%%%%%%%%%%%%%%%%%%%%%%%%%%%%%%%%%%%
\section{Introduction}\label{sect_Introduction}
Accurate prediction of nonlinear vibrational responses, including backbone curves, isolated resonance branches, and internal resonances, has become increasingly important for engineering systems ranging from robots and flexible mechanisms to rotating machinery and lightweight structures \cite{touze2021model,han2025model,martin2023reducedorder}.

Many such systems are naturally formulated as constrained multibody systems on matrix Lie groups, where finite rigid-body motions, elastic deformation, and kinematic constraints coexist within a differential-algebraic equation (DAE) formulation \cite{Bruls2012,Sonneville2014,MuellerTerze2016}. 
These formulations provide geometrically exact kinematics, i.e., finite displacements and rotations are described without kinematic approximation, avoid the singularities of global rotation parameterizations, and admit efficient Lie group time integration schemes \cite{Bruls2012,Holzinger2025}. 
Despite these advances, the resulting models remain high-dimensional, making their repeated simulation computationally prohibitive for vibration analysis, optimization, and real-time applications \cite{han2025model,Pozzi2025}.

This motivates reduced-order modeling (ROM). 
Classical approaches rely on component mode synthesis within the floating frame of reference or motion formalism \cite{Sonneville2021}, and have been extended to geometric nonlinearity by modal derivatives \cite{WuTiso2016}; both approaches project onto linear subspaces. 
Among nonlinear reduction approaches, the parameterization method constructs reduced-order models on invariant manifolds of nonlinear dynamical systems \cite{haller2016nonlinear,touze2021model}.
This method \cite{haro2016parameterization} builds polynomial expansions whose coefficients follow from a recursive sequence of linear cohomological equations. %\cite{jainhaller2021how}. 
%Periodic responses are recovered as fixed points of the reduced dynamics, allowing entire forced response curves to be computed at negligible cost \cite{ponsioen2020model,li2022nonlinear}. 
Reduced dynamics also predicts backbone curves, isolas, internal resonances, and parametric resonances \cite{ponsioen2020model,li2022nonlinear,thurnherr2023nonautonomous}.
Formulated directly in physical coordinates, the method has been applied directly to large nonlinear finite element models \cite{vizzaccaro2021direct,jainhaller2021how,vizzaccaro2022high}, rotating structures \cite{martin2023reducedorder,gao2026nonintrusive}, and design optimization \cite{Pozzi2025}.
Since approximate kinematics can misrepresent the dynamics at large amplitudes \cite{Bagheri2024}, such reduced models must retain the geometrically exact description of the full-order model.
In the mathematics community, Cabré et al.~\cite{cabre2003parameterization} proved that there exists a unique smoothest nonlinear continuation of the spectral subspace of the linearized dynamics, in a small neighborhood of the equilibrium if strict non-resonance conditions hold. Later, Haller et al. named it the primary spectral submanifold (SSM)~\cite{haller2016nonlinear, haller2023nonlinear}, when introducing the concept to the engineering community.

Regardless of whether the primary SSM exists (at the amplitudes of interest), engineers have repeatedly and successfully applied the parameterization method in a variety of problems. 
To facilitate communication, henceforth we use the term SSM to refer to both the parameterization method and the manifolds it constructs.

Most existing SSM formulations, however, assume ordinary differential equations whose states evolve in finite-dimensional vector spaces. 
In contrast, the configuration variables of multibody systems on matrix Lie groups evolve on manifolds such as $\RealNumbers^3 \times SO(3)$ and $SE(3)$, while kinematic joint constraints restrict the motion to lower-dimensional submanifolds \cite{Bruls2012,Arnold2016,MuellerTerze2016}.
Whenever these configurations are expressed with a minimal set of coordinates, non-polynomial terms enter the formulation. 
For instance, the exponential map, together with its tangent operator and its inverse, arises in the configuration update of Lie-group integrators \cite{Bruls2012,Holzinger2023}. 
Higher-order derivatives of these operators may additionally appear in the equations of motion \cite{Cardona1989,Todesco2023}, while the logarithmic map appears in relative configuration and interpolation formulations for geometrically exact beams \cite{Sonneville2014}.
These operators introduce non-polynomial terms into the governing equations, whether through transcendental functions in exponential-coordinate formulations or rational terms in Cayley-map-based motion parameterizations.
The parameterization method, in turn, requires analytic vector fields on finite-dimensional vector spaces, and existing SSM implementations additionally require the nonlinearities in polynomial form.

Recent work has extended SSMs to constrained mechanical systems. 
Li and co-workers \cite{li2023model} extended SSM theory and computation to mechanical systems with holonomic constraints formulated as differential-algebraic equations. 
Their formulation assumes constant mass, damping, and stiffness matrices, and non-polynomial terms, such as trigonometric functions of rotation angles, are recast into polynomial form through auxiliary variables and additional algebraic constraints.
Han et al. \cite{han2025model} subsequently applied this methodology to flexible multibody systems with large deformations modeled using the absolute nodal coordinate formulation.
For free-floating geometrically exact beams and rigidly connected beam assemblies, including attached rigid bodies, Qi et al. \cite{qi2026invariantmanifoldbased} combined a floating frame decomposition with invariant manifold reduction, approximating the non-polynomial inertial and elastic terms by truncated Taylor expansions. 
Joints permitting relative motion, such as revolute or spherical joints, were not treated, although the authors consider the approach readily extensible to them.
An SSM framework for multibody systems on matrix Lie groups in which rigid bodies \cite{Bruls2011}, nodal-based floating frame of reference formulation (FFRF) bodies \cite{Zwoelfer2019}, and geometrically exact beam finite elements \cite{Sonneville2014} are connected by joints permitting relative motion therefore remains unavailable.
Existing implementations of the parameterization method require the nonlinear terms in polynomial form, decomposed into multilinear maps, which are supplied either as pre-assembled coefficient tensors \cite{jainhaller2021how,li2023model,han2025model} or as multilinear functions evaluated at the element or quadrature-point level \cite{vizzaccaro2022high, Martins2026MORFE}.
Deriving such decompositions for every element type, joint, and force law becomes increasingly cumbersome for assembled finite element and multibody models. 
Non-intrusive procedures avoid this derivation. Data-free variants evaluate the nonlinear forces of a finite element code as a black box, but assume nonlinearities of at most cubic order \cite{li2025datafree,gao2026nonintrusive}, whereas data-driven and data-assisted variants learn the nonlinear part of the SSM from simulated or measured trajectories \cite{Cenedese2022, Cenedese2025, 2025_Morsy}.

Lie groups have also been exploited as a tool for model reduction itself. 
Wotte et al. \cite{Wotte2026} approximate high-dimensional dynamical systems on manifolds by low-dimensional dynamical systems on Lie groups acting on the state space, identified by intrusive or data-driven optimization, to overcome the Kolmogorov $N$-width barrier of linear-subspace methods. 
In the present work, by contrast, the Lie group is the configuration space of the full-order multibody model, and the reduced-order model is the dynamics on an invariant manifold of that model, computed directly from its equations of motion without training data.

A central observation of this work is that assembled multibody residuals provide a natural interface between geometrically exact multibody formulations and the parameterization method. 
Each group element is represented redundantly by its entries, and group membership is enforced by configuration-space constraints with Lagrange multipliers, in the spirit of the constrained formulations of Betsch and Steinmann \cite{Betsch2001,Betsch2003}. 
Such redundant representations are established in time integration \cite{Betsch2001,Betsch2003} and in frequency-domain continuation of geometrically exact beams \cite{debeurre2024quaternionbased}. 
In invariant-manifold-based reduction, they have so far been restricted to planar rotations, represented by the cosine and sine of the rotation angle subject to a unit-circle constraint \cite{li2023model,Grolet2025}, where the latter additionally recast the remaining nonlinearities into quadratic form through auxiliary variables \cite{Grolet2025}.
For spatial rotations, the floating-frame motion has instead been kept outside the reduced model \cite{qi2026invariantmanifoldbased}. 
To the best of our knowledge, redundant representations of three-dimensional rotations have not yet been used for SSM computation.
With the proposed representation, the equations of motion of rigid and FFRF bodies become polynomial \cite{Zwoelfer2019,Holzinger2025FFRF}.
Polynomial formulations of geometrically exact beams exist, based on directors \cite{Betsch2003} or quaternions \cite{debeurre2024quaternionbased}. 
There, the rotational field is interpolated additively, so that orthonormality or the unit norm holds only at the nodes \cite{Betsch2003} or in a weak sense along the element \cite{debeurre2024quaternionbased}. 
We instead derive a log-free polynomial counterpart of the $SE(3)$ beam element of Sonneville et al. \cite{Sonneville2014}, which retains its screw-motion interpolation, so that every interpolated configuration lies exactly on $SE(3)$, using a few lifted auxiliary variables and minimax polynomials accurate to machine precision.

Many multibody and finite element codes evaluate a residual vector of the assembled equations of motion during implicit time integration and static analysis. 
The parameterization method is reformulated directly in terms of this residual, avoiding prior conversion into an explicit first-order or descriptor formulation. 
The required Taylor coefficients are computed by truncated Taylor arithmetic, i.e., higher-order forward automatic differentiation \cite{Berz1989,Griewank2008}. 
Taylor arithmetic is well established in physics \cite{Berz1989}, in the parameterization method for explicitly given maps and vector fields \cite{haro2016parameterization}, and in Taylor series-based continuation of nonlinear residual equations \cite{Koutsawa2008}.
To the best of our knowledge, however, it has not yet been used for SSM-based reduction of mechanical systems.

\paragraph{Contributions.}
Building on the above observations, this work develops an SSM framework for
constrained multibody systems on matrix Lie groups, covering assembled models
comprising rigid bodies, nodal-based FFRF bodies, and geometrically exact beam
finite elements.
A polynomial computational representation reformulates the Lie-group kinematics
without altering the underlying geometrically exact formulation. For rigid and
FFRF bodies, polynomial kinematics are obtained through rotation-matrix
components or Euler parameters. For geometrically exact beam elements, an
Euler--Rodrigues representation eliminates the logarithmic map while retaining
screw-motion interpolation of the $SE(3)$ element.
A group-structure-preserving reduction is developed together with an analysis of constraint preservation in the reduced-order model. 
An ambient polynomial parameterization, in which the configuration variables are expanded directly in the vector space embedding the Lie group, preserves linear joint constraints exactly, while nonlinear constraints are satisfied up to the truncation order.
For the rotation-preserving variant, the group block is
parameterized in the tangent space and the resulting configuration is
subsequently corrected along the kinematic tree, yielding exact Lie-group
membership and exact joint-constraint satisfaction.
A residual-based reformulation of the parameterization method computes SSMs for
autonomous and periodically forced systems directly from the residual vector of
the assembled equations of motion, avoiding the explicit construction of
multilinear force tensors.
The methodology is realized through truncated Taylor arithmetic, which computes
all Taylor coefficients of a given order from a single residual evaluation. For
codes that provide real-valued residual evaluations only, a non-intrusive
coupling reconstructs the residual as an explicit polynomial, as demonstrated
with the open-source software Exudyn \cite{Gerstmayr2024}.

\medskip

The remainder of this paper is organized as follows. 
\sect{sect_LieGroupMBD} develops the polynomial formulation of constrained multibody systems on matrix Lie groups and introduces the computational representation required for SSM computation.
\sect{sect_ResidualRepresentation} presents the residual-based reformulation of the parameterization method, and \sect{sect_ConstraintPreservation} analyzes constraint preservation in the reduced-order model. 
\sect{sec:implementation-aspects} describes the implementation based on truncated Taylor arithmetic and non-intrusive coupling to existing multibody codes. 
\sect{sec:numerical-examples} demonstrates the proposed methodology on representative constrained multibody systems, and \sect{sec_Conclusions} concludes the paper.

%%%%%%%%%%%%%%%%%%%%%%%%%%%%%%%%%%%%%%%%%%%%%%%%%%%%%%%%%%%%%%%%%%%%%%%%%%%%%%%%
\section{Multibody dynamics on matrix Lie groups in polynomial form}\label{sect_LieGroupMBD}
This section addresses the formulation of the equations of motion of constrained multibody systems in polynomial form. The section is organized as follows.

\sect{sect_LieGroupEOM} summarizes the equations of motion of periodically forced constrained multibody systems in a Lie group setting.
\sect{sect_CompRep} introduces a computational representation of the configuration variables describing body displacements and rotations, compatible with the parameterization method.
Building on this representation, \sect{sect_RigidBodyPoly} derives polynomial kinematic equations for rigid-body motion, thereby enabling polynomial equations of motion for constrained multibody systems comprising rigid and floating frame of reference formulation bodies.
Subsequently, \sect{sect_PolyInterpolation} derives a polynomial interpolation of rigid-body motion, providing the basis for a polynomial formulation of geometrically exact beam finite elements.

Within the Lie group framework, the configuration variables are represented by matrix Lie-group elements, such as homogeneous transformation matrices in the special Euclidean group $SE(3)$ that describe rigid-body poses. Once these variables are replaced by a truncated polynomial parameterization in the reduced coordinates, they are no longer guaranteed to remain on the underlying Lie group. The final part of this section addresses this issue and develops the foundation for the group-structure-preserving parameterization introduced later in \sect{sec_LieParameterization}.

\subsection{Equations of motion on Lie groups}\label{sect_LieGroupEOM}
The dynamics of multibody systems consisting of $\nB$ bodies can be described on a product manifold with Lie group structure $G=\cSpaceBody_1\times\cSpaceBody_2\times\ldots\times\cSpaceBody_{\nB}$, where each factor $\cSpaceBody_i$ denotes the configuration space of a single rigid or flexible body, cf. \cite{Bruls2012,Sonneville2014,Arnold2016}.
Depending on the underlying multibody or finite-element formulation, the body configuration space $\cSpaceBody_i$ may be given by $\RealNumbers^{3}\times SO(3)$, $SE(3)=\RealNumbers^{3}\rtimes SO(3)$, or other combinations of translational and rotational variables constructed by direct ($\times$) or semi-direct ($\rtimes$) products, including Euler-parameter\footnote{Euler parameters denote the four components of a unit quaternion and form a Lie group isomorphic to the unit sphere in $\RealNumbers^4$, defined by $S^3:=\{\ttheta\in\RealNumbers^4\mid\|\ttheta\|=1\}$. Consequently, $\cSpaceBody_i$ may also be given by $\RealNumbers^{3}\times S^3$ or $\RealNumbers^{3}\rtimes S^3$.} or other rotation parameter-based formulations \cite{Holzinger2021}.
In the absolute coordinate approach adopted in this paper, an element $q\in G$ collects the absolute displacement and rotational variables of all bodies, which are initially treated as independent. 
The corresponding velocity is represented by a vector $\vb\in\RealNumbers^{\dimG}$, where $\dimG=\sum_{i=1}^{\nB}\dim(\cSpaceBody_i)$.
These variables are generally subject to $\dimJnt$ kinematic joint constraints
{\setlength{\abovedisplayskip}{5pt}
 \setlength{\belowdisplayskip}{5pt}
\be\label{eq_jntConstrMapping}
\phib\jnt:G\rightarrow\RealNumbers^{\dimJnt}\eqComma
\ee}
representing joints such as spherical, revolute, and prismatic joints
\cite{Shabana2010}. The constraints may be imposed either at the position
level through \eq{eq_jntConstrMapping} or at the velocity level,
{\setlength{\abovedisplayskip}{5pt}
 \setlength{\belowdisplayskip}{5pt}
\be\label{eq_jntConstrMapping_t}
\phibDot\jnt:G\times\RealNumbers^{\dimG}\rightarrow\RealNumbers^{\dimJnt}\eqComma
\ee}
thereby restricting the dynamics to a submanifold of dimension
$\dimG-\dimJnt$.

Accordingly, we consider periodically forced constrained multibody systems of the form
\begin{align}
\Mb(q)\vbDot + \gb\ine(q, \vb) + \gb\int(q, \vb) + \gb\aux(q,\vb,\auxVarb,t) + \Bm\jnt\tp(q)\lagMulb &=  \epsilon\,\gb_\textrm{ext}(q, \vb, \Omega t) 
\eqComma \label{eq_DynEOM} \\
\dot{q} &= DL_{q}(e)\cdot \vbSkew \label{eq_ODELieGroupe}\eqComma \\
\dot\auxVarb &= \ab(q,\vb,\auxVarb,t) \eqComma \label{eq_auxVarODE} \\
\phib\jnt(q) &= \zero \eqComma \label{eq_phi}\\
0 \leq\;\; & \epsilon \ll 1 \eqComma \label{eq_epsilon}
\end{align}
where, $\epsilon$ denotes the amplitude of the external excitation and
$\Omega$ its forcing frequency. 
The autonomous unforced system is recovered for $\epsilon=0$.

\eq{eq_DynEOM} represents the balance of forces and momenta. The mass matrix $\Mb(q)\in\RealNumbers^{\dimG\times\dimG}$ is assumed to be symmetric and positive definite, while $\gb\ine(q,\vb),\gb\int(q,\vb)\in\RealNumbers^{\dimG}$ denote the complementary inertia and internal force contributions, respectively.

The term $\gb\aux(q,\vb,\auxVarb,t)$ captures the contribution of auxiliary variables $\auxVarb\in\RealNumbers^{\dimAux}$ to the system dynamics. These variables are not directly represented by the configuration variables $q$ but describe additional quantities of engineering interest.
Examples include relative distances, relative rotation angles, constitutive variables, and viscoelastic state variables. Their evolution is governed by \eq{eq_auxVarODE}. %, where $\ab(q,\vb,\auxVarb,t)\in\RealNumbers^{\dimAux}$ is assumed to be sufficiently smooth.

The term $\Bm\jnt\tp(q)\lagMulb$ represents the constraint forces associated with the kinematic joint constraints \eqs{eq_jntConstrMapping}{eq_jntConstrMapping_t}, where $\lagMulb\in\RealNumbers^{\dimJnt}$ is the vector of Lagrange multipliers.
The associated constraint gradient matrix is given by
\begin{equation}
\Bm\jnt(q)
=
\dfrac{\partial\phibDot\jnt(q,\vb)}{\partial\vb}
\in
\RealNumbers^{\dimJnt\times\dimG}.
\label{eq_Bq}
\end{equation}
The term $\gb_\textrm{ext}(q,\vb,\Omega t)\in\RealNumbers^{\dimG}$ denotes external forces, which are assumed to be harmonic in time.

\eq{eq_ODELieGroupe} denotes the kinematic reconstruction equation, which relates the configuration rate $\dot q$ to the Lie algebra element $\vbSkew$ through the differential $DL_q(e)$ of the left translation map $L_q$ evaluated at the neutral element $e \in G$ \cite{Bruls2012}.
The Lie algebra element $\vbSkew$ is obtained by mapping the velocity vector $\vb$ to its matrix Lie algebra representation through the invertible linear mapping
\be\label{eq_DefTildeOperator}
\tildeOperator:\RealNumbers^{\dimG}\rightarrow\LieAlg\eqComma\quad
\elLieAlgVec\mapsto\elLieAlg\eqComma
\ee
where $\LieAlg:=T_eG$ is the tangent space of $G$ at $e$ \cite{Arnold2016}.

The kinematic reconstruction \eq{eq_ODELieGroupe} evolves on the $\dimG$-dimensional matrix Lie group $G$. 
However, the parameterization method requires analytic mappings between finite-dimensional vector spaces.
Consequently, the formulation of the kinematic reconstruction in the form of \eq{eq_ODELieGroupe} prevents the direct computation of invariant manifolds from \eqs{eq_DynEOM}{eq_epsilon}.

As can be seen from \eqs{eq_DynEOM}{eq_epsilon}, all terms contributing to the equations of motion depend, in general, on the configuration variables
$q\in G$. 
Consequently, a polynomial formulation of the equations of motion requires the rotation matrix to be represented as a polynomial in the configuration variables, since it enters the kinematic joint constraints, the kinematic reconstruction equations, and potentially the computation of forces and momenta.
As shown later in \sects{sect_RigidBodyPoly}{sect_PolyInterpolation}, two representations achieve this naturally. 

%%%%%%%%%%%%%%%%%%%%%%%%%%%%%%%%%%%%%%%%%%%%%%%%%%%%%%%%%%%%%%%%%%%%%%%%%%%%%%%%
\subsection{Computational representation}
\label{sect_CompRep}
%%%%%%%%%%%%%%%%%%%%%%%%%%%%%%%%%%%%%%%%%%%%%%%%%%%%%%%%%%%%%%%%%%%%%%%%%%%%%%%%
%
To obtain a vector-space representation, we introduce the vectorization
operator $\vecOp(\bullet)$, which stacks the entries of an $m\times n$ matrix
into an $mn$-dimensional column vector,
\be\label{eq_DefVecOperator}
\vecOp:\RealNumbers^{m\times n}\rightarrow\RealNumbers^{mn}\eqComma
\qquad
\Um\mapsto\vec{\Um}\eqDot
\ee
Applied to $q\in G$, it gives
\be\label{eq_qVec}
\qVec=\vec q\in\RealNumbers^{\dimq}\eqComma
\ee
where $\dimq=\#G$ is the number of entries of $q$.

Expressed in terms of the computational variable $\qVec$, the kinematic reconstruction \eq{eq_ODELieGroupe} becomes
\be\label{eq_ODELieGroupVectorized}
\dot\qVec=\kinRecMat(\qVec)\vb
\quad\xleftrightarrow{\vecOp}\quad
\vec{\dot q}=\vec{DL_q(e)\cdot\vbSkew}\eqComma
\ee
which is the equivalent matrix-vector form of \eq{eq_ODELieGroupe} in $\RealNumbers^{\dimq}$ and is therefore compatible with the parameterization method.
The matrix $\kinRecMat(\qVec)$ is given by
\be\label{eq_CompRecMat}
\kinRecMat(\qVec)
=
(e\otimes q)\Em
\in\RealNumbers^{\dimq\times\dimG}\eqComma
\ee
where
\be\label{eq_LieAlgebraBasisMatrix}
\setlength{\arraycolsep}{2pt}
\Em=
\begin{bmatrix}
\vec{\widetilde{\ev}_1}&
\vec{\widetilde{\ev}_2}&
\cdots&
\vec{\widetilde{\ev}_{\dimG}}
\end{bmatrix}
\in\RealNumbers^{\dimq\times\dimG}\eqDot
\ee
Here, $\otimes$ denotes the Kronecker product, and
$\widetilde{\ev}_i\in\LieAlg$ is the skew-symmetric matrix associated with
the $i$-th standard basis vector $\ev_i\in\RealNumbers^{\dimG}$.
The basis matrices $\widetilde{\ev}_i$ are ordered in \eq{eq_LieAlgebraBasisMatrix} such that $\widetilde{\vb}=\sum_{i=1}^{\dimG}v_i\,\widetilde{\ev}_i$, where $v_i \in \RealNumbers$ denotes the $i$-th component of $\vb$.

The next section shows that the computational representation \eqs{eq_qVec}{eq_ODELieGroupVectorized} yields a polynomial formulation of the rigid-body kinematics.
Throughout the remainder of this work, all $q$-dependent equations are formulated directly in terms of the computational variable $\qVec$, with the inverse vectorization implied. 
Thus, expressions such as $\Mm(\qVec)$, $\Bm\jnt(\qVec)$, $\phib(\qVec)$, and $\gb\int(\qVec)$ denote the corresponding mappings evaluated at $q=\vecOp^{-1}(\qVec)$.

%%%%%%%%%%%%%%%%%%%%%%%%%%%%%%%%%%%%%%%%%%%%%%%%%%%%%%%%%%%%%%%%%%%%%%%%%%%%%%%%
\subsection{Polynomial formulation of rigid-body kinematics}\label{sect_RigidBodyPoly}
%%%%%%%%%%%%%%%%%%%%%%%%%%%%%%%%%%%%%%%%%%%%%%%%%%%%%%%%%%%%%%%%%%%%%%%%%%%%%%%%
%
This subsection shows that the computational representation introduced in \sect{sect_CompRep} yields a polynomial formulation of rigid-body kinematics for two computational representations of $\SOThree$: rotation-matrix components and Euler parameters.
Both representations yield polynomial equations and may be used interchangeably. 
They differ primarily in the polynomial degree of the resulting rotation matrix and in the dimension of the group block of the computational variables $\qVec$.
These differences and their implications for reduced-order modeling are discussed later in \sects{sec_LinearConstraints}{sec_LieParameterization}.
The extension to $\RealNumbers^3\times\SOThree$ and $\SEThree$ introduces three translational components whose governing equations are already polynomial and therefore require no further treatment.

%%%%%%%%%%%%%%%%%%%%%%%%%%%%%%%%%%%%%%%%%%%%%%%%%%%%%%%%%%%%%%%
\subsubsection{Rotation-matrix components}
\label{sect_RotationMatrixComponents}
%%%%%%%%%%%%%%%%%%%%%%%%%%%%%%%%%%%%%%%%%%%%%%%%%%%%%%%%%%%%%%
%
Let $\RotMat\in\SOThree$ denote the rotation matrix of a rigid body.
Vectorization yields the computational variable
{\setlength{\abovedisplayskip}{5pt}
 \setlength{\belowdisplayskip}{5pt}
\be
\qVec_\RotMatScalar = \vec{\RotMat}\in\RealNumbers^9 \eqComma
\ee}
while the velocity vector reduces to the angular velocity
$\vb=\tomega\in\RealNumbers^3$.
The vectorized kinematic reconstruction equation then takes the form
{\setlength{\abovedisplayskip}{5pt}
 \setlength{\belowdisplayskip}{5pt}
\be
\dot{\qVec}_\RotMatScalar=\kinRecMat_\RotMatScalar(\qVec_\RotMatScalar)\tomega \eqComma
\qquad
\textrm{with} \qquad
\kinRecMat_\RotMatScalar(\qVec_\RotMatScalar)
=(\Im_3\otimes\RotMat)\Em
\in\RealNumbers^{9\times3} \eqDot
\label{eq_app_Q_SO3}
\ee}
Unlike the original kinematic reconstruction ODE on $SO(3)$ $\dot{\RotMat}=\RotMat\widetilde{\tomega}$ (cf. \eq{eq_ODELieGroupe}), the representation \eq{eq_app_Q_SO3} evolves the vectorized rotation matrix directly in $\RealNumbers^9$, see e.g., \cite{Betsch2001}. 

Two properties of this representation are worth mentioning explicitly. 
The rotation matrix is recovered from $\qVec_\RotMatScalar$ by reshaping $\RotMat=\vecOp^{-1}(\qVec_\RotMatScalar)$ and is therefore
\emph{linear} --- indeed the identity --- in the configuration variables.
And every entry of $\kinRecMat_\RotMatScalar(\qVec_\RotMatScalar)=(\Im_3\otimes\RotMat)\Em$ is a single component of $\qVec_\RotMatScalar$, so $\kinRecMat_\RotMatScalar(\qVec_\RotMatScalar)\vb$ is bilinear in $(\qVec_\RotMatScalar,\vb)$.

%%%%%%%%%%%%%%%%%%%%%%%%%%%%%%%%%%%%%%%%%%%%%%%%%%%%%%%%%%%%%%
\subsubsection{Euler parameters}
\label{sect_EulerParameters}
%%%%%%%%%%%%%%%%%%%%%%%%%%%%%%%%%%%%%%%%%%%%%%%%%%%%%%%%%%%%%%
%
Using the Euler parameters
{\setlength{\abovedisplayskip}{5pt}
 \setlength{\belowdisplayskip}{5pt}
\be
\qVec_\EpSca =(\EpSca_0,\EpVec)\in\SThree
=\{\qVec_\EpSca\in\RealNumbers^4\;\vert\;\Vert\qVec_\EpSca\Vert=1\}\eqComma
\label{eq_EpDef}
\ee}
as computational coordinates, the resulting rotation matrix $\RotMat(\qVec_\EpSca)\in SO(3)$ becomes a quadratic polynomial in $\qVec_\EpSca$,
{\setlength{\abovedisplayskip}{5pt}
 \setlength{\belowdisplayskip}{5pt}
\be
\RotMat(\qVec_\EpSca)
=(\EpSca_0^2-\EpVec\tp\EpVec)\identity
+2\,\EpVec\EpVec\tp
+2\,\EpSca_0\widetilde{\EpVec}\eqComma
\label{eq_RfromER}
\ee}
while the velocity vector reduces to the angular velocity $\vb=\tomega\in\RealNumbers^3$.

Using \eq{eq_ODELieGroupVectorized}, the vectorized kinematic reconstruction equations on $\SThree$ become
\be
\dot{\qVec}_\EpSca=\kinRecMat_\EpSca(\qVec_\EpSca)\,\tomega  \eqComma
\qquad 
\textrm{with} \qquad
\kinRecMat_\EpSca(\qVec_\EpSca)
=\tfrac12
\begin{bmatrix}
-\EpVec\tp\\[2pt]
\EpSca_0\identity+\widetilde{\EpVec}
\end{bmatrix} \in\RealNumbers^{4\times3} \eqComma
\label{eq_QmatER}
\ee
which are bilinear in the Euler parameters and angular velocity.

%%%%%%%%%%%%%%%%%%%%%%%%%%%%%%%%%%%%%%%%%%%%%%%%%%%%%%%%%%%%%%%%%%%%%%%%%%%%%%%%
\subsection{Floating frame of reference bodies}
\label{sect_FFRF}
%%%%%%%%%%%%%%%%%%%%%%%%%%%%%%%%%%%%%%%%%%%%%%%%%%%%%%%%%%%%%%%%%%%%%%%%%%%%%%%%
%
Building on \sect{sect_RigidBodyPoly}, this subsection shows that the floating frame of reference formulation preserves polynomiality whenever the underlying rigid-body motion is represented polynomially.
\begin{figure}
\centering
\includegraphics[width=\textwidth]{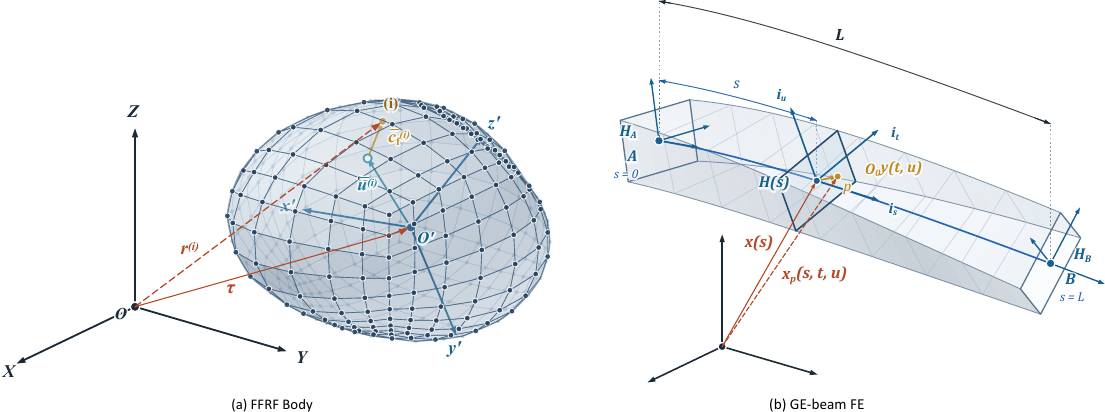}
\caption{The flexible bodies considered in this work. (a) A floating frame of reference body carries a pose together with a deformation field measured relative to that frame. 
(b) A geometrically exact beam finite element of length $L$ carries a pose at each node and interpolates the relative configuration between them. }
\label{fig_geomExactBodies}
\end{figure}

The position of finite-element node $(i)$ with respect to the inertial frame is
\be
\rv^{(i)}
=\ttau+\RotMat\left(\overline{\uv}^{(i)}+\overline{\cv}\indf^{(i)}\right)\eqComma
\label{eq_FFRFPosition}
\ee
where $\ttau\in\RealNumbers^3$ and $\RotMat\in SO(3)$ denote the position and orientation of the floating frame (FF), respectively, $\overline{\uv}^{(i)}$ is the undeformed nodal position in the FF, and $\overline{\cv}\indf^{(i)}$ is the corresponding flexible nodal displacement, cf.~\fig{fig_geomExactBodies}a.

Polynomiality follows directly from \eq{eq_FFRFPosition}. 
Whenever the underlying rigid-body kinematics is polynomial, the FFRF equations of motion inherit this polynomial structure because all coupling terms are assembled from the rigid-body transformation and its derivatives through polynomial matrix operations, see, e.g., \cite{Zwoelfer2019,Holzinger2025FFRF}.
Consequently, the mass matrix and complementary inertia forces are polynomial in the state variables when the rigid-body motion of the FF is represented by the polynomial formulations presented in \sect{sect_RigidBodyPoly}.
The internal elastic forces are typically assumed to be linear, i.e., $\gv\indf = \overline{\Km}\, \overline{\cv}$, where $\overline{\Km}$ is the constant finite-element stiffness matrix.
Of course, internal force expressions based on higher-order polynomial stiffness tensors preserve polynomiality while extending the formulation to large deformations.

%%%%%%%%%%%%%%%%%%%%%%%%%%%%%%%%%%%%%%%%%%%%%%%%%%%%%%%%%%%%%%%%%%%%%%%%%%%%%%%%
\subsection{Geometrically exact beam finite elements}
\label{sect_PolyInterpolation}
%%%%%%%%%%%%%%%%%%%%%%%%%%%%%%%%%%%%%%%%%%%%%%%%%%%%%%%%%%%%%%%%%%%%%%%%%%%%%%%%
%
Building on the rigid-body polynomial kinematics established in \sect{sect_RigidBodyPoly}, this subsection derives a polynomial interpolation of rigid-body motion for geometrically exact beam finite elements.

A geometrically exact beam \cite{Reissner1972,Reissner1981,Simo1985,SimoVuQuoc1986} models a slender body as a curve of rigid cross-sections undergoing arbitrarily large displacements and rotations, see \fig{fig_geomExactBodies}b. 
The term geometrically exact refers to the fact that the relations between the configuration and the strain measures are fully consistent with the virtual work principle and the equations of motion, regardless of the magnitude of displacements, rotations, and strains \cite{JelenicCrisfield1999}.
The only kinematic assumption is the rigid cross-section \cite{Simo1985}.
Formulated on the special Euclidean group $\SEThree$, the rigid-body motion between the nodal frames $\Hm_A,\Hm_B\in SE(3)$ is commonly interpolated by a screw motion parameterized by $s\in[0,L]$, where $A$ and $B$ denote the end nodes at $s=0$ and $s=L$, respectively,
\begin{equation}
\Hm(s)=\Hm_A\exp_{SE(3)}\!\bigl((s/L)\widetilde{\dv}\bigr) \in SE(3) \eqComma
\qquad
\widetilde{\dv}=\log_{SE(3)}\!\bigl(\Hm_A^{-1}\Hm_B\bigr) \in T_e SE(3) \eqComma
\label{eq:sonneville}
\end{equation}
where $\dv = [\dv_U\tp \quad \dv_\Omega\tp]\tp=\vecOp(\widetilde{\dv})\in\RealNumbers^6$ denotes the relative configuration vector with $\dv_U, \dv_\Omega \in\RealNumbers^3$; see, e.g. \cite{Sonneville2014}.

The exponential and logarithmic maps in \eq{eq:sonneville} are non-polynomial and therefore incompatible with the parameterization method. 
This section derives a polynomial representation of the geometrically exact beam kinematics using Euler--Rodrigues motion parameters. 
First, the nodal unknowns of the beam element are defined and a polynomial representation of the relative configuration vector $\dv$ is derived.
Second, this representation is used to construct a polynomial interpolation of the rigid-body motion along the beam element.
Together, these steps form the basis for the geometrically exact polynomial beam formulation summarized in Appx.~\ref{app_MinimaxCoefficients}-\ref{app_PolyBeam}.

%%%%%%%%%%%%%%%%%%%%%%%%%%%%%%%%%%%%%%%%%%%%%%%%%%%%%%%%%%%%%%%%%%%%%%%%%%%%%%%%
\subsubsection{Nodal unknowns}
\label{sect_BeamNodalUnknowns}
%%%%%%%%%%%%%%%%%%%%%%%%%%%%%%%%%%%%%%%%%%%%%%%%%%%%%%%%%%%%%%%%%%%%%%%%%%%%%%%%
%
Each beam node $i\in\{A,B\}$ is treated as a rigid cross-section in $\RealNumbers^3\rtimes S^3$, described by its global position $\rv_i$ and by the Euler parameters $\qVec_{\EpSca,i}$ of \sect{sect_EulerParameters}.
Collecting both yields the computational variables of \sect{sect_CompRep} for node $i$,
\be
\qVec_i=
\begin{bmatrix}
\rv_i\\
\qVec_{\EpSca,i}
\end{bmatrix}
\in\RealNumbers^{7}\eqComma
\qquad
\vv_i=
\begin{bmatrix}
\vv_{i,U}\\
\vv_{i,\Omega}
\end{bmatrix}
\in\RealNumbers^{6}\eqComma
\label{eq_nodalUnknowns}
\ee
where both velocity components are expressed in the body frame of the node.
The nodal kinematic reconstruction equations are given by $\dot{\qVec}_{\EpSca,i}=\kinRecMat_\EpSca(\qVec_{\EpSca,i})\vv_{i,\Omega}$ (cf. \eq{eq_QmatER}) for rotation and by $\dot{\rv}_i=\RotMat(\qVec_{\EpSca,i})\vv_{i,U}$ for translation.

%%%%%%%%%%%%%%%%%%%%%%%%%%%%%%%%%%%%%%%%%%%%%%%%%%%%%%%%%%%%%%%%%%%%%%%%%%%%%%%%
\subsubsection{Log-free relative configuration}
%%%%%%%%%%%%%%%%%%%%%%%%%%%%%%%%%%%%%%%%%%%%%%%%%%%%%%%%%%%%%%%%%%%%%%%%%%%%%%%%
%
The beam interpolation \eq{eq:sonneville} requires the configuration of node $B$ relative to node $A$, expressed through the relative configuration vector $\dv$.
In this section, the relative configuration is first expressed in terms of the Euler--Rodrigues motion vector
\be
\pV =
\begin{bmatrix}
\pU\\
\pO
\end{bmatrix}
\in\RealNumbers^6\eqComma
\label{eq_thetaDef}
\ee
with translational part $\pU\in\RealNumbers^3$ and rotational part $\pO\in\RealNumbers^3$.
Unlike $\dv$, the vector $\pV$ follows from the nodal unknowns without the logarithmic map, and $\dv$ is subsequently recovered from $\pV$ by an algebraic relation.

The rotational part $\pO$ consists of the Euler--Rodrigues parameters, which are obtained from the vectorial part of the Euler parameters $(\EpSca_0,\EpVec)$ by
\be
\pO=2\,\EpVec\eqDot
\label{eq_erpDef}
\ee
The relative rotation is given by the relative Euler parameters, i.e., the quaternion product of the conjugate of $\qVec_{\EpSca,A}$ with $\qVec_{\EpSca,B}$.
The relative translation $\tv$ is the position of node $B$ with respect to node $A$. %, expressed in the frame of node $A$.
With \eq{eq_erpDef}, the relative configuration reads
\be
\begin{aligned}
\pU&=\left(q_0\identity-\frac12\vecToSkew{\pO}\right)\tv\eqComma
&\qquad
\tv&=\RotMat\tp(\qVec_{\EpSca,A})\left(\rv_B-\rv_A\right)\\[3pt]
\pO&=2\left(\EpSca_{0,A}\,\EpVec_B-\EpSca_{0,B}\,\EpVec_A-\EpVec_A\times\EpVec_B\right)\eqComma
&\qquad
q_0&=\EpSca_{0,A}\,\EpSca_{0,B}+\EpVec_A\tp\EpVec_B \eqDot
\end{aligned}
\label{eq_relativeErpFromNodes}
\ee
Since the quaternion product preserves the norm, the pair $(q_0,\pO)$ satisfies 
\be
q_0^2+\frac14\,\pO\tp\pO=1\eqComma
\label{eq_erpUnit}
\ee
and the relative rotation matrix follows from \eq{eq_RfromER} as
\be
\RotMat(q_0,\pO)=\identity+q_0\,\vecToSkew{\pO}+\frac12\vecToSkew{\pO}^{\,2}\eqDot
\label{eq_RotErp}
\ee
Note, every entry of \eq{eq_relativeErpFromNodes} is polynomial in the nodal unknowns.

The rates of the relative configuration involve the nodal velocities only through the velocity of node $B$ relative to node $A$, expressed in the frame of node $B$ \cite{Sonneville2014},
\be
\Delta\vv=
\begin{bmatrix}
\Delta\vv_U\\
\Delta\vv_\Omega
\end{bmatrix}
=
\begin{bmatrix}
\vv_{B,U}-\RotMat\tp\left(\vv_{A,U}-\vecToSkew{\tv}\,\vv_{A,\Omega}\right)\\[4pt]
\vv_{B,\Omega}-\RotMat\tp\vv_{A,\Omega}
\end{bmatrix}\eqComma
\label{eq_relVelocity}
\ee
with $\RotMat=\RotMat(q_0,\pO)$ of \eq{eq_RotErp} and $\tv$ of \eq{eq_relativeErpFromNodes}.
In particular, the pair $(q_0,\pO)$ obeys the same vectorized kinematic reconstruction \eq{eq_QmatER} as a nodal pair, with the relative angular velocity $\Delta\vv_\Omega$, so that
\be
\dot q_0=-\frac14\,\pO\tp\Delta\vv_\Omega\eqDot
\label{eq_q0Dot}
\ee

Next, $\dv$ is recovered from $\pV$.
Let $\phi\in[0,\pi)$ and $\mathbf{n}$ be the angle and the unit axis of the relative rotation $\RotMat(q_0,\pO)$, and let $a=(\phi/2)^2$.
In terms of $\phi$, $\mathbf{n}$, and $\tv$, both $\dv$ and $\ttheta$ are given explicitly by
\be
\begin{aligned}
\pO &= 2\sin\sqrt a\;\mathbf{n}\eqComma
&\qquad
\pU &= \left(q_0\identity-\frac12\vecToSkew{\pO}\right)\tv
\eqComma\\[3pt]
\dv_\Omega &= 2\sqrt a\;\mathbf{n}\eqComma
&\qquad
\dv_U &= \left(\identity-\frac12\vecToSkew{\dv_\Omega}
+F(4a)\,\vecToSkew{\dv_\Omega}^{\,2}\right)\tv\eqComma
\end{aligned}
\label{eq_bothCharts}
\ee
with $q_0=\cos\sqrt a$.
The first line is \eq{eq_relativeErpFromNodes} written in terms of $\phi$ and $\mathbf{n}$.
The second line is the exponential coordinate $\dv$ of \eq{eq:sonneville}, in which the matrix multiplying $\tv$ is the inverse translational block of $\exp_{SE(3)}$ \cite{Sonneville2014}.
Both translational parts are thus linear maps of the same relative translation $\tv$.
The scalar functions appearing here and in the remainder of \sect{sect_PolyInterpolation} are
\be
F(4a)=\frac{1-\sqrt a\cot\sqrt a}{4a}\eqComma
\qquad
S(a)=\frac{\sin\sqrt a}{\sqrt a}\eqComma
\qquad
C(a)=\cos\sqrt a\eqComma
\label{eq_SFdef}
\ee
so that $q_0=C(a)$.

The rotational parts of \eq{eq_bothCharts} have the lengths $2\sin\sqrt a$ and $2\sqrt a$ and therefore differ by the scalar $\rho=1/S(a)$, i.e., $\dv_\Omega=\rho\,\pO$.
The translational parts are related by eliminating $\tv$ between the two lines of \eq{eq_bothCharts}.
With
\be
r=\frac{1}{q_0}\eqComma
\qquad
\mu=F(4a)\,r\eqComma
\label{eq_rDef}
\ee
and using $\dv_\Omega\tp\vecToSkew{\dv_\Omega}=\zero\tp$, $\dv_\Omega\dv_\Omega\tp=\vecToSkew{\dv_\Omega}^{\,2}+4a\,\identity$ and the scalar identities $\rho\,S(a)=1$, $\rho\,q_0=1-4aF(4a)$ and $\mu\,q_0=F(4a)$, a direct computation gives
\be
\identity-\frac12\vecToSkew{\dv_\Omega}
+F(4a)\,\vecToSkew{\dv_\Omega}^{\,2}
=\left(\rho\identity+\mu\,\dv_\Omega\dv_\Omega\tp\right)
\left(q_0\identity-\frac12\vecToSkew{\pO}\right)\eqDot
\label{eq_verify}
\ee
Applying both sides of \eq{eq_verify} to $\tv$ and comparing with \eq{eq_bothCharts} yields $\dv_U=\left(\rho\identity+\mu\,\dv_\Omega\dv_\Omega\tp\right)\pU$.
Together with $\dv_\Omega=\rho\,\pO$, the relative configuration vector is recovered from the Euler--Rodrigues motion vector through
\be
\dv=\Pbridge(\pV)\,\pV\eqComma
\qquad
\Pbridge(\pV)=
\begin{bmatrix}
\Pbridge_U(\pV) & \zero\\
\zero & \rho\identity
\end{bmatrix}\eqComma
\qquad
\Pbridge_U(\pV)=\rho\identity+\mu\,\dv_\Omega\dv_\Omega\tp\eqComma
\label{eq_bridge}
\ee
with
\be
\dv_\Omega=\rho\,\pV_\Omega \eqComma \qquad \rho=\frac{1}{S(a)} \eqComma
\label{eq_bridgeScalars}
\ee
where $q_0=C(a)$ determines $a=(\phi/2)^2$ from the nodal unknowns, consistent with $a=\tfrac14\dv_\Omega^\top\dv_\Omega$.
Since $q_0$ and $\pV$ are polynomial in the nodal unknowns, \eq{eq_bridge} contains two non-polynomial quantities. The first is the reciprocal $r$ in $\mu$, see \eq{eq_rDef}. The second is the argument $a$, which enters $\rho$ and $\mu$ and is determined by the nodal unknowns only implicitly through $q_0=C(a)$, whose inversion $a=(\arccos q_0)^2$ is not polynomial. 
Both are therefore represented by lifted states, pinned by the identities
\be
q_0\,r=1\eqComma
\qquad
q_0=C(a)\eqDot
\label{eq_thetaConstraints}
\ee
Differentiating \eq{eq_thetaConstraints} with \eq{eq_q0Dot} and $C'(a)=-\frac12S(a)$ gives their rates,
\begin{align*}
\dot r=\frac14\left(\pO\tp\Delta\vv_\Omega\right)r^2\eqComma
%\label{eq_rDot}\\ 
\qquad
\dot a=\frac12\,\rho\left(\pO\tp\Delta\vv_\Omega\right)\eqDot
%\label{eq_aDot}
\end{align*}
The scalar $\rho=1/S(a)$ is bounded by $1\le\rho<\pi/2$ over the admissible range and needs no lift of its own.
The scalar functions $\rho$ and $F(4a)$ are not evaluated from \eq{eq_SFdef}: the numerator and denominator of $F$ both vanish at $a=0$, so the closed form cannot be evaluated at the undeformed configuration, although the function itself is analytic there.
Both are instead evaluated from minimax polynomials in $a$, constructed in \appx{app_MinimaxCoefficients} such that their error over the admissible range $0\le a<(\pi/2)^2$ stays below double-precision round-off.

This matters for the parameterization method, because the invariant manifold is expanded about the equilibrium, where the relative rotation of an element vanishes for a straight, unloaded beam and is small otherwise.
Implementations that keep the closed form switch to a truncated series for small $a$, but the switch makes the expression non-differentiable at the switching point, and the truncated Taylor arithmetic of \sect{sec:implementation-aspects} cannot pass through it.
A polynomial needs no such switch, so the residual can be differentiated to any order at any equilibrium in the admissible range.
The relative configuration vector required for the beam interpolation is thus obtained through polynomial operations only.

%%%%%%%%%%%%%%%%%%%%%%%%%%%%%%%%%%%%%%%%%%%%%%%%%%%%%%%%%%%%%%%%%%%%%%%%%%%%%%%%
\subsubsection{Polynomial interpolation of rigid-body motion}
\label{sect_PolyInterpolationMotion}
%%%%%%%%%%%%%%%%%%%%%%%%%%%%%%%%%%%%%%%%%%%%%%%%%%%%%%%%%%%%%%%%%%%%%%%%%%%%%%%%
%
Following \cite{Sonneville2014}, the motion between the two nodes is interpolated by letting the relative configuration vector vary linearly along the element,
\be
\dv(s)=\frac{s}{L}\,\dv\eqComma
\label{eq_screwInterp}
\ee
so that $\Hm(s)=\Hm_A\exp_{SE(3)}\!\left(\dv(s)\right)$ is the screw motion from $\Hm_A$ to $\Hm_B$.
The strain is then constant along the element and equal to $(\dv-\dv_0)/L$, with $\dv_0$ the relative configuration vector of the undeformed element\footnote{This is why it is measured by $\dv$ and not by $\pV$.}.
 
Since the rotation angle at $s$ is $(s/L)\phi$, the argument of the scalar functions at $s$ is
\be
a(s)=\tfrac14\,\dv_\Omega(s)\tp\dv_\Omega(s)=\frac{s^2}{L^2}\,a\eqDot
\ee
The Euler--Rodrigues parameters at $s$ relative to node $A$ are therefore obtained by inverting \eq{eq_bridge} for the exponential coordinate $\dv(s)=(s/L)\,\Pbridge(\pV)\pV$, with all scalar functions evaluated at $a(s)$,
\be
\pO(s)=S\!\left(a(s)\right)\dv_\Omega(s)\eqComma
\qquad
\pU(s)=\Pbridge_U^{-1}\!\left(a(s),\dv_\Omega(s)\right)\dv_U(s)\eqComma
\qquad
\Pbridge_U^{-1}(a,\dv_\Omega)= S(a)\left(\identity - F(4a)\,\dv_\Omega\dv_\Omega\tp\right)\eqDot
\label{eq_invBridge}
\ee
Using $\dv_\Omega=\rho\,\pO$ and the identity $\rho+4a\mu=r$, the Euler--Rodrigues parameters of the configuration at $s$ relative to node $A$ are
\be
\begin{aligned}
\pO(s)&=\sigma(s)\,\pO\eqComma\\[3pt]
\pU(s)&=\sigma(s)\left[\pU
+r\rho\left(F(4a)-\frac{s^2}{L^2}F\!\left(4a(s)\right)\right)
\left(\pO\tp\pU\right)\pO\right]\eqComma\\[3pt]
q_0(s)&=C\!\left(a(s)\right)\eqComma
\end{aligned}
\label{eq_scaling}
\ee
with
\be
\sigma(s)=\frac{s}{L}\,\rho\,S\!\left(a(s)\right)\eqComma
\qquad
a(s)=\frac{s^2}{L^2}\,a\eqDot
\label{eq_scalingScalars}
\ee
The pair $\left(q_0(s),\pO(s)\right)$ is again an Euler--Rodrigues pair, since $C^2+aS^2=1$ holds at $a(s)$ as it does at $a$, so it satisfies \eq{eq_erpUnit} identically and no renormalisation is required.

The rotation at $s$ follows directly from this pair.
Its position additionally requires the translation $\tv(s)$ of the cross-section at $s$ relative to node $A$.
Unlike the nodal translation $\tv$ in \eq{eq_relativeErpFromNodes}, it cannot be computed from nodal positions and is therefore recovered from $\pU(s)$ by inverting the last relation of \eq{eq_relativeErpFromNodes},
\be
\tv(s)=\Tm_\Omega\tp\,\pU(s)\eqComma
\qquad
\Tm_\Omega=r(s)\,\identity-\frac12\vecToSkew{\pO(s)}+\frac{r(s)}{4}\vecToSkew{\pO(s)}^{\,2}\eqComma
\qquad
r(s)=\frac{1}{q_0(s)}\eqDot
\label{eq_Rtheta}
\ee
The rotation and the position along the element are thus
\be
\begin{aligned}
\RotMat(s)&=\RotMat_A\,\RotMat\!\left(q_0(s),\pO(s)\right)\eqComma\\[3pt]
\rv(s)&=\rv_A+\RotMat_A\,\tv(s)\eqComma
\end{aligned}
\label{eq_field}
\ee
where $\rv_A$ is the position of node $A$ and $\RotMat_A=\RotMat(\qVec_{\EpSca,A})$ its rotation matrix.
The field is exact at the nodes, reducing to node $A$ at $s=0$ and to node $B$ at $s=L$.
If the rotation of the cross-section at $s$ is required in terms of global Euler parameters, these follow from the quaternion product of $\qVec_{\EpSca,A}$ with the relative Euler parameters $\left(q_0(s),\tfrac12\pO(s)\right)$,
\be
\EpSca_0(s)=\EpSca_{0,A}\,q_0(s)-\frac12\,\EpVec_A\tp\pO(s)\eqComma
\qquad
\EpVec(s)=\frac12\,\EpSca_{0,A}\,\pO(s)+q_0(s)\,\EpVec_A+\frac12\,\EpVec_A\times\pO(s)\eqDot
\label{eq_EpAtStation}
\ee

No exponential and no logarithm is evaluated in \eqs{eq_scaling}{eq_field}.
The remaining scalar functions $S$, $C$ and $F$ are evaluated at $a(s)$ for each required cross-section and at $a$ once per element, using the minimax polynomials in \appx{app_MinimaxCoefficients}.
The velocity field required for the inertia terms follows from the same relative configuration vector through the strain tangent of \cite{Sonneville2014}, which is likewise polynomial in $a$.
The element residual is collected in \appx{app_PolyBeam}.

Each element carries the lifted scalars $a$ and $r$ of \eq{eq_thetaConstraints}, which enter \eq{eq_bridge}.
In addition, the position field \eqs{eq_Rtheta}{eq_field} requires the reciprocal $r_k=r(s_k)$ at each quadrature point $s_k$, which is lifted in the same way.
They are collected in $\mathbf{y}_e=(a,r,r_1,\ldots,r_{n_g})$.
Each lifted scalar is added to the configuration variables as a coordinate without inertia and is pinned at position level by its identity,
\be
\mathbf{c}_e\bigl(\qVec_{\EpSca,A},\qVec_{\EpSca,B},\mathbf{y}_e\bigr)=
\begin{bmatrix}
q_0-C(a)\\
q_0\,r-1\\
C\!\left(a(s_1)\right)r_1-1\\
\vdots\\
C\!\left(a(s_{n_g})\right)r_{n_g}-1
\end{bmatrix}=\zero\eqComma
\label{eq_liftedIdentities}
\ee
with $q_0$ from \eq{eq_relativeErpFromNodes} and $a(s_k)$ from \eq{eq_scalingScalars}.
These identities are included in the joint constraints \eq{eq_phi}, each with its own Lagrange multiplier, so that the lifted scalars are determined algebraically by the nodal unknowns.

%%%%%%%%%%%%%%%%%%%%%%%%%%%%%%%%%%%%%%%%%%%%%%%%%%%%%%%%%%%%%%%%%%%%%%%%%%%%%%%%
\subsection{Configuration-space constraints}
\label{sec:csp-constraint}
%%%%%%%%%%%%%%%%%%%%%%%%%%%%%%%%%%%%%%%%%%%%%%%%%%%%%%%%%%%%%%%%%%%%%%%%%%%%%%%%
%
The vectorization \eq{eq_DefVecOperator} represents $q\in G$ by its vectorized form $\qVec=\vec q\in\RealNumbers^{\dimq}$.
Once $\qVec$ is replaced by a truncated polynomial parameterization $\qVec(\pv)$ in the reduced coordinates $\pv$, finite-order truncation no longer guarantees that $\vecOp^{-1}(\qVec)\in G$.
Inspired by \cite{Betsch2001}, we introduce $\dimCsp = \dimq - \dimG$ independent configuration space constraints
\be
\phib\csp(\qVec):=\operatorname{vec}\!\big(q^{-1}q-e\big) \in \RealNumbers^{\dimCsp } \eqComma
\label{eq_Phi-csp-def}
\ee
which are coupled by a correction term added directly to the kinematic reconstruction \eq{eq_KinRecEqCorr} using the Lagrange multiplier $\lagMulKinb\in\RealNumbers^{\dimCsp}$:
\begin{equation}
\dot{\qVec} = \Qm(\qVec)\vv + \Bm\tp\csp(\qVec)\lagMulKinb \in \RealNumbers^{\dimq} \eqComma
\label{eq_KinRecEqCorr}
\end{equation}
where
\be
\Bm_{\mathrm{csp}}(\qVec) = \frac{\partial \phib\csp(\qVec)}{\partial \qVec} \in \RealNumbers^{\dimCsp \times \dimq} \eqDot
\ee

From this point onward, \eq{eq_KinRecEqCorr} supersedes \eq{eq_ODELieGroupe} as the kinematic reconstruction equation used in the remainder of this paper.

Following the previous section, we illustrate the configuration-space
constraints \eq{eq_Phi-csp-def} for the Lie group $SO(3)$ and $S^3$.
The rotation-matrix representation yields six quadratic constraints, the independent entries of the symmetric matrix $\RotMat\tp\RotMat-\Im_3$,
\begin{equation}
\phib\csp(\qVec_\RotMatScalar) = \vec{\RotMat\tp\RotMat-\Im_3} = \zero \eqDot
\label{eq_app_Phi_SO3}
\end{equation}
For Euler parameters, these six constraints reduce to the single quadratic
constraint
\be
\phib\csp(\qVec_\EpSca) = \EpSca_0^2+\EpVec\tp\EpVec-1 \eqComma
\label{eq_cspER}
\ee
so that four configuration variables are constrained by a single equation,
whereas the rotation-matrix representation requires six constraints for nine
variables.

%%%%%%%%%%%%%%%%%%%%%%%%%%%%%%%%%%%%%%%%%%%%%%%%%%%%%%%%%%%%%%%%%%%%%%%%%%%%%%%%
\section{Computing invariant manifolds based on the residual formulation}\label{sect_ResidualRepresentation}
%%%%%%%%%%%%%%%%%%%%%%%%%%%%%%%%%%%%%%%%%%%%%%%%%%%%%%%%%%%%%%%%%%%%%%%%%%%%%%%%
%
This section formulates the parameterization method directly on the residual of the constrained equations of motion, without converting them into an explicit first-order form and without decomposing the nonlinear forces into multilinear maps.
\sect{sect_residualFormulation} assembles the governing equations of \sect{sect_LieGroupMBD} into a single residual, which is the only object required by the subsequent construction. 
\sect{sec:fixed-points} computes the equilibrium and introduces shifted coordinates.
\sect{sec:linearization} linearizes the residual and identifies the master spectral subspace from the finite eigenvalues of the resulting generalized eigenvalue problem.
\sect{sec:computing-parameterizations} states the invariance equation and derives the conventional cohomological equation directly from the residual, showing that its
forcing coincides with that of the explicit formulation \cite{jainhaller2021how}.
Consequently, the established recursive solution procedure applies without modification, as summarized in \sect{sec:cohomologicalSolve}.
Finally, \sect{sect_NonAutonomous} extends the formulation to periodically forced systems.

\subsection{Assembled residual}\label{sect_residualFormulation}
%%%%%%%%%%%%%%%%%%%%%%%%%%%%%%%%%%%%%%%%%%%%%%%%%%%%%%%%%%%%%%%%%%%%%%%%%%%%%%%%
%
%
Conventional computational methods for computing invariant manifolds of fixed points in physical coordinates assume that the governing equations already sit in an explicit form and require, in particular, explicit expressions of nonlinear force terms in the form of multilinear maps, see \cite{jainhaller2021how, li2023model, han2025model}. 
For general three-dimensional constrained multibody systems with rigid and flexible bodies, different types of kinematic joints, potentially user-defined force laws and auxiliary dynamics, obtaining the governing equations non-intrusively in a form that these computational methods can consume, is both laborious and error-prone.
The approach developed in the remainder of this section relies solely on the evaluation of residuals and their Fréchet derivatives, eliminating the need for explicit expressions of nonlinear force terms in the form of multilinear maps. 
As a consequence, the proposed computation approach is independent of the underlying multibody formulation, finite-element discretization, and matrix Lie group used to construct the governing equations. 
We therefore begin by assembling the governing equations \eqs{eq_DynEOM}{eq_phi} and \eq{eq_KinRecEqCorr} of \sect{sect_LieGroupMBD} into a single residual, which is the only object the subsequent construction requires. 

The computational representation introduced in \sect{sect_LieGroupMBD} lets us collect every unknown quantity of the multibody system into a single computational state vector 
\be
\xv = \begin{bmatrix}
\qVec\tp &
\vv\tp &
\tmu\tp &
\tlambda\tp &
\teta\tp &
\end{bmatrix}\tp
\in X
\eqComma
\label{eq:ComputationalState}
\ee
where $X =
\RealNumbers^{\dimq}
\times
\RealNumbers^{\dimG}
\times
\RealNumbers^{\dimAux}
\times
\RealNumbers^{\dimJnt}
\times
\RealNumbers^{\dimCsp}$.
We first collect the equilibrium equations of forces and momenta \eq{eq_DynEOM} into the mechanical residual
\be
\MRes( \qVec,\vv,\dot\vv,\lagMulb,\auxVarb) = \Mb(\qVec)\vbDot + \gb\ine(\qVec, \vb) + \gb\int(\qVec, \vb) + \gb\aux(\qVec,\vb,\auxVarb,t) + \Bm\jnt\tp(\qVec)\lagMulb  - \epsilon\,\gb_\textrm{ext}(\qVec, \vb, \Omega t) = \Null \eqDot
\label{eq_MechanicalResidual}
\ee
Together with the configuration-space constraints \eq{eq_Phi-csp-def}, the formulation in \sect{sect_LieGroupMBD} consists of five components: the mechanical residual \eq{eq_MechanicalResidual}, the kinematic reconstruction \eq{eq_KinRecEqCorr}, the auxiliary-variable evolution \eq{eq_auxVarODE}, the joint constraint \eq{eq_phi}, and the configuration-space constraint \eq{eq_Phi-csp-def}. 
These are collected into the assembled residual
\be\label{eq_CompleteResidual}
\CRes(\xv,\dot{\xv})
=
\begin{bmatrix}
\MRes( \qVec,\vv,\dot\vv,\lagMulb,\auxVarb)
\\[1.0em]
\dot{\qVec} - \Qm(\qVec)\vv - \Bm\tp\csp(\qVec)\lagMulKinb
\\[0.8em]
\dot\auxVarb - \ab(\qVec,\vb,\auxVarb,t)
\\[0.8em]
\phib\jnt(\qVec)
\\[0.8em]
\phib\csp(\qVec)
\end{bmatrix} =\Null \eqDot
\ee

The residual \eq{eq_CompleteResidual} provides a compact representation of the governing equations. 
The computational method developed below depends exclusively on evaluations of \eq{eq_CompleteResidual} together with its Fréchet derivatives
$D_{\xv}\CRes$ and $D_{\dot{\xv}}\CRes$.
Since all spaces considered are finite-dimensional, the Fréchet derivatives coincide with the corresponding Jacobian matrices, i.e,
\be
D_{\xv}\CRes = \frac{\partial \CRes(\xv,\dot{\xv})}{\partial \xv} \eqComma \qquad D_{\dot{\xv}}\CRes = \frac{\partial \CRes(\xv,\dot{\xv})}{\partial \dot{\xv}} \eqDot
\ee

The particular collection of residual equations adopted in \eq{eq_CompleteResidual} merely reflects the structure of the present multibody formulation.
More generally, additional algebraic or differential equations describing arbitrary physical subsystems may be incorporated by augmenting the assembled residual. 
Since the proposed parameterization method depends exclusively on evaluations of this assembled residual and its Fréchet derivatives, all such equations are treated uniformly, irrespective of their mathematical origin or physical interpretation.

%%%%%%%%%%%%%%%%%%%%%%%%%%%%%%%%%%%%%%%%%%%%%%%%%%%%%%%%%%%%%%%%%%%%%%%%%%%%%%%%

%%%%%%%%%%%%%%%%%%%%%%%%%%%%%%%%%%%%%%%%%%%%%%%%%%%%%%%%%%%%%%%%%%%%%%%%%%%%%%%%
\subsection{Computation of fixed points}
\label{sec:fixed-points}
%%%%%%%%%%%%%%%%%%%%%%%%%%%%%%%%%%%%%%%%%%%%%%%%%%%%%%%%%%%%%%%%%%%%%%%%%%%%%%%%
%
We seek to locally approximate the invariant manifolds associated with the fixed point $\zero\in\RealNumbers^{\dimSS}$ of \eq{eq_CompleteResidual}, where $\dimSS = \textrm{dim}(X)$.
To obtain this fixed point, we first compute an equilibrium configuration $\xv\eqp$ satisfying $\CRes(\xv\eqp,\zero)=\zero$ and subsequently shift the coordinates such that $\xv\eqp$ is mapped to the origin.
Therefore, we define the shifted residual
\be
\CResS(\zv,\dot{\zv}) = \CRes(\xv\eqp + \zv,\dot \zv) \eqComma
\label{eq_ResCompleteShifted}
\ee
which satisfies $\CResS(\zero,\zero)=\CRes(\xv\eqp,\zero)=\zero$.
Hence, $\zv=\zero$ is a fixed point of $\CResS$, corresponding to the equilibrium state $\xv=\xv\eqp$ through the translation~\mbox{$\xv = \xv\eqp + \zv$}.
The subsequent linearization and the recursive computation of the invariant manifold's Taylor coefficients are carried out on $\CResS$ directly in terms of $\zv$.

%%%%%%%%%%%%%%%%%%%%%%%%%%%%%%%%%%%%%%%%%%%%%%%%%%%%%%%%%%%%%%%%%%%%%%%%%%%%%%%%
\subsection{Linearization and master spectral subspace}
\label{sec:linearization}
%%%%%%%%%%%%%%%%%%%%%%%%%%%%%%%%%%%%%%%%%%%%%%%%%%%%%%%%%%%%%%%%%%%%%%%%%%%%%%%%
%
Introducing the abbreviations
\begin{equation}
\SSA =D_{\zv}\CResS(\zero,\zero)\eqComma \quad\quad \SSB = D_{\dot\zv}\CResS(\zero,\zero)\eqComma
\label{eq_residual_linearization_AB}
\end{equation}
the linearization of $\CResS$ at the fixed point $\zv=\zero$ reads~\mbox{$\CResS(\zv,\dot\zv) = \SSA\zv + \SSB\dot\zv + \text{higher order terms}$}, and defines a linear system whose invariant subspaces are given by the eigenspaces of the generalized eigenvalue problem\footnote{This work focuses only on non-defective master eigenspaces. Namely, the linear reduced dynamics are diagonalizable.}
\be\label{eq_genEigProblem}
(\SSA + \lambda_j \SSB)\wv_j=\zero \eqComma \quad\quad j = 1 \ldots \dimSS \eqDot
\ee
For each distinct eigenvalue $\lambda_j$, $\SpecSubSpace_j \subset \RealNumbers^{\dimSS}$ denotes the eigenspace spanned by the real and imaginary parts of the generalized eigenvector $\wv_j\in\ComplexNumbers^{\dimSS}$.
We construct the SSM around a master spectral subspace of dimension $\dimIVM \ll \dimSS$.
In general, the linear dynamics on this subspace is described by a matrix whose eigenvalues form the master spectrum. 
We represent it in the eigenvector basis, which gives
\be\label{eq_LambdaMatrix}
\mathbf{\Lambda} = \mathrm{diag}(\lambda_1,\ldots,\lambda_{\dimIVM}) \eqDot
\ee
For mechanical systems with algebraic constraints, matrix $\SSB$ is singular, cf. \cite{li2023model}.
Rows of $\CResS$ corresponding to such constraints depend on $\zv$ but never on $\dot\zv$, so the pencil $\SSA +\lambda \SSB$ carries eigenvalues at infinity alongside the finite ones.
Only the finite ones enter the master spectrum.

%%%%%%%%%%%%%%%%%%%%%%%%%%%%%%%%%%%%%%%%%%%%%%%%%%%%%%%%%%%%%%%%%%%%%%%%%%%%%%%%
\subsection{Invariance equation and its recursive solution}
\label{sec:computing-parameterizations}
%%%%%%%%%%%%%%%%%%%%%%%%%%%%%%%%%%%%%%%%%%%%%%%%%%%%%%%%%%%%%%%%%%%%%%%%%%%%%%%%
%
We seek an $\dimIVM$-dimensional invariant manifold $\ivm(\SpecSubSpace)$ tangent to the spectral subspace $\SpecSubSpace$, parameterized by
\be
\Wm:\ComplexNumbers^{\dimIVM}\rightarrow X \eqComma \qquad \zv(t)=\Wm(\pv(t)) \eqComma
\label{eq_StatesFromW}
\ee
where $\pv\in\ComplexNumbers^{\dimIVM}$ collects the reduced coordinates, whose evolution is governed by the reduced dynamics
\be
\Rm:\ComplexNumbers^{\dimIVM}\rightarrow\ComplexNumbers^{\dimIVM} \eqComma \qquad \dot{\pv}=\Rm(\pv) \eqDot
\label{eq_ReducedDynamics}
\ee
Substituting \eq{eq_StatesFromW} and its time derivative
\be\label{eq_States_t}
\dot \zv = \left(D_{\pv}\Wm(\pv)\right)\dot{\pv} = \left(D_{\pv} \Wm(\pv)\right)\Rm(\pv)
\ee
into \eq{eq_ResCompleteShifted} yields the invariance equation
\be
\CResS\bigl(\Wm(\pv), \left(D_{\pv}\Wm(\pv)\right)\Rm(\pv)\bigr)=\zero \eqComma
\label{eq_invariance}
\ee
a single functional equation for the pair $(\Wm,\Rm)$.

%%%%%%%%%%%%%%%%%%%%%%%%%%%%%%%%%%%%%%%%%%%%%%%%%%%%%%%%%%%%%%
\subsubsection{Multi-index expansion}
%%%%%%%%%%%%%%%%%%%%%%%%%%%%%%%%%%%%%%%%%%%%%%%%%%%%%%%%%%%%%%
%
The multi-index~\mbox{$\alpha=(\alpha_1,\ldots,\alpha_{\dimIVM})\in\NaturalNumbers_0^{\dimIVM}$} defines the monomial~\mbox{$\pv^\alpha=p_1^{\alpha_1} \dotsm p_{n_r}^{\alpha_{n_r}}$} of total degree~\mbox{$|\alpha|=\alpha_1 + \dotsm + \alpha_{n_r}$}. We expand the unknown functions as polynomials, as follows
{\setlength{\abovedisplayskip}{3pt}
 \setlength{\belowdisplayskip}{3pt}
\be
\Wm(\pv)=\sum_{|\alpha|\ge1}\Wm_\alpha\pv^\alpha \eqComma \qquad
\Rm(\pv)=\sum_{|\alpha|\ge1}\Rm_\alpha\pv^\alpha \eqComma
\label{eq:multiindex_expansions}
\ee}
with $\Wm_\alpha\in\ComplexNumbers^{\dimSS}$ and $\Rm_\alpha\in\ComplexNumbers^{\dimIVM}$.
Given the $j$-th standard unit vector $\ev_j$, the velocity field \eq{eq_States_t} follows as
{\setlength{\abovedisplayskip}{3pt}
 \setlength{\belowdisplayskip}{3pt}
\be
D_{\pv}\Wm(\pv)\,\Rm(\pv)
=\sum_{|\alpha|\ge1}\sum_{j=1}^{\dimIVM}\alpha_j\,\Wm_\alpha\,\pv^{\alpha-\ev_j}\,R_j(\pv) \eqComma
\label{eq:multiindex_velocity}
\ee}
where $R_j(\pv)$ is the $j$-th component of $\Rm(\pv)$.
Substituting \eqs{eq:multiindex_expansions}{eq:multiindex_velocity} into \eq{eq_invariance} and using analyticity of the residual, the coefficient of every monomial $\pv^\alpha$ must vanish independently.

%%%%%%%%%%%%%%%%%%%%%%%%%%%%%%%%%%%%%%%%%%%%%%%%%%%%%%%%%%%%%%
\subsubsection{First-order coefficients}
%%%%%%%%%%%%%%%%%%%%%%%%%%%%%%%%%%%%%%%%%%%%%%%%%%%%%%%%%%%%%%
%
For $|\alpha|=1$ all nonlinear residual terms vanish, yielding
{\setlength{\abovedisplayskip}{3pt}
 \setlength{\belowdisplayskip}{3pt}
\be
\SSA\Wm_{\ev_j} + \lambda_j\SSB\Wm_{\ev_j} = \zero \eqComma \qquad j=1,\ldots,\dimIVM \eqComma
\label{eq:first_order_multiindex}
\ee}
which is identical to the generalized eigenvalue problem \eq{eq_genEigProblem}.

%%%%%%%%%%%%%%%%%%%%%%%%%%%%%%%%%%%%%%%%%%%%%%%%%%%%%%%%%%%%%%
\subsubsection{Cohomological equation}
\label{sec:cohomologicalSolve}
%%%%%%%%%%%%%%%%%%%%%%%%%%%%%%%%%%%%%%%%%%%%%%%%%%%%%%%%%%%%%%
%
Consider now an arbitrary coefficient with $|\alpha|=k\ge2$.
Only two contributions contain unknowns of order $k$: the linear action of $\SSA$ on $\Wm_\alpha$, and the product of the first-order reduced dynamics with $D_{\pv}\Wm_\alpha$.
The latter contributes, by \eq{eq:multiindex_velocity} with the leading order term $\Rm_1=\mathbf{\Lambda}$,
{\setlength{\abovedisplayskip}{3pt}
 \setlength{\belowdisplayskip}{3pt}
\be
\sum_{j=1}^{\dimIVM}\alpha_j\lambda_j\Wm_\alpha=\lambda_\alpha\Wm_\alpha \eqComma \qquad
\lambda_\alpha=\sum_{j=1}^{\dimIVM}\alpha_j\lambda_j \eqDot
\label{eq:weighted_eigenvalue}
\ee}
All remaining contributions involve only monomials that precede (are computed before)~$\alpha$ and are therefore known recursively.
Collecting them into the forcing vector $\Nm_\alpha$ yields
\be
(\SSA+\lambda_\alpha\SSB)\Wm_\alpha + \SSB\Wm_1\Rm_\alpha = -\Nm_\alpha \eqDot
\label{eq:cohomological_multiindex}
\ee
Denoting the truncations of the two expansions to orders strictly below $\alpha$ by
\be
\Wm^{<\alpha}(\pv) = \sum_{0<\beta<\alpha} \Wm_\beta\pv^\beta \eqComma \qquad
\Rm^{<\alpha}(\pv) = \sum_{0<\beta<\alpha} \Rm_\beta\pv^\beta
\ee
where the multi-index inequalities hold component-wise, the corresponding forcing vector becomes
\be
\Nm_\alpha = \Bigl[\;\CResS\!\bigl(\Wm^{<\alpha}(\pv),\;\bigl(D_{\pv}\Wm^{<\alpha}(\pv)\bigr)\Rm^{<\alpha}(\pv)\bigr)\;\Bigr]_\alpha \eqComma
\label{eq:DefinitionNalpha}
\ee
where $[\,\bullet\,]_\alpha$ extracts the coefficient of $\pv^\alpha$.
Note that \eq{eq:DefinitionNalpha} uses the complete residual, not only its nonlinear part.
Because $\Wm^{<\alpha}$ contains no terms of multi-index $\alpha$, the contribution $\SSA\Wm^{<\alpha}$ vanishes at that monomial.
The velocity argument, however, does contribute, since \eq{eq:multiindex_velocity} generates terms of multi-index $\alpha$ from products of strictly lower-order coefficients.
The single evaluation \eq{eq:DefinitionNalpha} therefore collects both the nonlinear compositions and these lower-order velocity products, and depends exclusively on coefficients preceding $alpha$, which is what renders the recursion triangular.

Crucially, $\Nm_\alpha$ is never assembled from explicitly derived multilinear coefficient arrays.
It is obtained by evaluating the assembled residual \eq{eq_CompleteResidual} once on truncated polynomial arguments and extracting the required coefficient, an operation performed in the implementation by automatic differentiation in a truncated polynomial algebra, cf.~\sect{sec:implementation-aspects}.
The arguments $\Wm_{<|\alpha|}$ and $\Rm_{<|\alpha|}$ of the residual in \eq{eq:DefinitionNalpha} involve only previously computed coefficients, so the residual polynomial is computed once per monomial. 
Since these arguments coincide for all multi-indices of the same order, the residual may alternatively be evaluated once per order, which supplies $\Nm_\alpha$ for all $\binom{k+n_r-1}{k}$ multi-indices with $|\alpha|=k$ from a single evaluation.

\eq{eq:cohomological_multiindex} is precisely the conventional cohomological equation of the parameterization method, cf.~\cite{jainhaller2021how, li2023model}, obtained here directly from the implicit residual formulation, without converting the governing equations into an explicit first-order system.
It remains to verify that the forcing \eq{eq:DefinitionNalpha} agrees with its counterpart in the explicit setting of \cite{jainhaller2021how}, where the system is written as $\Bm\dot{\zv}=\Am\zv+\Fm(\zv)$.
Since the shifted residual $\CResS$ and that first-order form describe the same equations of motion, they agree identically,
\be
\SSA\zv + \SSB\dot\zv + \sum_{k\ge2}\CResS_k(\zv,\dot\zv) \equiv \Bm\dot\zv - \Am\zv - \Fm(\zv) \eqComma \quad \forall(\zv,\dot\zv) \eqDot
\ee
Identifying $\Am=-\SSA$ and $\Bm=\SSB$ and cancelling the shared linear part leaves
\be
\sum_{k\ge2}\CResS_k(\zv,\dot\zv) \equiv -\Fm(\zv) \eqComma \quad \forall(\zv,\dot\zv) \eqDot
\label{eq_CResSFIdentity}
\ee
Evaluating \eq{eq_CResSFIdentity} at $\zv=\Wm(\pv)$ and $\dot\zv=(D_{\pv}\Wm(\pv))\Rm(\pv)$ identifies two power series in $\pv$, which therefore have identical coefficients at every multi-index.
Moreover, $\CResS_k$ is homogeneous of degree $k\ge2$ in arguments that both start at degree one in $\pv$, so any term containing $\Wm_\alpha$ or $\Rm_\alpha$ has degree at least $|\alpha|+k-1>|\alpha|$.
The coefficient at $\alpha$ is thus unaffected by the truncation to $\Wm^{<\alpha}$, $\Rm^{<\alpha}$, and the nonlinear part of \eq{eq:DefinitionNalpha} equals $-(\Fm\circ\Wm)_\alpha$.
Hence $\Nm_\alpha$ coincides with the forcing of \cite{jainhaller2021how}, and no separate equivalence proof is required.
Therefore, the recursive solution procedure of \cite{jainhaller2021how} applies without modification.

\subsection{Solving the cohomological Equation}
%\label{sec:cohomologicalSolve}
%
By the equivalence established in \sect{sec:cohomologicalSolve}, the coefficients $\Wm_\alpha$ and $\Rm_\alpha$ satisfying the cohomological \eq{eq:cohomological_multiindex} can be obtained by the same recursive procedure proposed in \cite{jainhaller2021how}.
We refer to~\cite{jainhaller2021how} for the full derivation and to \cite{haller2016nonlinear} for the underlying existence and uniqueness theory of the spectral submanifold itself.

\subsection{Extension to periodically forced systems}
\label{sect_NonAutonomous}
The recursive solution procedure of the invariance \eq{eq_invariance} developed in \sect{sec:computing-parameterizations}
can be extended to periodically forced systems in accordance to \cite{jainhaller2021how,li2023model,thurnherr2023nonautonomous}.

Periodic forcing is incorporated by introducing an excitation phase $\phi$ governed by the trivial ODE
\be\label{eq_trivialODEExcitation}
\dot{\phi}=\Omega \eqComma
\ee
which turns the explicitly time-dependent forcing $\gv\ext(\qVec,\vv,\Omega t)$ into the autonomous form $\gv\ext(\qVec,\vv,\phi)$. 
Unlike $\qVec$ or $\auxVarb$, $\phi$ is not parametrized by $\Wm$. 
Instead it is carried as an independent argument alongside $\pv$, i.e.,
\be
\zv=\Wm(\pv,\phi) \eqComma
\qquad
\dot{\pv}=\Rm(\pv,\phi) \eqComma
\ee
with $\phi$ itself evolving directly according to \eq{eq_trivialODEExcitation}, cf.~\cite{jainhaller2021how,li2023model,thurnherr2023nonautonomous}.
The residual row $\dot\phi-\Omega=0$ this augmentation would otherwise add is therefore satisfied identically by construction. 
Applying the chain rule to $\zv=\Wm(\pv,\phi)$ gives
\be
\dot{\zv} = \left(D_{\pv}\Wm(\pv,\phi)\right)\Rm(\pv,\phi) + \Omega D_{\phi}\Wm(\pv,\phi) \eqDot
\label{eq:NonAutonomousVelocity}
\ee
Substituting \eq{eq:NonAutonomousVelocity} into the assembled residual
\eq{eq_ResCompleteShifted} gives the non-autonomous invariance equation
\be
\CResS \!\left( \Wm,
\left(D_{\pv}\Wm\right)\Rm + \Omega D_{\phi}\Wm \right) = \zero \eqComma
\label{eq:NonAutonomousInvariance}
\ee
where the only change relative to the autonomous case is the additional term $\Omega D_\phi\Wm$.

$\CResS(\zv,\dot\zv)$ and the non-autonomous first-order form $\Bm\dot \zv= \Am \zv + \Fm(\zv) + \epsilon \Fm^{\textrm{ext}}(\zv,\phi)$ of \cite{jainhaller2021how} describe equivalent representations of the same forced equations of motion, with the forcing terms identified as $\gv\ext(\qVec,\vv,\phi)=\Fm^{ext}(\zv,\phi)$, alongside the identification $\Am=-\SSA,\Bm=\SSB$ of \sect{sec:cohomologicalSolve}. 
Comparing the two expressions and canceling the shared linear and forcing terms leaves the identity $\sum_{k=2}^\infty\CResS_k(\zv,\dot\zv)\equiv-\Fm(\zv)$ unchanged. 
The residual identity underlying that section's equivalence proof therefore extends, without modification, to the non-autonomous case.
Consequently, the leading-order-in-$\epsilon$ solution of the invariance \eq{eq:NonAutonomousInvariance} can be obtained as shown in \cite[Eqs.~70--83]{jainhaller2021how}.
The resulting reduced dynamics provide the basis for established SSM-based techniques for forced-response prediction, backbone-curve computation and bifurcation analysis. 
Since these developments are well established, we omit the details and refer the interested reader to \cite{jainhaller2021how,li2023model,thurnherr2023nonautonomous}.

\section{Preservation of constraints by the reduced model}\label{sect_ConstraintPreservation}
The joint constraints \eqs{eq_jntConstrMapping}{eq_jntConstrMapping_t} encode the kinematic relationships between connected bodies.
The configuration-space constraints \eq{eq_Phi-csp-def} enforce that the computational coordinates remain on the Lie group throughout the motion.
A reduced-order model that violates these constraints no longer represents the kinematics of the original mechanical system exactly.
It is therefore important to understand whether the constraints remain satisfied after truncation of the polynomial parameterization.
\sect{sec_LinearConstraints} and \sect{sec_NonLinearConstraints} investigate the extent to which linear and nonlinear constraints are preserved by the resulting reduced-order model.
\sect{sec_LieParameterization} then shows how exactness can be restored for the configuration-space constraints, which are nonlinear, by parameterizing the group block intrinsically.

\subsection{Linear algebraic constraints}\label{sec_LinearConstraints}
Consider a linear algebraic constraint of the form
\begin{equation}
\phib_\bullet(\xv)=\Lm\xv=\zero \eqComma
\label{eq_linear_constraint}
\end{equation}
where $\Lm$ is a constant matrix. 
By definition, the equilibrium satisfies
\begin{equation}
\phib_\bullet(\xv\eqp)=\zero\eqDot
\label{eq_linear_constraint_equilibrium}
\end{equation}
Introducing the shifted coordinates
\begin{equation}
\zv = \xv - \xv\eqp = \sum_{|\alpha|\ge1}\Wm_\alpha\pv^{\alpha},
\label{eq_shifted_coordinates}
\end{equation}
the constraint becomes
\begin{equation}
\phib_\bullet(\xv\eqp + \zv) = \Lm\zv = \sum_{|\alpha|\ge1}\Lm\Wm_\alpha\pv^{\alpha} \eqDot
\label{eq_linear_constraint_shifted}
\end{equation}
Since the corresponding residual row is linear, each recursive coefficient equation is linear in the unknown coefficient $\Wm_\alpha$.
Consequently, no coupling terms involving previously computed coefficients are generated, and the forcing vector $\Nm_\alpha$ of \eq{eq:cohomological_multiindex} vanishes identically in these rows.
The recursive coefficient equation therefore reduces to
\begin{equation}
\Lm\Wm_\alpha=\zero\eqComma\qquad |\alpha|\ge1 \eqDot
\label{eq_linear_constraint_coefficients}
\end{equation}
Hence, every finite-order polynomial approximation of the parameterization satisfies the linear algebraic constraint exactly, irrespective of the truncation order.

The above result has direct implications for constrained multibody systems.
Consider a spherical joint connecting two bodies. 
The corresponding position-level constraint can be written as
\be\label{eq_SphericalJoint}
\phib_{\textrm{sph}}= \yv_1+\RotMat_1\uv_1-\yv_2-\RotMat_2\uv_2=\zero \eqComma
\ee
where $\yv_i \in \RealNumbers^3$ denotes the position of body $i$, $\RotMat_i \in SO(3)$ its rotation matrix, and $\uv_i\in \RealNumbers^3$ the local coordinates of the joint attachment point.
When the entries of the rotation matrices are chosen directly as computational variables (cf. \sect{sect_RotationMatrixComponents}), the spherical-joint constraint is linear in the computational state. 
Consequently, it is preserved exactly by every finite-order polynomial parameterization.
In fact, several classical lower-pair joints can be modeled using linear constraint equations. 
We show this briefly in \appx{sect_LinearLowerPairJoints}.
Redundant linear formulations may likewise be employed whenever they simplify the residual representation, as also discussed in \appx{sect_LinearLowerPairJoints}.

By contrast, if body orientations are represented by Euler parameters (cf. \sect{sect_EulerParameters}), \eq{eq_SphericalJoint} constitutes a quadratic algebraic constraint. 
Quadratic constraints are no longer preserved exactly by a truncated polynomial parameterization but only up to the approximation order established in \sect{sec_NonLinearConstraints}.

\subsection{Nonlinear constraints}\label{sec_NonLinearConstraints}
For notational simplicity, let $\phib_{\bullet}(\xv)=\zero$ denote an arbitrary nonlinear constraint represented by a residual row of the assembled invariance equation, and let $\dot{\phib}_{\bullet}(\xv,\dot\xv)=\zero$ denote the corresponding velocity-level constraint.
\sect{sect_ResidualRepresentation} established that the recursive solution satisfies every equation of the assembled residual for all multi-indices with $|\alpha|\le k$. 
Hence, the residual of every nonlinear constraint also vanishes through polynomial degree $k$.
Consequently,
\begin{align}
&\phib_{\bullet}\left(
\xv\eqp+\Wm(\pv) \right) = \bigO\!\left(|\pv|^{k+1}\right) \eqComma \\
&\dot{\phib}_{\bullet}\left(\xv\eqp+\Wm(\pv), \left(D_{\pv}\Wm(\pv)\right)\Rm(\pv) \right) = \bigO\!\left(|\pv|^{k+1}\right) \eqDot
\label{eq:NonlinearConstraintError}
\end{align}
Therefore, every analytic nonlinear constraint represented by a residual row of the assembled invariance equation is satisfied with a residual of $\mathcal O(|\pv|^{k+1})$.

\subsection{Rotation-preserving parameterization}
\label{sec_LieParameterization}
%%%%%%%%%%%%%%%%%%%%%%%%%%%%%%%%%%%%%%%%%%%%%%%%%%%%%%%%%%%%%%%%%%%%%%%%%%%%%%%%
%
\sect{sec_NonLinearConstraints} showed that nonlinear constraints are satisfied only up to a residual of $\bigO(|\pv|^{k+1})$.
The configuration-space constraints \eq{eq_Phi-csp-def} are of this kind, since a reduced model of order $k$ returns configuration variables that leave the matrix Lie group at $\bigO(|\pv|^{k+1})$.
This subsection removes that residual by parameterizing the group block intrinsically, so that group membership holds to machine precision at every order, while the exact preservation of linear joint constraints established in \sect{sec_LinearConstraints} is retained.

%%%%%%%%%%%%%%%%%%%%%%%%%%%%%%%%%%%%%%%%%%%%%%%%%%%%%%%%%%%%%%
\subsubsection{Tangent-space parameterization}\label{sec_LieTangentParam}
%%%%%%%%%%%%%%%%%%%%%%%%%%%%%%%%%%%%%%%%%%%%%%%%%%%%%%%%%%%%%%
%
In what follows, $G$ denotes a general matrix Lie group, of which $\SOThree$ and $\SThree$ are the instances used in this section.
Up to this point, the configuration coordinates $\qVec$, such as the nine entries of each rotation matrix, have been expanded as polynomials in $\pv$, so that group elements are represented as points of the vector space $\RealNumbers^{\dimq}$ in which $G$ is embedded.
We refer to this parameterization as the ambient formulation, in contrast to the rotation-preserving parameterization developed below.

Let $q_0\in G$ denote the group element at the fixed point, with computational vector $\qVec_0$, and let $\kinRecMat_0:=\kinRecMat(\qVec_0)$ and $\Bm_{\textrm{csp},0}:=\Bm_{\textrm{csp}}(\qVec_0)$ denote the kinematic reconstruction matrix \eq{eq_CompRecMat} and the configuration-space constraint Jacobian \eq{eq_Phi-csp-def} evaluated at $q_0$.
The columns of $\kinRecMat_0$ are $\vecOp(q_0\widetilde{\ev}_1),\ldots,\vecOp(q_0\widetilde{\ev}_{\dimG})$, so its range is the vectorized tangent space of $G$ at $q_0$, which coincides with the null space of $\Bm_{\textrm{csp},0}$,
\be
\ker(\Bm_{\textrm{csp},0}) = \mathrm{range}(\kinRecMat_0) \eqDot
\label{eq_TangentBasis}
\ee
\eq{eq_TangentBasis} is fundamental to the linearized formulation.
Without the configuration-space constraints \eq{eq_Phi-csp-def}, the linearized kinematic reconstruction leaves the $\dimq-\dimG$ directions transverse to $G$ unrestrained, which introduces spurious zero eigenvalues into the pencil $\SSA+\lambda\SSB$.
The multipliers $\lagMulKinb$ eliminate these directions, and \eq{eq_TangentBasis} identifies the admissible directions with the tangent space spanned by the columns of $\kinRecMat_0$.
This observation motivates the construction below, in which the configuration coefficients are replaced by $\dimG$ intrinsic coefficients in the tangent space and the multiplier is eliminated by projecting the kinematic-reconstruction rows with $\kinRecMat_0^{\top}$.

Let $\qVec_\alpha$ denote the $\qVec$-block of the unknown coefficient $\Wm_\alpha$ in \eq{eq:cohomological_multiindex}, i.e.\ the order-$\alpha$ coefficient of the configuration.
We replace $\qVec_\alpha$ by intrinsic coefficients $\wv_\alpha\in\ComplexNumbers^{\dimG}$ and reconstruct the configuration as
{\setlength{\abovedisplayskip}{4pt}
 \setlength{\belowdisplayskip}{4pt}
\be
\qVec(\pv)=\vecOp\Bigl(q_0\,\exp_G\!\bigl(\widetilde{\wv(\pv)}\bigr)\Bigr) \eqComma \qquad
\wv(\pv)=\sum_{|\alpha|\ge1}\wv_\alpha\pv^\alpha \eqComma
\label{eq_LieReconstruction}
\ee}
where
{\setlength{\abovedisplayskip}{3pt}
 \setlength{\belowdisplayskip}{3pt}
\be
\exp_G=\operatorname{blockdiag}\left(\exp_{\cSpaceBody_1},
\exp_{\cSpaceBody_2},\ldots,\exp_{\cSpaceBody_{\nB}}\right)
\ee}
is the exponential map of the product manifold $G=\cSpaceBody_1\times\cSpaceBody_2\times\ldots\times\cSpaceBody_{\nB}$, each factor $\cSpaceBody_i$ carrying its own exponential map, such as $\expSOThree$ or $\expSThree$, and translational factors $\RealNumbers^3$ the identity.
The vectorization in \eq{eq_LieReconstruction} acts factor-wise, stacking $\vecOp(q_i)$ over the bodies as in \eq{eq_CompRecMat}.
By construction, $\vecOp^{-1}(\qVec(\pv)) \in G$ for every $\pv$ and every truncation order.

Matching homogeneous orders in \eq{eq_LieReconstruction} and using $\vecOp\bigl(q_0\widetilde{\wv}_\alpha\bigr)=\kinRecMat_0\wv_\alpha$ gives
\be
\qVec_\alpha = \kinRecMat_0\wv_\alpha + \boldsymbol{\Gamma}_\alpha \eqComma
\qquad
\boldsymbol{\Gamma}_\alpha=\vecOp\Bigl(q_0\sum_{n=2}^{|\alpha|}\Bigl[\tfrac{1}{n!}\,\widetilde{\wv}(\pv)^n\Bigr]_\alpha\Bigr) \eqComma
\label{eq_LieReconstructionOrder}
\ee
where $[\,\bullet\,]_\alpha$ extracts the coefficient of $\pv^\alpha$.
Since $\wv(\zero)=\zero$, every power satisfies $\widetilde{\wv}(\pv)^n=\bigO(|\pv|^n)$, so the coefficient of $\pv^\alpha$ in $\exp_G(\widetilde{\wv}(\pv))$ receives contributions only from powers $n\le|\alpha|$, and the series ends.
Hence $\boldsymbol{\Gamma}_\alpha\in\ComplexNumbers^{\dimq}$ is a finite sum of products of $\widetilde{\wv}_\beta$ with $|\beta|<|\alpha|$; no transcendental function is evaluated, and the offline recursion remains polynomial.
\eq{eq_LieReconstructionOrder} is affine in the unknown $\wv_\alpha$ and triangular in the recursion, so the structure of the cohomological equation \eq{eq:cohomological_multiindex} is preserved.

The kinematic-reconstruction block of \eq{eq:cohomological_multiindex}, with $\dimq$ rows and $\dimq$ unknowns inherited from $\qVec_\alpha$, is no longer square once $\qVec_\alpha$ is replaced by $\wv_\alpha$.
Restricted to this block, and using $\vv\eqp=\zero$ and $\lagMulKinb\eqp=\zero$ at the fixed point, \eq{eq:cohomological_multiindex} reads
\be
\lambda_\alpha \qVec_\alpha - \kinRecMat_0 \vv_\alpha - \Bm\tp_{\csp,0}\lagMulKinb_\alpha = -\mathbf{n}_\alpha \eqComma
\label{eq_KinRow}
\ee
where $\vv_\alpha$ and $\lagMulKinb_\alpha$ are the order-$\alpha$ coefficients of $\vv(\pv)$ and $\lagMulKinb(\pv)$, and $\mathbf{n}_\alpha$ collects the known lower-order forcing.
Substituting \eq{eq_LieReconstructionOrder} and left-multiplying by $\kinRecMat_0\tp$ gives
\be
\kinRecMat_0\tp\kinRecMat_0\,(\lambda_\alpha\wv_\alpha-\vv_\alpha) - \kinRecMat_0\tp\Bm\tp_{\csp,0}\lagMulKinb_\alpha= -\kinRecMat_0\tp(\mathbf{n}_\alpha+\lambda_\alpha\boldsymbol{\Gamma}_\alpha) \eqDot
\label{eq_ProjectedKinRow}
\ee
Because $\kinRecMat_0\tp\Bm_{\textrm{csp},0}\tp=\zero$ by \eq{eq_TangentBasis}, the multiplier $\lagMulKinb_\alpha$ drops out of \eq{eq_ProjectedKinRow}.
The projection is thus the classical null-space method \cite{Shabana2010} applied to the saddle-point system \eq{eq:cohomological_multiindex}, in which the multiplier is eliminated rather than solved for.
Since $\kinRecMat_0$ has full column rank $\dimG$, the Gram matrix $\kinRecMat_0\tp\kinRecMat_0$ is invertible, and \eq{eq_ProjectedKinRow} is a square $\dimG\times\dimG$ system for $\wv_\alpha$, coupled to the mechanical residual through $\vv_\alpha$.
Once $\wv_\alpha$ is known, $\qVec_\alpha$ follows from \eq{eq_LieReconstructionOrder}.
The configuration-space constraint block, i.e.\ the $\dimCsp$ rows $\phib_{\textrm{csp}}=\zero$ and the $\dimCsp$ multiplier columns of $\lagMulKinb_\alpha$, is dropped entirely, since the tangent-space condition it enforces is built into \eq{eq_LieReconstruction}.

The reduction defined by \eq{eq_LieReconstructionOrder} and \eq{eq_ProjectedKinRow} is applied to the assembled system \eq{eq:cohomological_multiindex} as a whole.
Both operations are encoded in the constant matrices $\Sm\in\RealNumbers^{\dimRed\times\dimSS}$ and $\Vm\in\RealNumbers^{\dimSS\times\dimRed}$, with $\dimRed=\dimSS-\dimq-\dimCsp+\dimG$, which are assembled once before the recursion begins.
The column map $\Vm$ replaces the $\dimq$ columns of the group block by $\kinRecMat_0$, thereby mapping $\wv_\alpha$ to the configuration coefficients, and omits the $\dimCsp$ columns of $\lagMulKinb$.
The row map $\Sm$ replaces the $\dimq$ kinematic reconstruction rows by their projection with $\kinRecMat_0^{\top}$ and omits the $\dimCsp$ configuration-space constraint rows.
On all other rows and columns, both maps are the identity.
At every order, the reduced system
\be
(\Sm\Lm\Vm)\,\WmRed_\alpha
=
\Sm\bigl(-\Nm_\alpha-\SSB\Wm_1\Rm_\alpha-\Lm\bar{\boldsymbol{\Gamma}}_\alpha\bigr)
\eqComma
\qquad
\Lm=\SSA+\lambda_\alpha\SSB
\eqComma
\label{eq_SLTSystem}
\ee
for the reduced coefficients $\WmRed_\alpha\in\ComplexNumbers^{\dimRed}$, which contain $\wv_\alpha$ in place of the configuration coefficients, is square by construction and is solved directly, as is the unreduced system \eq{eq:cohomological_multiindex}.
Here, $\bar{\boldsymbol{\Gamma}}_\alpha\in\ComplexNumbers^{\dimSS}$ denotes $\boldsymbol{\Gamma}_\alpha$ padded by zeros outside the group block.
In the kinematic reconstruction rows, $\Lm\bar{\boldsymbol{\Gamma}}_\alpha$ reduces to $\lambda_\alpha\boldsymbol{\Gamma}_\alpha$, cf.\ \eq{eq_ProjectedKinRow}.
The full state coefficients are recovered as $\Wm_\alpha=\Vm\WmRed_\alpha+\bar{\boldsymbol{\Gamma}}_\alpha$, which generalizes \eq{eq_LieReconstructionOrder} to the full state.

\subsubsection{Constrained evaluation of the reduced model}
\label{sec_LieConstrainedEmbedding}
%%%%%%%%%%%%%%%%%%%%%%%%%%%%%%%%%%%%%%%%%%%%%%%%%%%%%%%%%%%%%%
%
In the rotation matrix representation, the joint constraints are linear in $\qVec$, and the coefficients $\qVec_\alpha$ of \eq{eq_LieReconstructionOrder} satisfy them exactly for $|\alpha|\le k$ (\sect{sec_LinearConstraints}).
The reconstruction \eq{eq_LieReconstruction}, however, applies the full exponential to the truncated polynomial $\wv(\pv)$ and thereby generates terms of order $|\pv|^{k+1}$ and higher that are not determined by the cohomological equations.
The reconstructed configuration therefore lies exactly on $G$ but violates the joint constraints by $\bigO(|\pv|^{k+1})$, whereas the polynomial evaluation of $\qVec(\pv)$ behaves conversely.

Both properties are obtained simultaneously by a correction map $\mathcal{C}$, defined in \eq{eq_CorrectionMap} below, which maps configurations in a neighborhood of the configuration manifold $\mathcal{M}=\{\qVec:\ \phib\csp(\qVec)=\zero,\ \phib\jnt(\qVec)=\zero\}$ onto $\mathcal{M}$ and reduces to the identity on $\mathcal{M}$.
Since the configuration of every point of the exact SSM lies in $\mathcal{M}$, $\mathcal{C}$ alters the reconstructed configuration by $\bigO(|\pv|^{k+1})$ only.
The corrected parameterization therefore has the same Taylor coefficients up to order $k$ and the same reduced dynamics, and it satisfies the configuration-space and joint constraints exactly for every $\pv$ and every truncation order.
The configuration of the reduced model thus follows from the polynomial computed by the recursion through two fixed maps that act beyond order $k$ only,
\be
\pv \;\longmapsto\; \wv(\pv) \;\xrightarrow{\;\vecOp\left(q_0\exp_G(\widetilde{\,\cdot\,})\right)\;}\; \vecOp(G) \;\xrightarrow{\;\mathcal{C}\;}\; \mathcal{M} \eqDot
\label{eq_EmbeddingChain}
\ee

The construction of $\mathcal{C}$ follows from the elementary constraints of which lower-pair joints are composed (\appx{sect_LinearLowerPairJoints}).
The point-coincidence constraint \eq{eq_SphericalJoint} involves the positions of both bodies, whereas the vector-coincidence constraint \eq{eq:VectorCoincidence} involves rotations only and, in the tangent coordinates of \eq{eq_LieReconstruction}, reads
\be
\RotMat_{a,0}\exp\bigl(\widetilde{\wv}_a\bigr)\,\sv_a = \RotMat_{b,0}\exp\bigl(\widetilde{\wv}_b\bigr)\,\sv_b \eqDot
\label{eq_VectorCoincidenceTangent}
\ee
The truncated polynomials $\wv_a(\pv)$ and $\wv_b(\pv)$ satisfy this nonlinear relation exactly only in special cases, namely for joints that lock the relative orientation, where it reduces to $\wv_b=\RotMat_{b,0}\tp\RotMat_{a,0}\wv_a$, and for revolute joints whose axes are parallel to all joint axes on the path to the ground, as in planar systems.
In general, it is violated by $\bigO(|\pv|^{k+1})$, and the violation cannot be removed by correcting positions.

The correction of each joint is therefore assigned to the body farther from the ground, the child $b$, while the parent $a$ is kept fixed, and the joints are processed such that every parent precedes its children.
This ordering only fixes the sequence of the corrections; no relative coordinates are introduced.
Since the vector-coincidence equations do not involve positions, the constraints of each joint are block-triangular in the rotation and the position of the child and are solved in two explicit steps.
First, the child of a revolute joint is rotated by the minimal rotation that maps its joint axis $\mathbf{d}_b=\RotMat_b\sv_b$ onto the joint axis $\mathbf{d}_a=\RotMat_a\sv_a$ of the parent,
{\setlength{\abovedisplayskip}{3pt}
 \setlength{\belowdisplayskip}{3pt}
\be
\RotMat_b \leftarrow \RotMat_\Delta\RotMat_b \eqComma
\qquad
\RotMat_\Delta = \identity_3 + \widetilde{\mathbf{c}} + \frac{1}{1+\gamma}\,\widetilde{\mathbf{c}}^{\,2} \eqComma
\qquad
\mathbf{c} = \mathbf{d}_b\times\mathbf{d}_a \eqComma
\quad
\gamma = \mathbf{d}_b\tp\mathbf{d}_a \eqComma
\label{eq_MinimalRotation}
\ee}
whereas for orientation locks the relative rotation is reset to its constant value, $\RotMat_b\leftarrow\RotMat_a\RotMat_{a,0}\tp\RotMat_{b,0}$, which by the above is the identity up to roundoff.
$\RotMat_\Delta$ is rational in the entries of $\RotMat_a$ and $\RotMat_b$, equals $\identity_3+\bigO(|\pv|^{k+1})$ since $\mathbf{d}_b-\mathbf{d}_a=\bigO(|\pv|^{k+1})$, and is singular only for opposite axes, $\gamma=-1$.
Second, with the global joint points $\mathbf{x}_a=\yv_a+\RotMat_a\uv_a$ and $\mathbf{x}_b=\yv_b+\RotMat_b\uv_b$ evaluated with the corrected rotations, the position of the child is updated by
\be
\yv_b \leftarrow \yv_b + \Sm_p\tp\Sm_p\,(\mathbf{x}_a-\mathbf{x}_b) \eqComma
\label{eq_ConstrainedEmbedding}
\ee
which yields $\Sm_p(\mathbf{x}_a-\mathbf{x}_b)=\zero$, since the selection matrices $\Sm_p$ of \appx{sect_LinearLowerPairJoints} have orthonormal rows.
For prismatic joints, $\Sm_p\tp\Sm_p$ leaves the free translation along the slide axis unchanged.

Let $\mathcal{C}_j$ denote the map that replaces the rotation and the position of the child $b$ of joint $j$ by the results of \eq{eq_MinimalRotation} and \eq{eq_ConstrainedEmbedding} and leaves all other coordinates unchanged.
The correction map is the composition
\be
\mathcal{C} = \mathcal{C}_{j_{\nJ}}\circ\cdots\circ\mathcal{C}_{j_2}\circ\mathcal{C}_{j_1} \eqComma
\label{eq_CorrectionMap}
\ee
where $j_1,\ldots,j_{\nJ}$ enumerate the $\nJ$ joints in traversal order, so that every parent is corrected before its children.
$\mathcal{C}$ acts on the configuration $\qVec$ only and is defined wherever $\gamma\neq-1$ for all joints.
After this single traversal, all joint constraints hold exactly, and group membership is retained because rotation matrices are only multiplied by rotation matrices.
For $\mathbf{d}_b=\mathbf{d}_a$ and $\mathbf{x}_a=\mathbf{x}_b$, both updates vanish, so $\mathcal{C}$ reduces to the identity on $\mathcal{M}$ as required. %, and the magnitude of the correction is a direct measure of the truncation error.
The velocity-level joint constraints $\Bm\jnt(\qVec)\vv=\zero$ are bilinear in the state and are therefore satisfied up to $\bigO(|\pv|^{k+1})$ (\sect{sec_NonLinearConstraints}), with and without the correction.
Joints that close kinematic loops cannot be enforced by correcting a single child and would require a simultaneous projection onto $\mathcal{M}$, which is not considered here.
The reduced dynamics is unaffected, since $\mathcal{C}$ acts only when reduced coordinates are mapped to the configuration of the full model.
The offline and online parts of the evaluation are summarized in \sect{sec:implementation-lie}.

\subsection{Residual representation and constraint preservation}
\label{sec_ConstraintRepresentation}
The preservation results established in \sects{sec_LinearConstraints}{sec_NonLinearConstraints} apply only to constraints represented explicitly within the assembled residual.
Consequently, the constraint preservation properties of the resulting SSM reduced-order model are determined by the chosen residual representation and, beyond the truncation order, by the evaluation of the parameterization.
\sect{sec_LieParameterization} shows how the latter can be chosen, without altering the residual, such that the configuration-space and the joint constraints are both satisfied exactly.

Equivalent formulations of the same mechanical system may employ different residual equations to represent the same constraints.
For example, a minimal parameterization of $G$, such as exponential coordinates, satisfies the configuration-space constraints identically by construction, thereby removing them from the assembled residual, at the price of joint constraints that are no longer linear in these coordinates and hence no longer preserved exactly.
In contrast to the tangent-space parameterization of \sect{sec_LieParameterization}, which parameterizes only the SSM, such a choice alters the residual itself.
Likewise, replacing position-level joint constraints by their differentiated velocity-level counterparts removes the position-level constraints from the assembled residual.\footnote{A stabilized index-2 formulation, see e.g.\ \cite{Arnold2016}, would inherit the preservation properties for both position- and velocity-level constraints.}
In both cases, only the remaining residual equations inherit the preservation properties.

\section{Implementation}\label{sec:implementation-aspects}
Modern multibody and finite element codes repeatedly evaluate the residual of the governing equations at prescribed states as part of nonlinear solution procedures, such as time integration and nonlinear static analyses. 
The proposed residual-based SSM reduction approach leverages this existing functionality by computing SSM reduced-order models directly from evaluations of the assembled residual \eq{eq_ResCompleteShifted}. 
In contrast to conventional approaches, it neither requires a reformulation of the governing equations into a dedicated descriptor representation nor the explicit construction of multilinear maps for the nonlinear forces, cf.~\sect{sec:cohomologicalSolve}. 
The practical implementation of the residual-based SSM reduction is discussed in this section. 

The remainder of this section is organized as follows. 
\sect{sec_ImplRecAlg} introduces the recursive solution procedure and truncated Taylor arithmetic.
\sect{sec:external-residual-coupling} presents the coupling to codes providing real-valued residual evaluations only, while \sect{sec:implementation-lie} describes the rotation-preserving implementation. 
These ingredients are assembled into the complete recursion in \sect{sec:complete-recursion}, and the coupling is demonstrated with Exudyn in \appx{app:exudyn}.

\subsection{Implementation of the recursive algorithm}\label{sec_ImplRecAlg}
The computational core of the proposed residual-based SSM reduction is the recursive solution of the cohomological \eq{eq:cohomological_multiindex}. % directly from the assembled residual \eq{eq_ResCompleteShifted}. 
Beginning with the linear solution, the manifold parameterization $\Wm$ and the reduced dynamics $\Rm$ are constructed order by order. 
At polynomial order $i$, the Taylor coefficients $\Wm_i$ and $\Rm_i$ are computed recursively from the coefficients of all lower polynomial orders $j<i$.
Two computational concepts underlie the numerical realization of this recursion. 
The first is the use of complex arithmetic. 
The second is the representation of intermediate quantities by truncated multivariate Taylor series rather than ordinary numbers. 
Together, these enable a simple recursive solution.
These two concepts are addressed in the following.

\subsubsection{Complex coordinates}
\label{sec:implementation-aspects-complex}
The generalized eigenvalue problem \eq{eq_genEigProblem} determines the master subspace through its eigenvectors and the leading-order coefficients $(\Wm_1,\Rm_1)$ through \eq{eq:first_order_multiindex}. 
Whenever a selected eigenvalue has a nonzero imaginary part, its complex conjugate $\bar\lambda$ must likewise be included in the master subspace $\SpecSubSpace$.
The real and imaginary parts of the corresponding eigenvector then span $\SpecSubSpace$ as a genuine real subspace.
Each conjugate pair of eigenvalues contributes a corresponding conjugate pair to the reduced coordinates $\pv$.

Representing such a pair by two complex-conjugate scalars $p,\bar p$, rather than by two independent real coordinates, keeps the leading-order reduced dynamics in its simplest possible form for that pair
\begin{equation}
  \dot p = \lambda p + \bigO(|\pv|^2) \eqComma
\end{equation}
a single scalar equation rather than a coupled pair of real equations.
This has a direct algebraic consequence for the cohomological equations themselves. 
At every order $i$, the linear system solved for $\Wm_i$ involves the scalar $\lambda_\alpha$ (\eq{eq:weighted_eigenvalue}), built entirely from the diagonal eigenvalue matrix $\mathbf{\Lambda}$ \eq{eq_LambdaMatrix}. 
For a single conjugate pair, this matrix's action reduces, component by component, to evaluating the resonance condition
\be
\alpha_1\lambda + \alpha_2\bar\lambda - \lambda_k \approx 0
\ee
for the eigenvalue $\lambda_k$ of the coefficient being solved for, which is
an ordinary arithmetic comparison of complex numbers. 
Any other coordinate representation of the same subspace recovers this structure only after additional algebraic or trigonometric manipulations and obscures near-resonances that appear directly as small complex denominators. 
The complex representation therefore provides the simplest algebraic formulation of the recursive solution.

\subsubsection{Taylor arithmetic}
\label{sec:implementation-aspects-taylor}
At order $k$, the cohomological \eq{eq:cohomological_multiindex} requires the coefficients $\Nm_\alpha$ of \eq{eq:DefinitionNalpha} for every multi-index with $|\alpha|=k$.
Since $\CResS$ is in general a nonlinear function of its arguments, the multivariate Fa\`a di Bruno formula \cite{Constantine1996} characterizes $\Nm_\alpha$ exactly as a combinatorial sum over all ways of distributing the total degree $|\alpha|$ across the nested structure of the composition.
However, this sum need never be formed explicitly. 
Every quantity entering the residual evaluation -- the reduced coordinate $\pv$ itself, and every intermediate quantity built from it through the parameterization and the residual's own internal operations -- can instead be represented not as a single number but as its own truncated Taylor series~\cite{Berz1989,Griewank2008} in $\pv=(\pv_1,\ldots,\pv_{n_r})$.
Specifically, a quantity $u$, viewed as a function of $\pv$, is represented by the finite coefficient set
{\setlength{\abovedisplayskip}{3pt}
 \setlength{\belowdisplayskip}{3pt}
\begin{equation}
  u \leftrightarrow \{u_\alpha\}_{|\alpha|\le N} \eqComma \quad
  u(\pv) = \sum_{|\alpha|\le N} u_\alpha\, \pv^\alpha
  + \bigO\!\left(|\pv|^{N+1}\right) ,
\end{equation}}
where $\alpha=(\alpha_1,\ldots,\alpha_{n_r})$ is a multi-index and $N$ denotes the polynomial order currently being resolved by the recursion.
The entries $u_\alpha$ are the Taylor coefficients of $u$.

The coefficient-set form is used from here on, for two reasons. 
First, it is not restricted to a single degree. 
The set $\{u_\alpha\}_{|\alpha|\le N}$ collects every degree up to order $N$ in one object, unlike $\Wm_i$, which is homogeneous of degree $i$ alone. 
A general quantity built from the residual's own operations mixes several degrees at once, so this is the form it actually needs. 
Second, the elementary arithmetic operations from which the residual is built are simple to define directly on this form. %, coefficient by coefficient. 
Let $u$ and $v$ be two such sets, representing two quantities. 
Their sum and product are again coefficient sets of the same form. 
Each is defined to be the exact truncated Taylor series of the corresponding sum or product. %, not just some combination of the individual coefficients. 
Addition acts coefficientwise: $(u+v)_\alpha=u_\alpha+v_\alpha$.
Multiplication follows the truncated Cauchy product
\begin{equation}
  (uv)_\alpha = \sum_{\beta\le\alpha} u_\beta\,v_{\alpha-\beta} \eqComma
  \qquad |\alpha|\le  N \eqDot
  \label{eq:cauchy-product-impl}
\end{equation}
Integer powers and other elementary nonlinear operations appearing in the residual admit an analogous, closed-form recursion for their own Taylor coefficients \citep{Berz1989,Griewank2008}.

As an illustration, tracking a single scalar coordinate $p$ to order $N=2$ around $p=p_0$ amounts to the coefficient triple $(p_0,1,0)$.
Squaring it via the Cauchy product gives $(p_0^2,2p_0,1)$, and adding $(3p_0,3,0)$ yields $(p_0^2+3p_0,\,2p_0+3,\,1)$ -- the value, first derivative and factorial-scaled second derivative of $p^2+3p$ at $p=p_0$, obtained without ever differentiating that expression symbolically.

A single residual evaluation with series-valued arguments directly computes every Taylor coefficient required for the current-order cohomological equation.
Consequently, every component of $\Nm_\alpha$ (\eq{eq:DefinitionNalpha}) is obtained exactly, up to working numerical precision, without approximation.

In contrast, estimating the same coefficients via repeated finite differencing of the residual carries a truncation error at every order that does not remain confined to that order.
As each order's own computation depends on the already-determined coefficients of every lower order, any imprecision introduced at one order propagates into, and can be amplified by, every subsequent order of the recursion.
The Taylor-series representation avoids this entirely, by construction, i.e., every coefficient at every order is exact, independent of how many recursive orders precede it.

\subsubsection{Automatic differentiation}
\label{sec:implementation-aspects-ad}
%%%%%%%%%%%%%%%%%%%%%%%%%%%%%%%%%%%%%%%%%%%%%%%%%%%%%%%%%%%%%%%%%%%%%%%%%%%%%%%%
%
The recursive solution requires higher-order derivatives of the governing residual. 
In general, derivatives can be obtained symbolically, numerically, or by automatic differentiation~\cite{Griewank2008,Baydin2018,Neidinger2010}. 
Symbolic differentiation requires an explicit analytical expression for the residual and rapidly becomes impractical for nonlinear compositions of increasing order. 
Moreover, it is unavailable when the residual is provided only by single or double precision numbers. 
Numerical differentiation avoids these limitations but introduces truncation errors that propagate through the recursive solution.

The proposed framework instead employs higher-order forward automatic differentiation. 
In particular, truncated Taylor arithmetic constitutes a higher-order forward automatic differentiation method \citep{Berz1989,Griewank2008}, as the elementary arithmetic operations are redefined to propagate Taylor coefficients according to the chain and product rules.
Consequently, differentiation is not performed as a separate computational step but occurs automatically during the residual evaluation, without symbolic expressions or finite-difference approximations.
Once the elementary arithmetic has been redefined for Taylor coefficients (cf. \sect{sec:implementation-aspects-taylor}), the recursive solution requires only repeated evaluations of the governing residual. 
The required derivatives are obtained automatically during these evaluations.

\subsubsection{Seeding and order-wise extraction}
\label{sec:implementation-aspects-algorithm}

The preceding discussion can be assembled into an explicit recursive procedure. 
Two steps need further explanation before the procedure is stated.
That is how the reduced coordinate $\pv$ enters the recursion as a Taylor-typed quantity, and how the known term $\Nm_\alpha$ of \eq{eq:DefinitionNalpha} is obtained from a single evaluation of $\CResS$.
\paragraph{Seeding.}
The recursion needs $\pv$ itself represented as a truncated Taylor series in $\pv$, before any computation begins, carrying coefficients up to the target order $\ordMax$ of the SSM approximation. 
Each component $p_j$ is seeded by the coefficient set
\begin{equation}
  (p_j)_\alpha =
  \begin{cases}
    1 & \alpha = e_j \eqComma \\
    0 & \text{otherwise,}
  \end{cases}
  \qquad |\alpha|\le \ordMax \eqComma
\end{equation}
where $p_j$, viewed as a function of $\pv$, is $p_j$ itself. 
Its coefficient in its own direction, $\alpha=e_j$, is $1$. 
Every other coefficient, including at the origin, is $0$.

\paragraph{Order-$k$ extraction.}
With $\pv$ seeded this way, evaluating any polynomial expression in $\pv$ again yields a truncated Taylor series in $\pv$.
This applies in particular to the truncated ansatz $\Wm^{<k}(\pv)$ and its velocity field $(D_{\pv}\Wm^{<k}(\pv))\Rm^{<k}(\pv)$, both built from the already-known coefficients by the coefficientwise and Cauchy-product rules above.
A single evaluation of $\CResS$ at this Taylor-typed argument then delivers $\Nm_\alpha$ for every multi-index of order $k$ at once, by \eq{eq:DefinitionNalpha}.

This last point is worth emphasising, because it is not merely a convenience.
The evaluation depends on $\alpha$ only through the truncation order $k$, so evaluating once per multi-index rather than once per order repeats the entire nonlinear composition unnecessarily.
The number of multi-indices at order $k$ grows as $\binom{k+\dimIVM-1}{k}$; for $\dimIVM=4$ and $\ordMax=13$ the batched evaluation performs twelve residual evaluations where the per-multi-index variant performs some $2375$, all but twelve of them redundant.

Neither construction avoids the combinatorial structure of $\Fm\circ\Wm$. 
They differ only in where this structure is resolved.
In \cite{jainhaller2021how} the nonlinearity is supplied as explicit tensor coefficients and the composition is expanded for that particular form at each order.
Here the Cauchy product and composition rules of \sect{sec:implementation-aspects-taylor} resolve it once, generically, independently of the form $\CResS$ takes.

\subsubsection{Order-wise extraction by contour integration}
\label{sec:implementation-aspects-contour}
Taylor-typed evaluation requires every operation of the residual to accept series-valued arguments, and its memory grows with the length of the composition.
The coefficients $\Nm_\alpha$ are therefore alternatively extracted by contour integration in the reduced coordinates, which requires only residual evaluations at complex arguments.
Along a line $\pv=t\,\mathbf{w}$ with a real unit direction $\mathbf{w}\in\RealNumbers^{\dimIVM}$, the residual evaluated at the truncated ansatz, $\boldsymbol{\gamma}(t):=\CResS\bigl(\Wm^{<k}(t\mathbf{w}),(D_{\pv}\Wm^{<k})\Rm^{<k}\bigr)$, is sampled at $N_c$ points on a circle of radius $r_c$, and its $k$-th coefficient follows from the trapezoidal discretization of Cauchy's integral formula, i.e., a discrete Fourier transform \cite{LynessMoler1967,Fornberg1981},
{\setlength{\abovedisplayskip}{3pt}
 \setlength{\belowdisplayskip}{3pt}
\begin{equation}
\frac{1}{N_c\,r_c^{\,k}}\sum_{j=0}^{N_c-1}\boldsymbol{\gamma}\bigl(r_c\,e^{2\pi \mathrm{i}j/N_c}\bigr)\,e^{-2\pi \mathrm{i}jk/N_c}
=\sum_{|\alpha|=k}\Nm_\alpha\,\mathbf{w}^\alpha
+\sum_{m\ge1}r_c^{\,mN_c}\,\boldsymbol{\gamma}_{k+mN_c}\eqComma
\label{eq_contourExtraction}
\end{equation}}
where $\boldsymbol{\gamma}_{k}$ denotes the $k$-th Taylor coefficient of $\boldsymbol{\gamma}(t)$.
Repeating \eq{eq_contourExtraction} for at least $\binom{k+\dimIVM-1}{k}$ directions $\mathbf{w}$, chosen such that the matrix of monomials $\mathbf{w}^\alpha$ has full rank, yields all $\Nm_\alpha$ with $|\alpha|=k$ by a least-squares solve.
For a polynomial residual, $\boldsymbol{\gamma}(t)$ is a polynomial, and the aliasing sum in \eq{eq_contourExtraction} vanishes once $N_c$ exceeds its degree minus $k$.
Otherwise it decays with the coefficients of $\boldsymbol{\gamma}$.
Round-off, in contrast, is amplified by $r_c^{-k}$, so that the radius is chosen of the order of the amplitudes of interest rather than small as a finite-difference step \cite{Bornemann2011}.

%%%%%%%%%%%%%%%%%%%%%%%%%%%%%%%%%%%%%%%%%%%%%%%%%%%%%%%%%%%%%%%%%%%%%%%%%%%%%%%%
\subsection{Coupling real-valued residual evaluations}
\label{sec:external-residual-coupling}
%%%%%%%%%%%%%%%%%%%%%%%%%%%%%%%%%%%%%%%%%%%%%%%%%%%%%%%%%%%%%%%%%%%%%%%%%%%%%%%%
%
The preceding sections implicitly assumed $\CRes$ is built from generic arithmetic primitives. 
These primitives dispatch correctly regardless of their arguments' own type. 
Evaluating $\CRes$ directly at complex- or Taylor-typed arguments then yields the derivative or Taylor-coefficient information required by the recursion. 
However, when the residual is provided by an established multibody or finite element code, this is typically not the case. 
This is because the residual assembly is performed by compiled software operating on real double-precision data and therefore does not support the complex- or Taylor-typed arguments required above.
The proposed coupling assumes only that the residual can be evaluated at arbitrary states in a neighborhood of the equilibrium. No access to symbolic expressions, derivative information, or internal data structures is required in general. 
Starting from these real-valued residual evaluations, all quantities required by the reduction procedure described in the previous sections are subsequently recovered.

Two complementary strategies are employed. First-order derivative information required for the equilibrium computation and linearization is recovered either from analytically assembled system matrices or, alternatively, by finite differences. 
The nonlinear dependence required by the recursive solution of the invariance equations is recovered by constructing an explicit polynomial representation of the residual from real-valued residual evaluations.

\subsubsection{Linearization}
Where available, analytically assembled system matrices provided by the underlying multibody or finite element code are used directly to construct the operators $\SSA$ and $\SSB$. 
These include the mass matrix $\Mb(\qVec)$, tangent stiffness and damping matrices, and the constraint Jacobians $\Bm_{\bullet}(\qVec)$ for example. % of \eq{eq_Bq}. 
Otherwise, the required Jacobians are approximated using real finite differences, introducing the corresponding truncation error.
Since this approximation is required only during the equilibrium computation and linearization stage, its influence is confined to these computations.

\subsubsection{Extracting the nonlinear coefficients}\label{sect_ExtractNonlinearCoeff}
The recursive solution also requires the coefficients of degree $\ge 2$ of the residual.
An existing multibody code typically evaluates the residual at arbitrary states, but only in real arithmetic, so the Taylor arithmetic of the recursion cannot be propagated through it.
We therefore reconstruct the residual once near the equilibrium from real-valued evaluations as an explicit polynomial, which is then evaluated with the Taylor arithmetic of the recursion.
The coupling to a specific code is described in \sect{sec:coupling-existing-code}.

The residual is separated into terms that are linear in $\dot{\vv}$ and $\boldsymbol{\lambda}$, i.e.\ $\Mm(\qVec)\dot{\vv}$ and $\Bm_\mathrm{jnt}\tp(\qVec)\boldsymbol{\lambda}$, a configuration-dependent force, and a velocity-dependent force that is exactly quadratic in $\vv$. 
The remaining configuration-dependent quantities, including the position-level constraints $\boldsymbol{\phi}_\mathrm{jnt}(\qVec)$, are fitted as polynomials in $\qVec-\qVec_\mathrm{eq}$ on a box of radius $\varrho$ around the equilibrium. 
Unit-vector probing first detects the variables on which each output depends, so that outputs with the same dependencies are fitted jointly.
We use the Chebyshev basis with Chebyshev--Gauss--Lobatto nodes~\cite{Trefethen2019}. 
For few active variables, the tensor-product interpolant is exact for polynomials of bounded degree per variable and requires one evaluation per grid point. For more active variables, we switch to total-degree Chebyshev least squares at points drawn from the Chebyshev measure~\cite{Cohen2013,Chkifa2015}. 
The basis is identical in both cases, only the sampling and truncation change. 
All fits are verified at held-out points.

Non-intrusive coefficient extraction has also been proposed for finite-element models, namely the stiffness evaluation procedure and enforced-displacement methods~\cite{MuravyovRizzi2003,Perez2014,Pereira2025}. 
These fit the internal polynomial force of the model, typically up to cubic order, from static solutions with prescribed displacements. 
Since we require all terms of the residual (mass, constraint, and velocity-dependent terms) up to the SSM order, we sample the residual itself instead.

\subsection{Implementation of the rotation-preserving reduction}
\label{sec:implementation-lie}
%%%%%%%%%%%%%%%%%%%%%%%%%%%%%%%%%%%%%%%%%%%%%%%%%%%%%%%%%%%%%%%%%%%%%%%%%%%%%%%%
%
The rotation-preserving parameterization of \sect{sec_LieParameterization} affects two places in the recursion and the evaluation of the reduced model, and is otherwise transparent to the recursion.
\begin{itemize}

\item \textsl{Projection of the cohomological solve.}
The matrices $\Sm$ and $\Vm$ and the reduced system \eq{eq_SLTSystem} are
assembled once after the fixed point is known and reused at every order. Both
projections are required to obtain a square kinematic block
(\sect{sec_LieTangentParam}).

\item \textsl{Recursive assembly of $\boldsymbol{\Gamma}_\alpha$.}
By \eq{eq_LieReconstructionOrder}, $\boldsymbol{\Gamma}_\alpha$ is the order-$\alpha$ coefficient of the exponential series built from intrinsic coefficients of order strictly below $|\alpha|$.
Two properties make its assembly inexpensive.
First, the series terminates, so that neither a convergence criterion nor a transcendental evaluation is involved.
Second, as for $\Nm_\alpha$, the assembly depends on $\alpha$ only through the order, so that the convolution is evaluated once per order and indexed rather than recomputed for every multi-index.
The intrinsic coefficients $\wv_\alpha$ obtained from each solve must be recorded for both the multi-index and its conjugate partner, since coefficients related by conjugate symmetry are filled without passing through the solve.

\item \textsl{Evaluation of the reduced model.}
The configuration is evaluated along the chain \eq{eq_EmbeddingChain}: the intrinsic polynomial $\wv(\pv)$, the group exponential, and the correction $\mathcal{C}$ of \sect{sec_LieConstrainedEmbedding}.
The kinematic tree, the traversal order, and the constant joint data $\sv_a$, $\sv_b$, $\uv_a$, $\uv_b$, and $\Sm_p\tp\Sm_p$ are determined offline, once after the recursion.
Online, the joints are processed in the stored order by the explicit updates \eq{eq_MinimalRotation} and \eq{eq_ConstrainedEmbedding}, at the cost of a few $3\times3$ products per joint, so that no solve, factorization, or residual evaluation occurs at evaluation time.
The group exponential is the only transcendental evaluation and can be computed efficiently using closed-form formulas, see, e.g., \cite{Sonneville2014}.
\end{itemize}

%%%%%%%%%%%%%%%%%%%%%%%%%%%%%%%%%%%%%%%%%%%%%%%%%%%%%%%%%%%%%%%%%%%%%%%%%%%%%%%%
\subsection{The complete recursion}
\label{sec:complete-recursion}
%%%%%%%%%%%%%%%%%%%%%%%%%%%%%%%%%%%%%%%%%%%%%%%%%%%%%%%%%%%%%%%%%%%%%%%%%%%%%%%%
%
Algorithm~\ref{alg:recursion} assembles the preceding subsections into the complete recursion, which is carried out offline.
Steps marked $\triangleright$ belong to the rotation-preserving variant of \sect{sec_LieParameterization}, described in \sect{sec:implementation-lie}.
Omitting them recovers the ambient formulation exactly.
Both the residual evaluation and the assembly of $\boldsymbol{\Gamma}_\alpha$ are placed outside the loop over multi-indices, since each depends on $\alpha$ only through the order.
Algorithm~\ref{alg:evaluation} summarizes the online evaluation of the configuration of the rotation-preserving reduced model at a given amplitude $\pv$, following the chain \eq{eq_EmbeddingChain}.
\begin{algorithm}%[htb]
\caption{Offline recursive computation of the manifold and reduced dynamics
coefficients. Steps marked $\triangleright$ are required only for the
rotation-preserving variant of \sect{sec_LieParameterization} and are omitted
in the ambient formulation.}
\label{alg:recursion}
\KwIn{Leading-order coefficients $\Wm_1$, $\Rm_1=\mathbf{\Lambda}$
(\eq{eq:first_order_multiindex}); truncation order $k$}
\KwOut{Coefficients $\Wm_\alpha$ and $\Rm_\alpha$ for $|\alpha|\le k$; $\triangleright$ intrinsic coefficients $\wv_\alpha$}
Seed $\pv$: for $j=1,\ldots,n_r$, set $(p_j)_\alpha=1$ if $\alpha=e_j$,
else $0$, for $|\alpha|\le k$\;
$\triangleright$ Assemble $\Sm,\Vm$ from $\kinRecMat_0$ at the fixed point (\eq{eq_TangentBasis})\;
$\triangleright$ Seed the intrinsic coefficients at first order:
$\wv_{e_j}\leftarrow(\kinRecMat_0\tp\kinRecMat_0)^{-1}\kinRecMat_0\tp\,\qVec_{e_j}$ for
$j=1,\ldots,n_r$\;
\For{$\ell=2,\ldots,k$}{
    Build $\Wm^{<\ell}(\pv)$ and $\Rm^{<\ell}(\pv)$ from the known
    coefficients, evaluated at the seeded $\pv$\;
    Evaluate
    $\CResS\!\left(\Wm^{<\ell}(\pv),(D_{\pv}\Wm^{<\ell}(\pv))\Rm^{<\ell}(\pv)\right)$,
    \emph{one} Taylor-typed evaluation for the whole order\;
    $\triangleright$ Assemble $\{\bar{\boldsymbol{\Gamma}}_\alpha:|\alpha|=\ell\}$ from
    the intrinsic coefficients of order $<\ell$
    (\eq{eq_LieReconstructionOrder}), \emph{one} convolution for the whole
    order\;
    \For{each $\alpha$ with $|\alpha|=\ell$}{
        Extract $\Nm_\alpha \leftarrow$ the coefficient of $\pv^\alpha$
        (\eq{eq:DefinitionNalpha})\;
        \eIf{ambient}{
            Solve
            $(\SSA+\lambda_\alpha\SSB)\Wm_\alpha+\SSB\Wm_1\Rm_\alpha=-\Nm_\alpha$
            for $(\Wm_\alpha,\Rm_\alpha)$ (\eq{eq:cohomological_multiindex})\;
        }{
            $\triangleright$ Solve the reduced system
            \eq{eq_SLTSystem} for $(\WmRed_\alpha,\Rm_\alpha)$ and set
            $\Wm_\alpha=\Vm\WmRed_\alpha+\bar{\boldsymbol{\Gamma}}_\alpha$\;
            $\triangleright$ Record the intrinsic coefficients $\wv_\alpha$, i.e., the rotation part of $\WmRed_\alpha$, and, if $\alpha$ has a
            distinct conjugate partner $\bar\alpha$, also
            $\wv_{\bar\alpha}=\overline{\wv_\alpha}$\;
        }
    }
}
\end{algorithm}
\begin{algorithm}%[htb]
\caption{Online evaluation of the configuration of the rotation-preserving
reduced model at an amplitude $\pv$ (\sect{sec_LieConstrainedEmbedding}).}
\label{alg:evaluation}
\KwIn{Coefficients from Algorithm~\ref{alg:recursion}; joints
$j_1,\ldots,j_{\nJ}$ in traversal order with parent $a$, child $b$, and
constant joint data $\sv_a,\sv_b,\uv_a,\uv_b,\Sm_p\tp\Sm_p$ and, for
orientation locks, the relative rotation $\RotMat_{a,0}\tp\RotMat_{b,0}$, all
determined offline; amplitude $\pv$}
\KwOut{Configuration $\qVec\in\mathcal{M}$}
Evaluate the polynomials $\wv(\pv)$ and $\Wm(\pv)$\;
Reconstruct the group block $\qVec=\vecOp\bigl(q_0\exp_G(\widetilde{\wv}(\pv))\bigr)$ (\eq{eq_LieReconstruction})\;
\For{$i=1,\ldots,\nJ$}{
    \uIf{$j_i$ is a revolute joint}{
        $\RotMat_b\leftarrow\RotMat_\Delta\RotMat_b$ (\eq{eq_MinimalRotation})\;
    }
    \ElseIf{$j_i$ locks the relative orientation}{
        $\RotMat_b\leftarrow\RotMat_a\RotMat_{a,0}\tp\RotMat_{b,0}$\;
    }
    $\mathbf{x}_a\leftarrow\yv_a+\RotMat_a\uv_a$,\quad
    $\mathbf{x}_b\leftarrow\yv_b+\RotMat_b\uv_b$\;
    $\yv_b\leftarrow\yv_b+\Sm_p\tp\Sm_p(\mathbf{x}_a-\mathbf{x}_b)$ (\eq{eq_ConstrainedEmbedding})\;
}
\end{algorithm}

%%%%%%%%%%%%%%%%%%%%%%%%%%%%%%%%%%%%%%%%%%%%%%%%%%%%%%%%%%%%%%%%%%%%%%%%%%%%%%%%
\subsection{Coupling to an existing multibody code}
\label{sec:coupling-existing-code}
%%%%%%%%%%%%%%%%%%%%%%%%%%%%%%%%%%%%%%%%%%%%%%%%%%%%%%%%%%%%%%%%%%%%%%%%%%%%%%%%
%
Nothing in the preceding construction presumes a particular multibody implementation.
The recursion requires from a host code only its assembled residual \eq{eq_CompleteResidual}, or the individual blocks from which it is formed, exposed as a callable object, together with the partition of the state \eq{eq:ComputationalState}.
The rotation-preserving variant additionally requires the indices and the group type of every factor of $G$, and the constrained evaluation of \sect{sec_LieConstrainedEmbedding} the joint data $\sv_a$, $\sv_b$, $\uv_a$, $\uv_b$, and $\Sm_p$ with which the joints are defined in the multibody model.
Since these data are index sets and geometric model parameters, only the residual depends on the arithmetic of the host.
Where the host evaluates its residual only in real arithmetic, so that no complex or Taylor-typed arguments can be passed through it, the coupling strategy of \sect{sec:external-residual-coupling} applies.
\appx{app:exudyn} demonstrates this for the multibody code Exudyn \cite{Gerstmayr2024}, whose residual blocks are exposed in exactly this form.
In this example, the rigid-body rotations are parameterized by Euler parameters, whose unit-norm constraints are contained in the algebraic equations of Exudyn and are not distinguished from the joint constraints, so that the reduced model is computed and evaluated in the ambient formulation.

\section{Numerical examples}\label{sec:numerical-examples}
This section presents several examples demonstrating the proposed SSM-based model reduction for constrained multibody systems formulated as DAEs on manifolds with Lie group structure.
Both free and forced vibrations are predicted by the SSM-based reduced-order models and validated against corresponding full-order reference solutions.
The proposed reduction approach and all numerical examples were implemented in 64-bit Python 3.13.
All simulations have been performed under Windows 10 on an Intel Core i7-6600U 2.60 GHz processor.

\subsection{Pendulum models with and without internal resonances}
\label{sec_PendulumExamples}
\begin{figure}%[h]
\centering
\begin{subfigure}{0.49\textwidth}
    \includegraphics[width=0.7\textwidth]{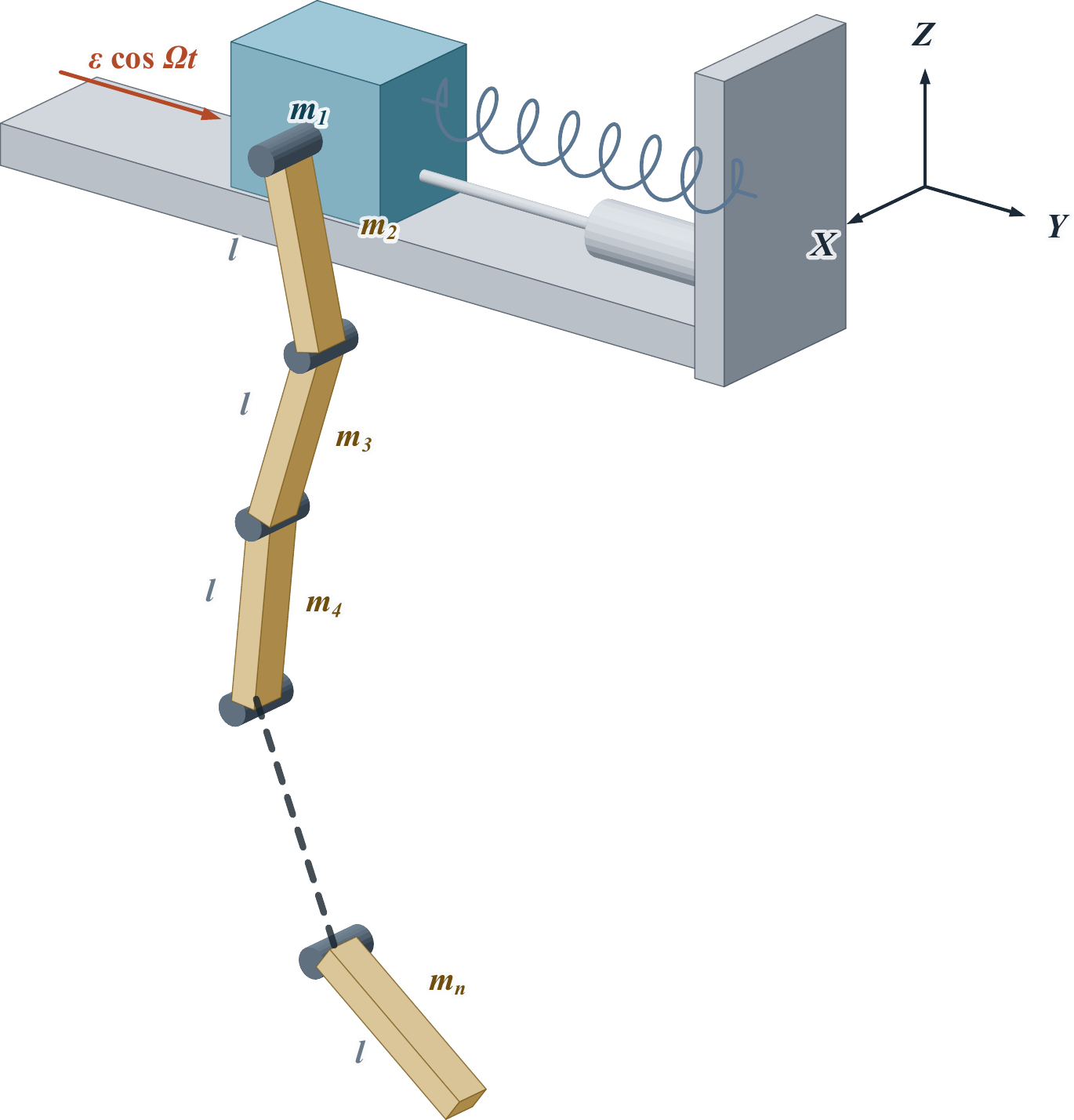}
    \caption{Rigid pendulum slider systems}
    \label{fig_pendulum_chain}
\end{subfigure}
\hfill
\begin{subfigure}{0.49\textwidth}
    \includegraphics[width=0.55\textwidth]{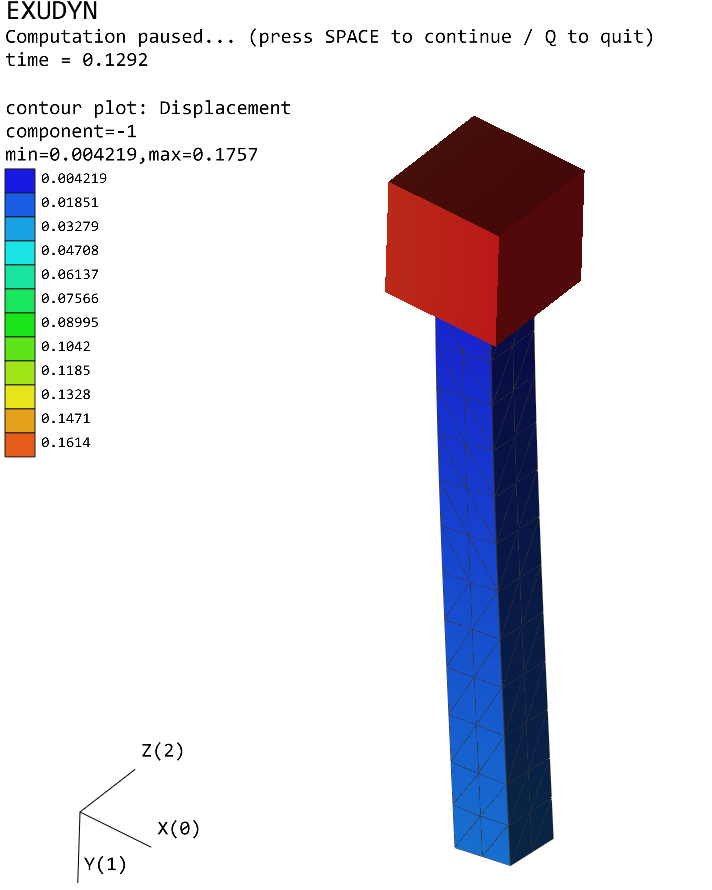}
    \caption{Flexible pendulum slider system modeled in Exudyn}
    \label{fig_flexiblePendulumSliderExu}
\end{subfigure}
\caption{Considered pendulum systems.}
\label{fig_pendulum_systems}
\end{figure}

In this section, we demonstrate the proposed SSM reduction using three different pendulum systems as example.
The first two pendulum systems considered have been proposed in \cite{li2023model} for evaluating SSM reduced-order models of constrained mechanical systems and are therefore considered here as well and are schematically illustrated in \fig{fig_pendulum_chain}.
All considered pendulum systems consist of a slider modeled as a rigid body that is attached to ground via a prismatic joint and a linear spring-damper with stiffness $k_1$ and velocity proportional damping $c_1$.
The first pendulum system attaches a single pendulum body modeled as rigid body to the slider via a revolute joint.
A linear torsional spring-damper element affects the relative rotation of pendulum with respect to the slider.
The second pendulum example consists of a chain of 41 rigid pendulum bodies attached to the slider, with torsional spring-damper elements acting on the relative rotations between adjacent bodies.
In the third pendulum system considered, we replace the rigid pendulum body of the first pendulum example with a flexible body modeled with the nodal-based FFRF with modal reduction \cite{Zwoelfer2019}.
While the two rigid body pendulum systems are modeled in python using the proposed implementation, the flexible pendulum slider system is modeled in Exudyn and is used to demonstrate the coupling to an existing multibody code.
The multibody model of the flexible pendulum slider system is schematically shown in \fig{fig_flexiblePendulumSliderExu}.

\subsubsection{Rigid slider-pendulum: 1:3 internal resonance}\label{sec_PendulumSliderRigid}
The first pendulum system that we consider was used in \cite[Sect. 5.2.2]{li2023model} consists of two rigid bodies (slider, pendulum) and was used to study its free and forced vibrations under SSM reduction.
The parameters of this two-body system were chosen to produce a 1:3 internal resonance between the first two modes. Detailed parameter values are provided in \cite{li2023model}.

To investigate the behavior of SSM reduced-order models when using different Lie group formulations for modeling the multibody system, we model the slider pendulum system using both direct and semi-direct products of $\RealNumbers^3$ and $SO(3)$, respectively $\RealNumbers^3$ and $\SThree$.
Here, $\RealNumbers^3\times\SOThree$ and $\SEThree$ denote the configuration spaces of the full model. %, differing only in the frame in which the translational velocity is expressed. 
In both cases, the rotation-preserving parameterization of \sect{sec_LieParameterization} acts on the rotational block.
Following \cite{li2023model}, we select the spectral subspace spanned by the first two modes as the master subspace for SSM computation.

\begin{figure}%[!h]
\centering
\begin{subfigure}{0.49\textwidth}
    \includegraphics[width=1\textwidth]{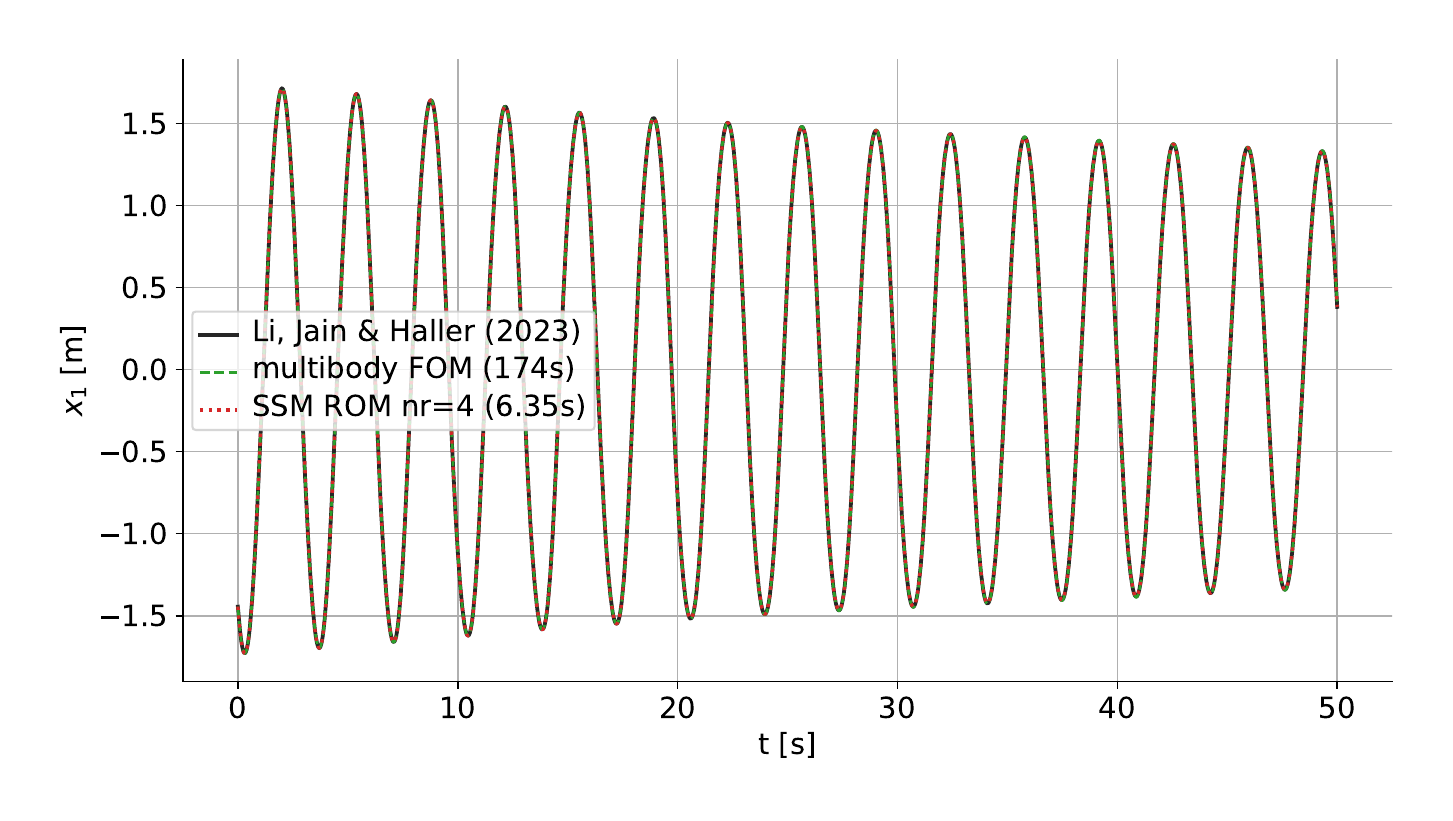}
    \vspace{-0.75cm}
    \caption{Slider position}
    \label{fig_SPR_Autonomous_x1_vs_Reference}
\end{subfigure}
\hfill
\begin{subfigure}{0.49\textwidth}
    \includegraphics[width=1\textwidth]{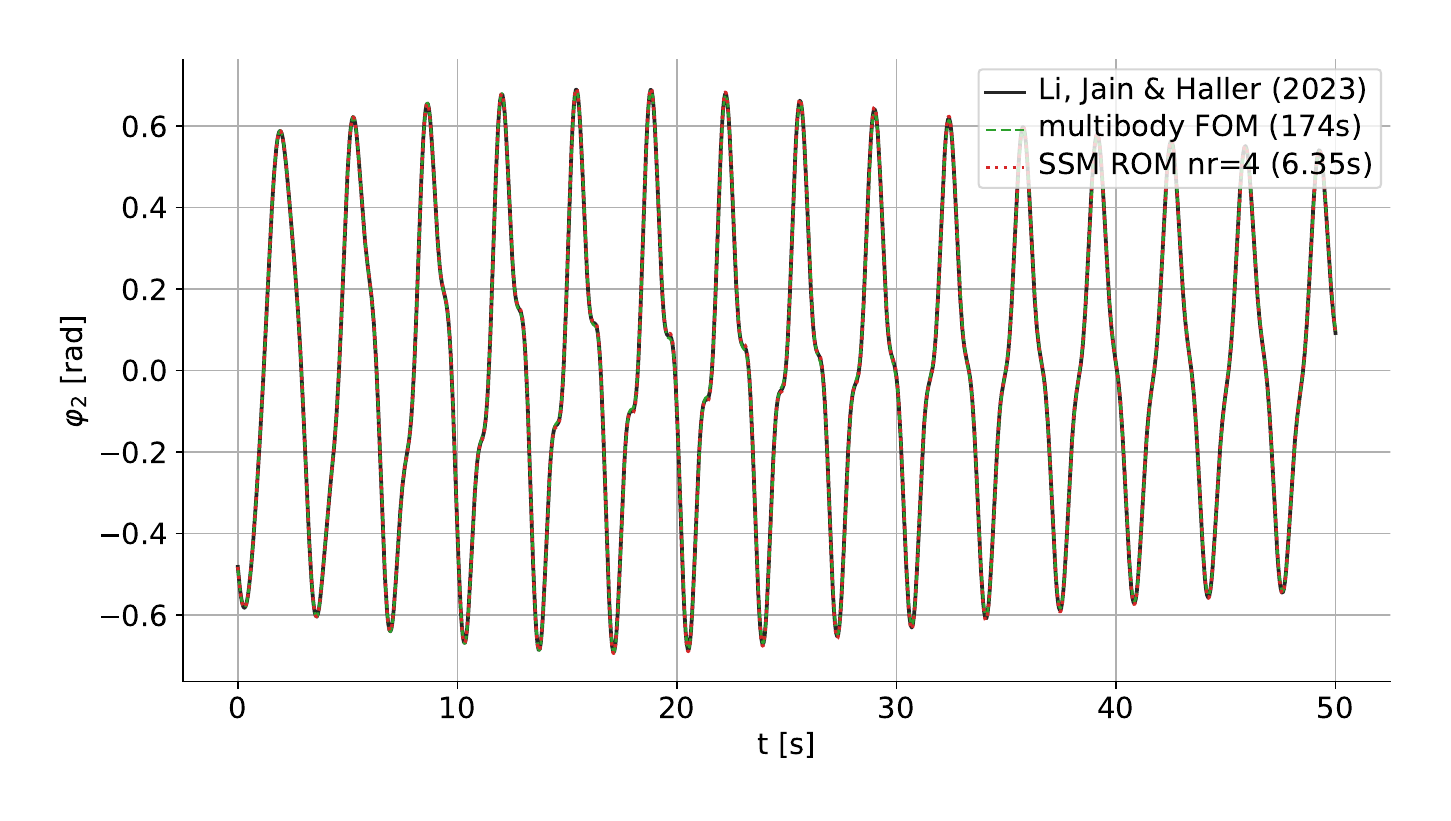}
    \vspace{-0.75cm}
    \caption{Pendulum angle}
    \label{fig_SPR_Autonomous_phi2_vs_Reference}
\end{subfigure}
\vspace{-0.1cm}
\caption{Time histories for the displacements $(x_1, \varphi_2$) of the pendulum slider system with 1:3 internal resonance.}
\label{fig_SPR_SliderPosPendulumAngle}
\end{figure}

First, we compare the trajectories from the simulations of a $\bigO(13)$-SSM ROM when elements of the Lie group $\RealNumbers^3 \times SO(3)$ are used as configuration space variables and compare with the trajectories of the full system.
Following \cite{li2023model}, we consider the initial condition $\pv_0=(3.5e^i, \, 3.5e^{-i}, \, 0, \, 0)$ for the trajectories of the autonomous system.
The time histories computed with the $\RealNumbers^3 \times SO(3)$-based $\bigO(13)$-SSM ROM and the full order model, which is based on a minimum coordinate formulation and given in \cite{li2023model}, are shown in \fig{fig_SPR_SliderPosPendulumAngle}.

As can be seen from \fig{fig_SPR_SliderPosPendulumAngle}, the ROM trajectories match the FOM trajectories well and are visually almost indistinguishable. 
The maximum ROM error in the slider position (\fig{fig_SPR_Autonomous_x1_vs_Reference}) is approximately $3\cdot10^{-3}$ m and the maximum error in the pendulum angle (\fig{fig_SPR_Autonomous_phi2_vs_Reference}) is approximately $8\cdot10^{-3}$ radians.

The slider response in \fig{fig_SPR_Autonomous_x1_vs_Reference} is dominated by the first mode, with the slider vibrating near its first natural frequency, whereas the pendulum response in \fig{fig_SPR_Autonomous_phi2_vs_Reference} is dominated by the second mode, with the pendulum oscillating at approximately three times that frequency.
The amplitude of $\varphi_2$ in \fig{fig_SPR_Autonomous_phi2_vs_Reference} first increases and then decreases, which indicates an exchange of energy between the slider vibration and the pendulum oscillation.

\begin{figure}%[h]
\centering
\begin{subfigure}{0.49\textwidth}
    \includegraphics[width=\textwidth]{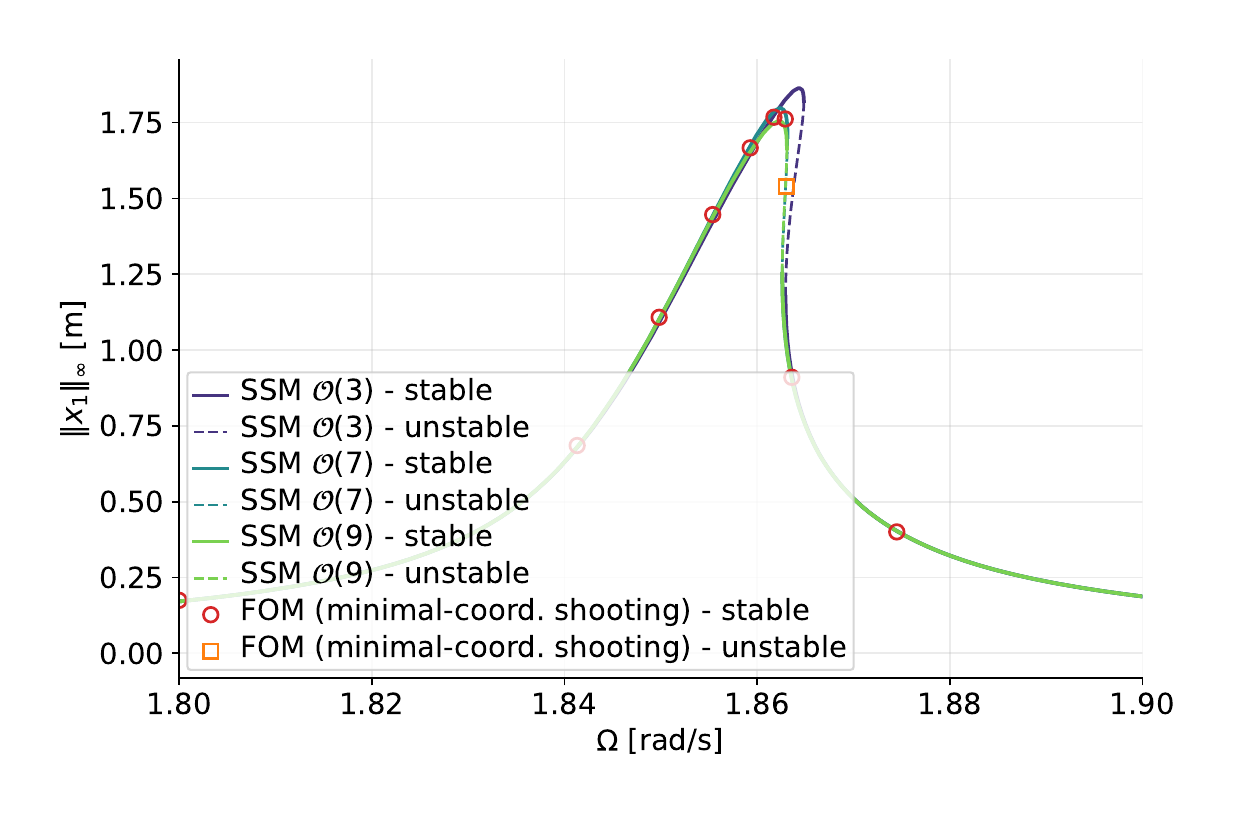}
    \vspace{-0.75cm}
    \caption{FRC: Slider position}
    \label{fig_SPR_FRC_x1}
\end{subfigure}
\hfill
\begin{subfigure}{0.49\textwidth}
    \includegraphics[width=\textwidth]{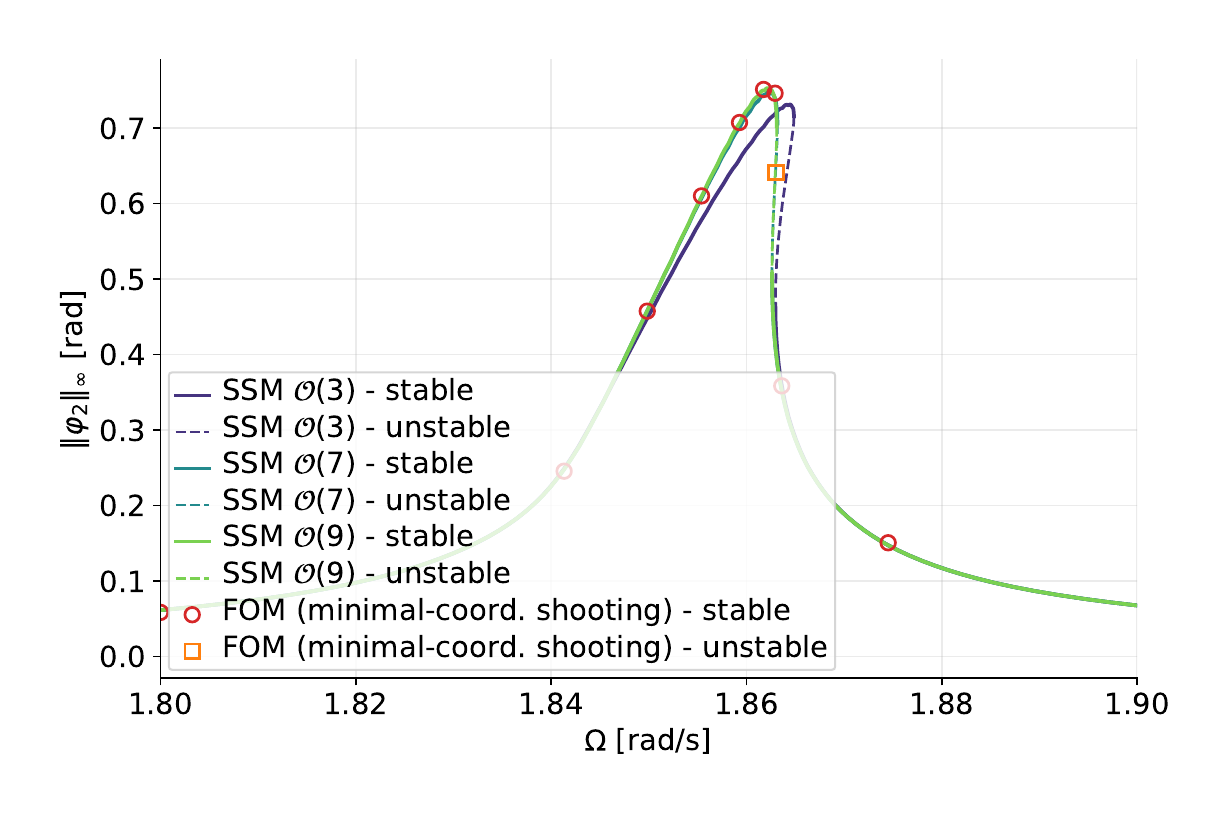}
    \vspace{-0.75cm}
    \caption{FRC: Pendulum angle}
    \label{fig_SPR_FRC_phi2}
\end{subfigure}
\caption{Forced response curves in vibration amplitude of the slider and the pendulum obtained by SSM reduction at various approximation orders.}
\label{fig_SPR_FRC}
\end{figure}
Adding a periodic forcing $\epsilon \cos(\Omega t)$ to the slider, we
now compute the FRC of the system with $\epsilon = 0.08$
and $\Omega \approx \textrm{Im}(\lambda_1)$. 
The FRC is obtained directly by analyzing the reduced dynamics of the four-
dimensional SSM in a normal-form parameterization style, where the periodic orbit is simply given by the solution to a fixed point problem, cf. \cite{li2022nonlinear}.
The FRCs for the horizontal vibration of the slider and the oscillation of the pendulum are shown in \fig{fig_SPR_FRC}.

We observe that the FRCs obtained from the SSM-based reduction converge toward the reference solution at $\bigO(7)$, consistent with the observations reported for the minimum-coordinate-based SSM ROM in \cite{li2022nonlinear}.

\begin{figure}%[h]
\centering
\begin{subfigure}{0.49\textwidth}
    \includegraphics[width=1.05\textwidth]{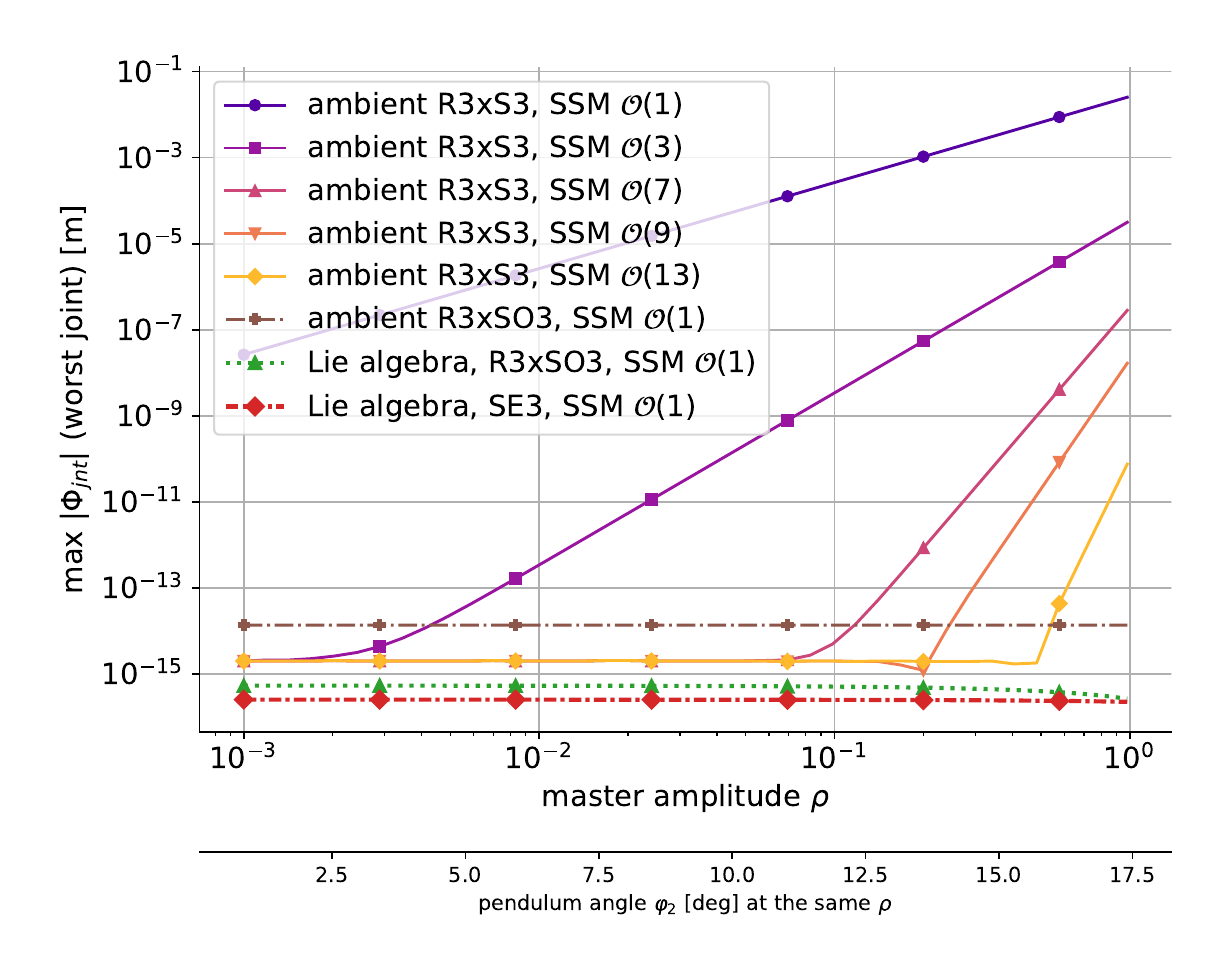}
    \vspace{-0.75cm}
    \caption{Error in the joint constraints}
    \label{fig_SPR_JointConstraint_OrderAndLie}
\end{subfigure}
\hfill
\begin{subfigure}{0.49\textwidth}
    \includegraphics[width=1.05\textwidth]{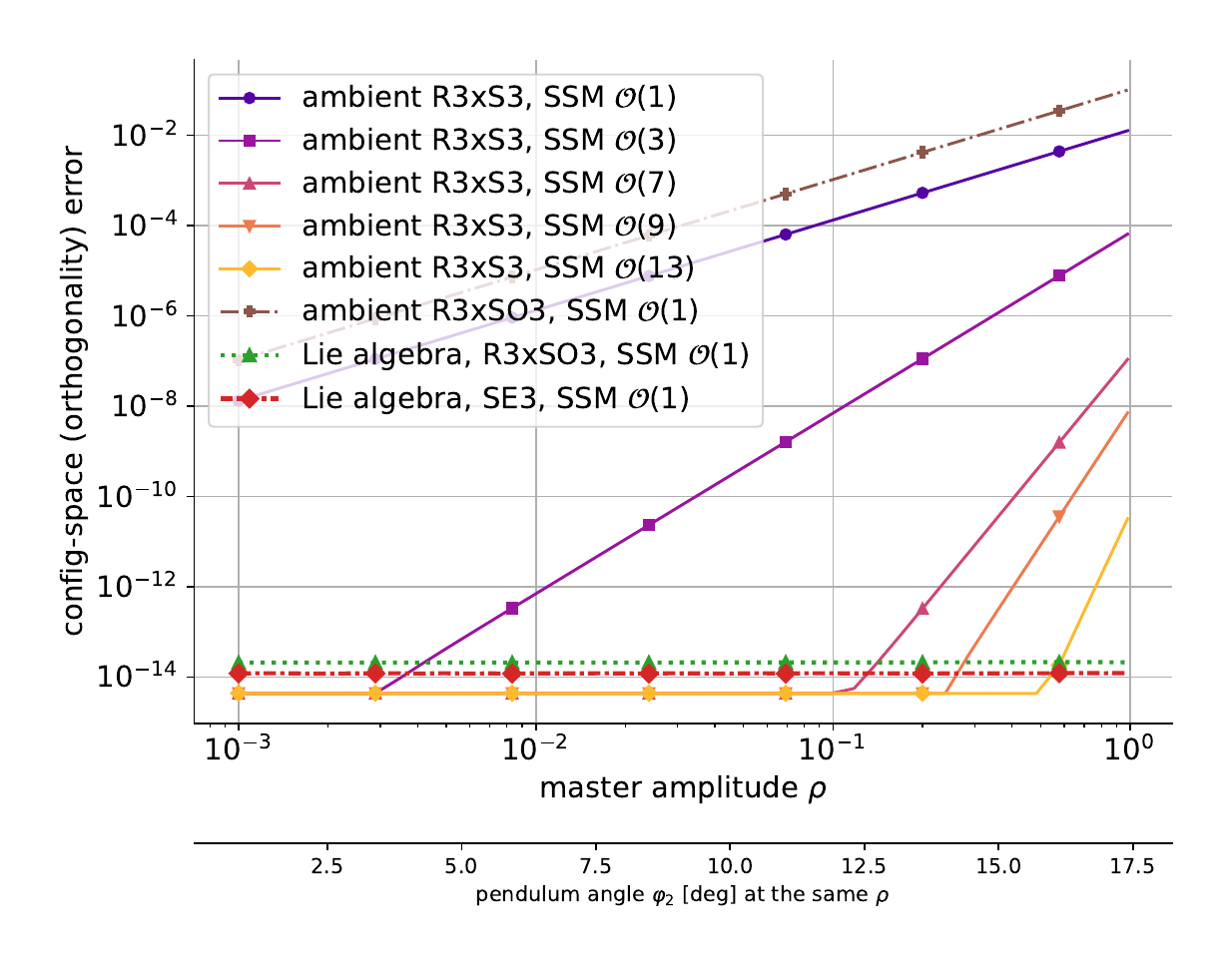}
    \vspace{-0.75cm}
    \caption{Error in the configuration space}
    \label{fig_SPR_ConfigSpaceConstraint_OrderAndLie}
\end{subfigure}
\caption{Error in the joint and configuration-space constraints along the backbone orbits of the pendulum slider system with 1:3 internal resonance for different SSM orders. ``Ambient'' denotes the polynomial parameterization of \sect{sec_LinearConstraints}, ``Lie algebra'' the rotation-preserving parameterization of \sect{sec_LieParameterization}; $\RealNumbers^3\times S^3$, $\RealNumbers^3\times\SOThree$, and $\SEThree$ denote the configuration space of the full model. }
\label{fig_SPR_Errors}
\end{figure}

To investigate constraint preservation in SSM reduced-order models, we consider two configuration-space representations of a rigid body based on the direct and semidirect products of $\RealNumbers^3$ with $SO(3)$ and $\SThree$, respectively.
The resulting errors in the kinematic joint and configuration-space constraints are shown in \fig{fig_SPR_Errors}.
As can be seen in \fig{fig_SPR_JointConstraint_OrderAndLie}, the $\RealNumbers^3 \times \SThree$-based SSM ROM violates both the joint and
configuration-space constraints at SSM order $\bigO(1)$, with the constraint error decreasing as the SSM order increases, as expected. 
In contrast, the $\RealNumbers^3 \times SO(3)$-based SSM ROM satisfies the joint constraints already at SSM order $\bigO(1)$ but violates the configuration space constraints, in agreement with the theoretical result in \sect{sec_LinearConstraints}.
The rotation-preserving parameterization-based SSM ROM, denoted as ``Lie algebra'' in \fig{fig_SPR_Errors}, satisfies both the joint and configuration-space constraints to machine precision at all considered SSM orders.

\subsubsection{A chain of pendulums}\label{sec_PendulumChain}
The next pendulum system we consider is the pendulum chain presented in \cite[Sect. 5.2.3]{li2023model} consisting of 41 pendulum bodies and hence represents a relatively large rigid multibody system.
\cite{li2023model} derived the explicit equations of motion for this systems, which we use as full order model here for comparison with our Lie group SSM reduced-order models.
Following \cite{li2023model}, we select the slowest spectral subspace spanned by the first mode as the master subspace for SSM reduction, owing to the absence of (near) internal resonances.

\begin{figure}
\centering
\includegraphics[width=0.5\textwidth]{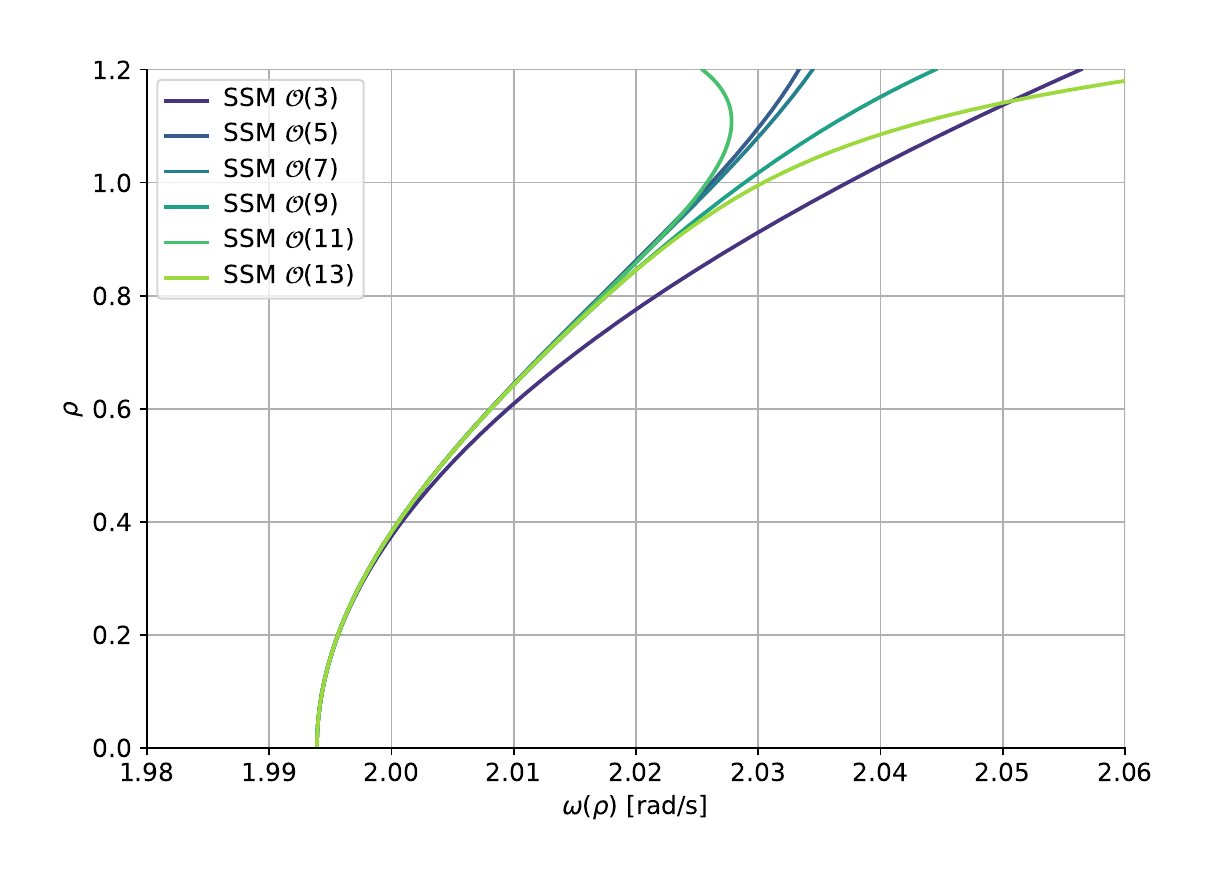}
\vspace{-0.5cm}
\caption{Backbone curves in reduced coordinates for the chain
of pendulums attached to a slider.}
\label{fig_Chain_Backbone_Fig13}
\end{figure}

The damped backbone curves obtained directly from the reduced dynamics for different orders of SSM approximation are shown in \fig{fig_Chain_Backbone_Fig13}, which agree well with the backbone curves reported in \cite[Fig.~13]{li2023model}.
Similarly to the previous examples, we observe that higher-order expansions are useful to obtain convergence in backbones at higher amplitudes. 
In particular, the backbone curves in \fig{fig_Chain_Backbone_Fig13} show convergence at $\bigO(3)$ for $\rho \leq 0.3$, at $\bigO(5)$ expansion for $\rho \leq 0.6$, and at $\bigO(13)$ expansion for $\rho \leq 0.8$. 
The backbone curves do not seem to converge for $\rho \geq 1$, indicating $\rho = 1$ is outside the domain of convergence for the power series $\omega(\rho)$ of the SSM ROM.

\begin{figure}%[h]
\centering
\begin{subfigure}{0.49\textwidth}
    \includegraphics[width=0.95\textwidth]{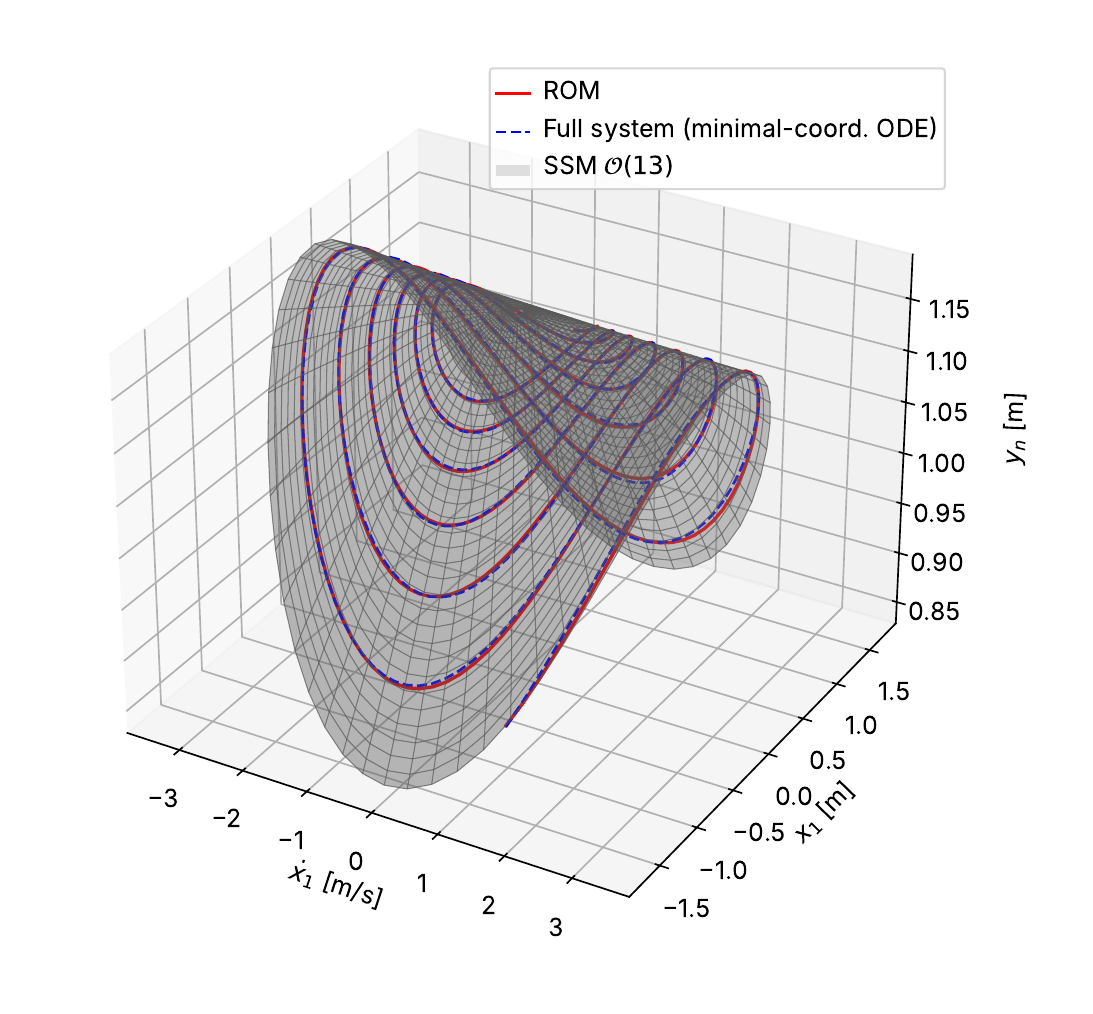}
    \vspace{-0.5cm}
    \caption{Slider position}
    \label{fig_Chain_SSMInvariance_yN}
\end{subfigure}
\hfill
\begin{subfigure}{0.49\textwidth}
    \includegraphics[width=0.95\textwidth]{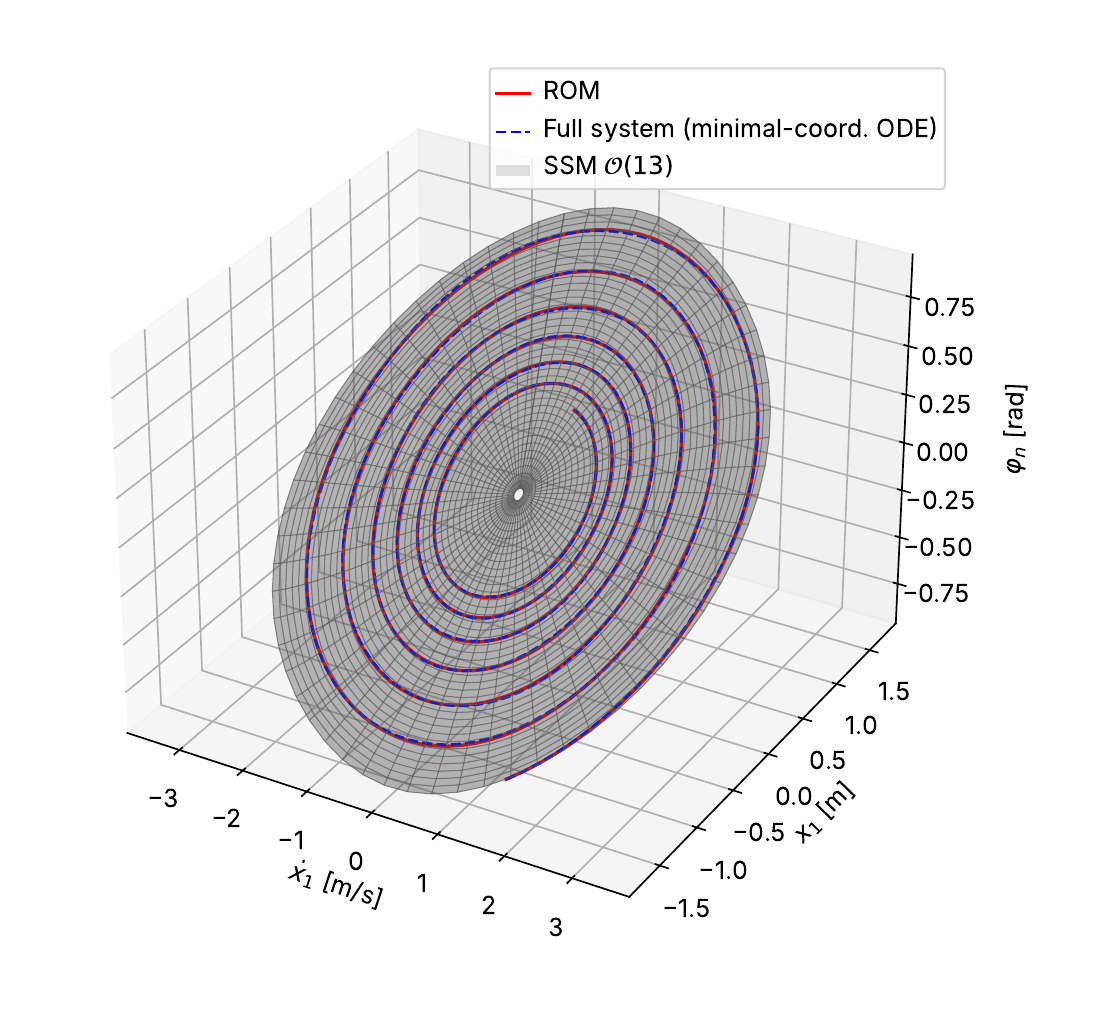}
    \vspace{-0.5cm}
    \caption{Pendulum angle}
    \label{fig_Chain_SSMInvariance_phiN}
\end{subfigure}
\caption{Projections of the SSM at $\bigO(13)$ approximation for the pendulum chain. 
The thin solid gray curves represent contour lines of
equal parameterized distance $\rho$ and $\vartheta$.}
\label{fig_Chain_SSM}
\end{figure}

To validate convergence, we illustrate the invariance of the SSM approximated up to $\bigO(13)$ within the convergence domain $\rho \leq 0.8$ in \fig{fig_Chain_SSM}.
To check invariance, we pick an initial condition $\pv_0 = (\rho_0 e^{i\vartheta_0},\, \rho_0 e^{-i\vartheta_0})$ on the SSM with $\rho_0 = 0.8$ and $\vartheta_0 = 3$, and perform time integration of both the SSM-based ROM and the full system in the index-1 formulation given in \cite{li2023model}. 
Indeed, the full-system trajectory stays on the computed SSM and overlaps well with the prediction of the SSM-based ROM as can be seen in \fig{fig_Chain_SSM}.

%%%%%%%%%%%%%%%%%%%%%%%%%%%%%%%%%%%%%%%%%%%%%%%%%%%%%%%%%%%%%%%%%%%%%%%%%%%%%%%%
\subsubsection{Flexible pendulum-slider demonstrating the Exudyn coupling}
\label{sec_PendulumSliderFlexible}
%%%%%%%%%%%%%%%%%%%%%%%%%%%%%%%%%%%%%%%%%%%%%%%%%%%%%%%%%%%%%%%%%%%%%%%%%%%%%%%%
%
The final pendulum system is modeled entirely in Exudyn, with a flexible pendulum based on the nodal-based floating frame of reference formulation, as shown schematically in \fig{fig_flexiblePendulumSliderExu}.
A rigid cubic slider of mass $1$\,kg and edge length $0.2$\,m moves along the global $x$ axis on a prismatic joint and is held by a linear spring-damper with $k_\textrm{s}=100$\,N/m and $c_\textrm{s}=0.5$\,N\,s/m.
The pendulum, a rod of length $L=1$\,m, square cross-section $0.1\times0.1$\,m, mass $1$\,kg, Young's modulus $E=3\cdot10^{7}$\,Pa, and Poisson's ratio $0.3$, is attached to the slider by a revolute joint about the $z$ axis and hangs in the direction of gravity, $g=9.81$\,m/s$^2$.
The rod is meshed with second-order tetrahedral elements in NGSolve \cite{Schoeberl2014}, using its mesh generator Netgen \cite{Schoeberl1997}, and the mass and stiffness matrices required for the floating frame of reference body are exported from NGSolve via Exudyn's FEM interface.
The body is then reduced by the Hurty--Craig--Bampton method with the pivot face as interface and a single fixed-interface eigenmode.
A force $-c\,\dot u-k_\textrm{L}u-k_\textrm{NL}u^3$ with $c=0.25$\,N\,s/m, $k_\textrm{L}=10^3$\,N/m, and $k_\textrm{NL}=10^4$\,N/m$^3$ acts at the mesh node closest to the centre of mass of the rod, where $u$ denotes the $x$ displacement of this node from its equilibrium position.
The resulting model has $15$ coordinates, namely Euler parameters for the slider and for the floating frame of the rod plus one modal coordinate, and $12$ algebraic equations.
The model, including the spring-dampers defined as user functions, exists only in Exudyn, and the reduction accesses its residual through the coupling of \sect{sec:external-residual-coupling}.
The linearization yields the slider mode at $9.43$\,rad/s, the swing mode of the pendulum at $27.96$\,rad/s with a damping ratio of $0.40\,\%$, and the first bending mode at $86.6$\,rad/s.
The complex-conjugate eigenvalue pair of the swing mode spans the two-dimensional master subspace ($n_r=2$), and the SSM is computed up to order~$9$.
 
\begin{figure}%[h]
\centering
\begin{subfigure}{0.49\textwidth}
    \includegraphics[width=0.95\textwidth]{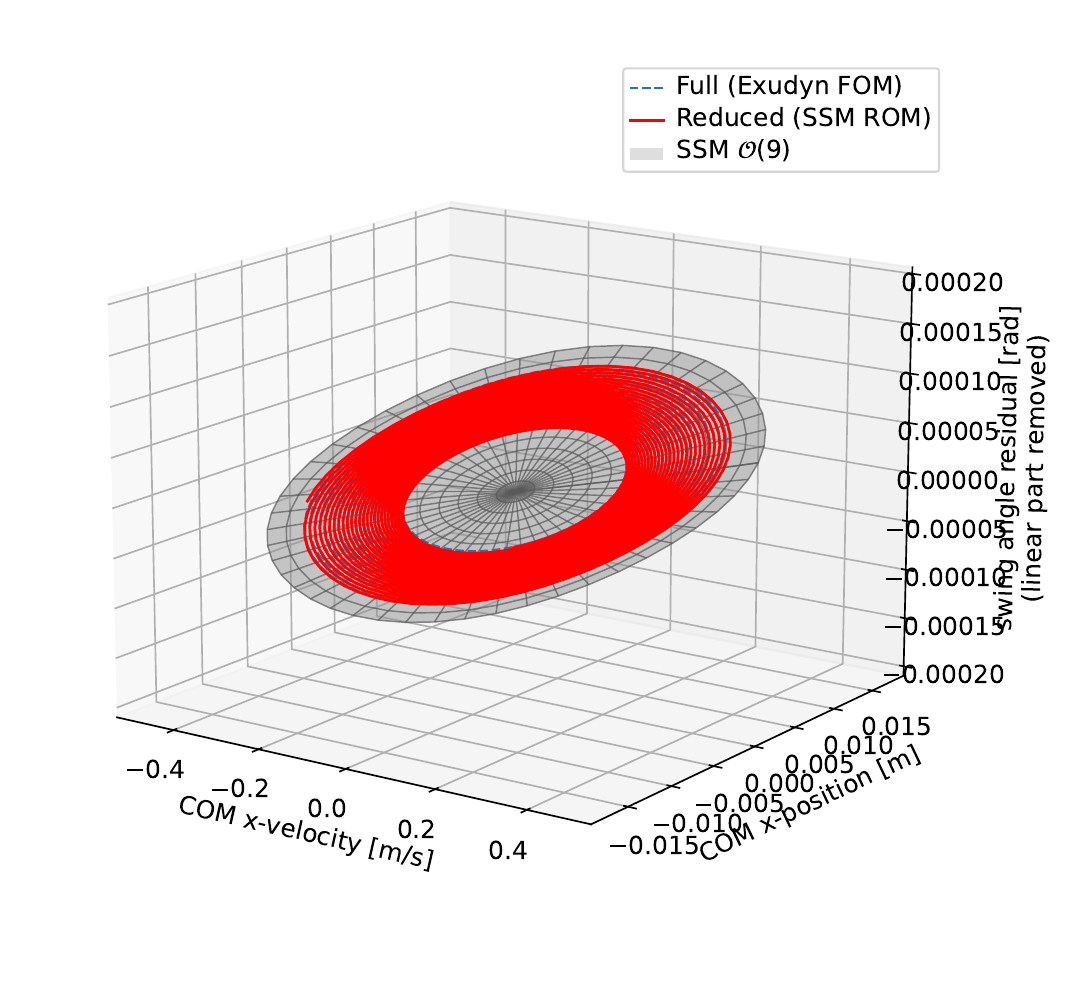}
    \vspace{-1cm}
    \caption{Swing angle}
    \label{fig_SPF_SSM_angle_zoom}
\end{subfigure}
\hfill
\begin{subfigure}{0.49\textwidth}
    \includegraphics[width=0.95\textwidth]{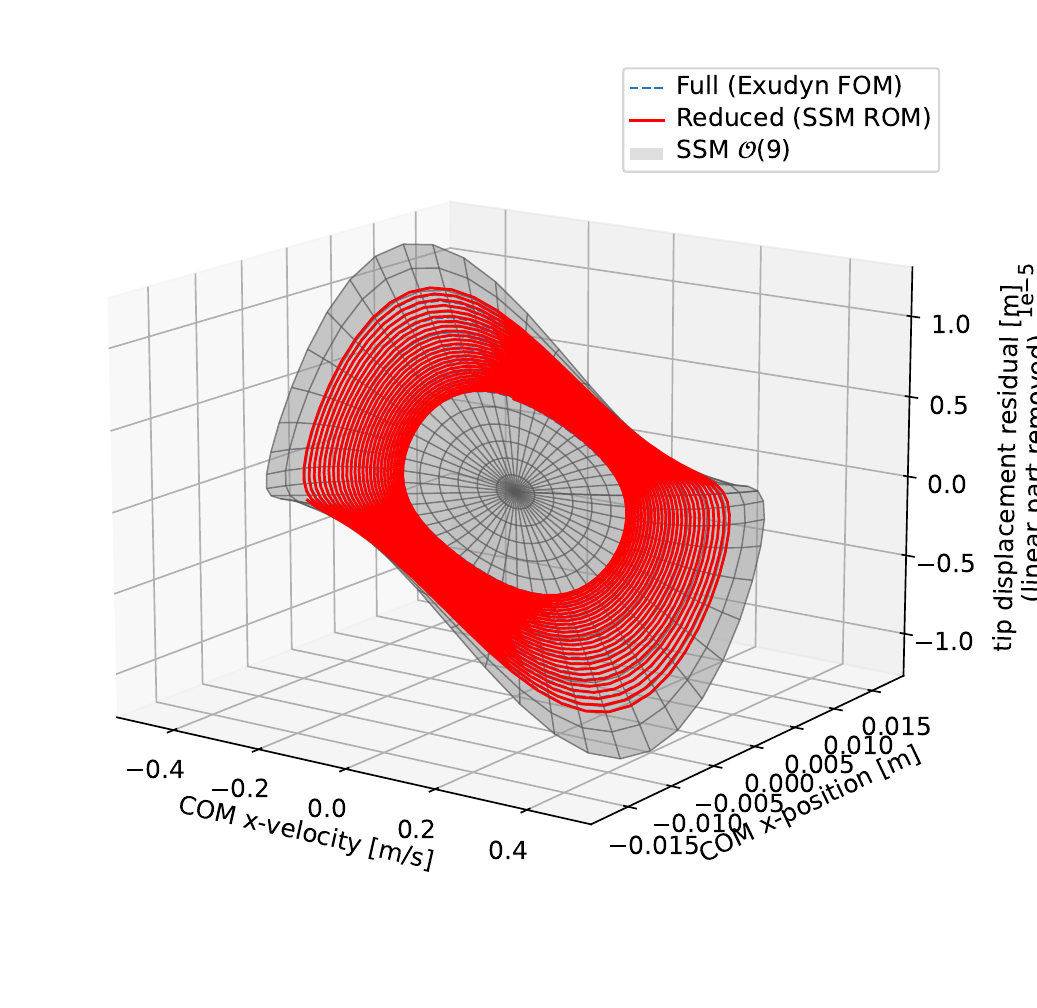}
    \vspace{-1cm}
    \caption{Elastic displacement of the tip node}
    \label{fig_SPF_SSM_defl_zoom}
\end{subfigure}
\caption{Projections of the SSM at $\bigO(9)$ approximation for the flexible pendulum-slider onto the velocity and position of the centre of mass along $x$ and (a) the swing angle and (b) the transverse elastic displacement of the tip node, each with its part linear in the centre-of-mass position removed. Full model (dashed) and reduced model (solid) start on the SSM at $\rho_0=1.5$. The thin solid gray curves represent contour lines of equal parameterized distance $\rho$ and $\vartheta$.}
\label{fig_SPF_SSMs}
\end{figure}

\begin{figure}%[h]
\centering
\includegraphics[width=0.6\textwidth]{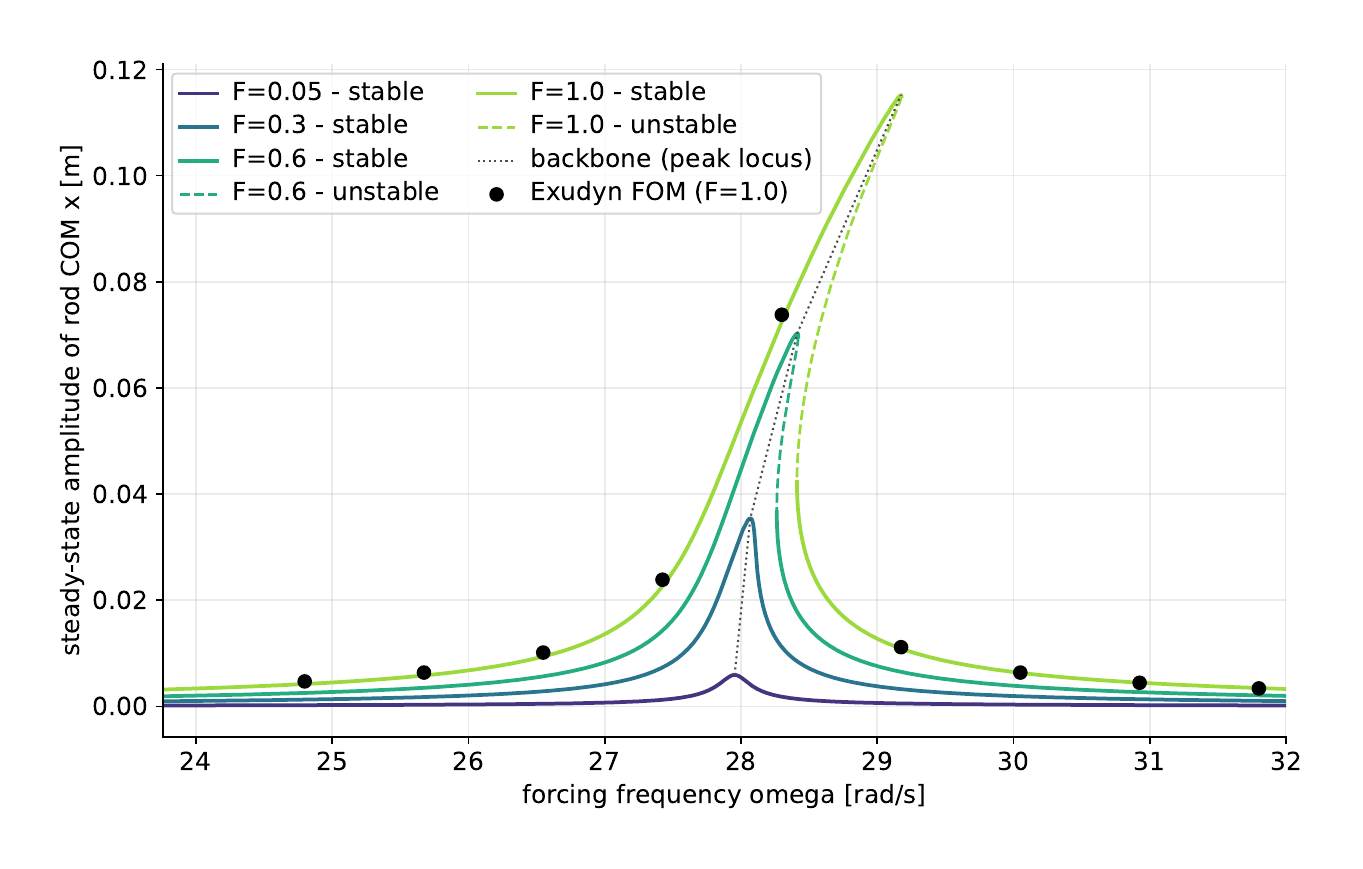}
\vspace{-0.35cm}
\caption{Frequency response curves of the flexible pendulum-slider: amplitude of the centre-of-mass position along $x$ of the reduced model at $\bigO(9)$ for four forcing amplitudes and full-order steady states at $F=1.0$\,N (dots).}
\label{fig_SPF_FRC}
\end{figure}

\fig{fig_SPF_SSMs} shows the SSM together with a free decay of $6$\,s started on the SSM at $\rho_0=1.5$, computed with the reduced model and with Exudyn's generalized-$\alpha$ integrator ($\rho_\infty=0.9$, step size $2\cdot10^{-4}$\,s).
Both projections use the position and velocity of the centre of mass along $x$ as the first two coordinates, and either the swing angle or the transverse elastic displacement of the tip node as the third.
From the third coordinate, its part linear in the centre-of-mass position is subtracted, which removes the tangent plane of the SSM and exposes its curvature.
The trajectories of the reduced and the full model coincide within line thickness in both projections, including the elastic tip displacement, which is of order $10^{-5}$\,m and therefore the most sensitive observable.
For harmonic forcing $F\cos(\Omega t)$ along $x$ at the same mesh node, \fig{fig_SPF_FRC} shows frequency response curves for $F=0.05$, $0.3$, $0.6$, and $1.0$\,N.
The response hardens, and the peak shifts from $27.96$\,rad/s at $F=0.05$\,N to $29.33$\,rad/s at $F=1.0$\,N, where the curve folds.
Full-order steady states at $F=1.0$\,N, each started from rest and integrated over $\max(60\,T,\,6/(\zeta\omega_0))$ with $200$ steps per period, agree with the reduced model and settle on the lower branch beyond the fold.

%%%%%%%%%%%%%%%%%%%%%%%%%%%%%%%%%%%%%%%%%%%%%%%%%%%%%%%%%%%%%%%%%%%%%%%%%%%%%%%%
\subsection{Floating satellite}
\label{sec:example-satellite}
\begin{figure}%[t]
\centering
\includegraphics[width=0.85\textwidth]{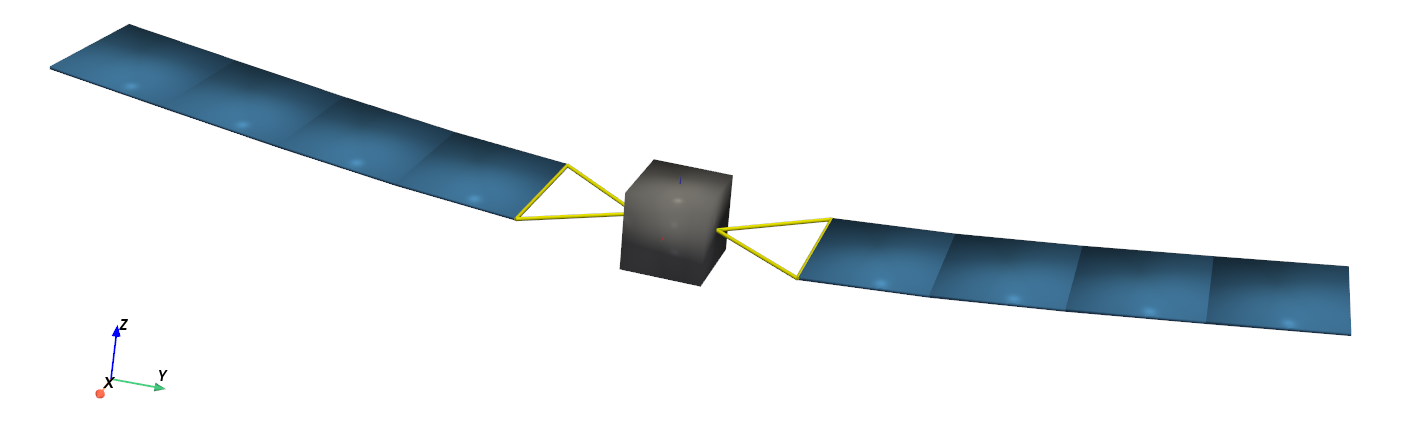}
\caption{Satellite with two solar-array wings: bus, yokes, panels, hinges, attitude controller, and forcing.}
\label{fig:satellite-model}
\end{figure}

%%%%%%%%%%%%%%%%%%%%%%%%%%%%%%%%%%%%%%%%%%%%%%%%%%%%%%%%%%%%%%%%%%%%%%%%%%%%%%%%
%
The satellite in \fig{fig:satellite-model} consists of a cubic bus and two solar-array wings along the global $\pm y$ axis, each composed of a yoke and four panels.
The model is motivated by the satellite of \cite{Fufa2010}, from which the masses of bus, yokes, and panels are adopted. 
The dimensions of yokes and panels are adapted, and the flexible panels are replaced by rigid panels connected by elastic hinges.
The eleven homogeneous cuboids are connected in series by ten revolute joints, whose axes are parallel to the global $x$ axis in the reference configuration.
The hinges lie at $y=\pm0.25$, $\pm0.85$, $\pm1.61$, $\pm2.37$, and $\pm3.13$\,m, and gravity is absent.
Each hinge carries the torque $\tau_j=k_j\theta_j+k_{3,j}\theta_j^3+c\,k_j\dot\theta_j$ in the relative hinge angle $\theta_j$, which is included as an additional coordinate whose rate equals the relative angular velocity about the hinge axis.
The cubic coefficient $k_{3,j}=\kappa k_j/\theta_\textrm{ref}^2$ represents the stiffening of deployment hinges.
The bus is held by a proportional-derivative attitude controller, $\boldsymbol{\tau}_\textrm{att}=-\omega_\textrm{att}^2\mathbf{J}\operatorname{vee}\bigl(\tfrac12(\RotMat_\textrm{b}-\RotMat_\textrm{b}\tp)\bigr)-2\zeta_\textrm{att}\omega_\textrm{att}\mathbf{J}\boldsymbol{\omega}_\textrm{b}$, which is polynomial in the bus rotation matrix $\RotMat_\textrm{b}$, and by a weak translational spring-damper.
Here, $\boldsymbol{\omega}_\textrm{b}$ denotes the inertial angular velocity of the bus and $\mathbf{J}$ the inertia of the rigid satellite about the bus centre.
All parameters are listed in \tab{tab:satellite-parameters}.
All bodies use the rotation-matrix representation of $\RealNumbers^3\times\SOThree$, which gives $\dimq=142$, $\dimG=76$, $\dimJnt=60$, $\dimCsp=66$, and $\dimSS=344$.
A harmonic force $F\cos(\Omega t)$ acts at the bus centre along the global $z$ axis, and the reported quantity is the angle $\theta_1$ of the hinge between yoke and first panel of the $+y$ wing.
\begin{table}%[t]
\centering
\caption{Parameters of the satellite. Hinge stiffnesses of the $-y$ wing are multiplied by $0.93$; panel hinges $i=1,\ldots,4$ are numbered outward.}
\label{tab:satellite-parameters}
\small
\begin{tabular}{ll}
\hline
bus & $680$\,kg, $0.5\times0.5\times0.5$\,m \\
yoke, panel ($x\times y\times z$) & $3.3$\,kg, $0.8\times0.6\times0.02$\,m; \ $5.0$\,kg, $0.8\times0.76\times0.02$\,m \\
hinge stiffness $k_j$ & root: $5$\,N\,m/rad; \ panel $i$: $2\,(1+0.02\,(i-1))$\,N\,m/rad \\
hinge nonlinearity, damping & $\kappa=0.3$, $\theta_\textrm{ref}=0.1$\,rad, $c=0.19$\,s \\
attitude controller & $\omega_\textrm{att}=2\pi\cdot0.01$\,rad/s, $\zeta_\textrm{att}=0.7$, $\mathbf{J}=\operatorname{diag}(286.01,30.82,288.49)$\,kg\,m$^2$ \\
translational support & $0.115$\,N/m, $9.13$\,N\,s/m \ ($2$\,mHz, $\zeta=0.5$, total mass $726.6$\,kg) \\
\hline
\end{tabular}
\end{table}

The support and attitude modes lie at $1.7$, $7.1$, and $9.2$\,mHz with damping ratios of at least $0.5$, followed by the elastic modes at $19.52$, $120.8$, and $126.7$\,mHz.
The first elastic mode, $\lambda_{1,2}=-0.00177\pm0.12263\,\mathrm{i}$\,rad/s with a damping ratio of $1.44\,\%$, is the symmetric flapping of both wings accompanied by a translation of the bus along $z$; it forms the master subspace ($n_r=2$).
SSMs of orders $3$, $5$, $7$, and $9$ are computed in the ambient formulation, and frequency response curves are obtained by continuation of the forced reduced dynamics.
The full-order reference uses the generalized-$\alpha$ Lie group integrator \cite{Holzinger2025} with $\rho_\infty=0.9$ and step size $0.25$\,s\footnote{Doubling the step size changes the forced amplitude by $0.4\,\%$.}.
A free decay of $1500$\,s, computed with step size $0.05$\,s, starts on the SSM of order~7 at $\rho_0=0.6$ after projection onto all position- and velocity-level constraints, and its backbone is extracted from the extrema of $\theta_1(t)$.
Stepped-sine sweeps at $F=0.1$\,N in both directions over $13$ frequencies between $18.0$ and $21.5$\,mHz start each frequency from the final state of the previous one and integrate over $2500$\,s ($3000$\,s for the first frequency); the amplitude is half the peak-to-peak value of $\theta_1$ over the last three periods.
 
\begin{figure}%[t]
\centering
\begin{subfigure}{0.49\textwidth}
\includegraphics[width=\textwidth]{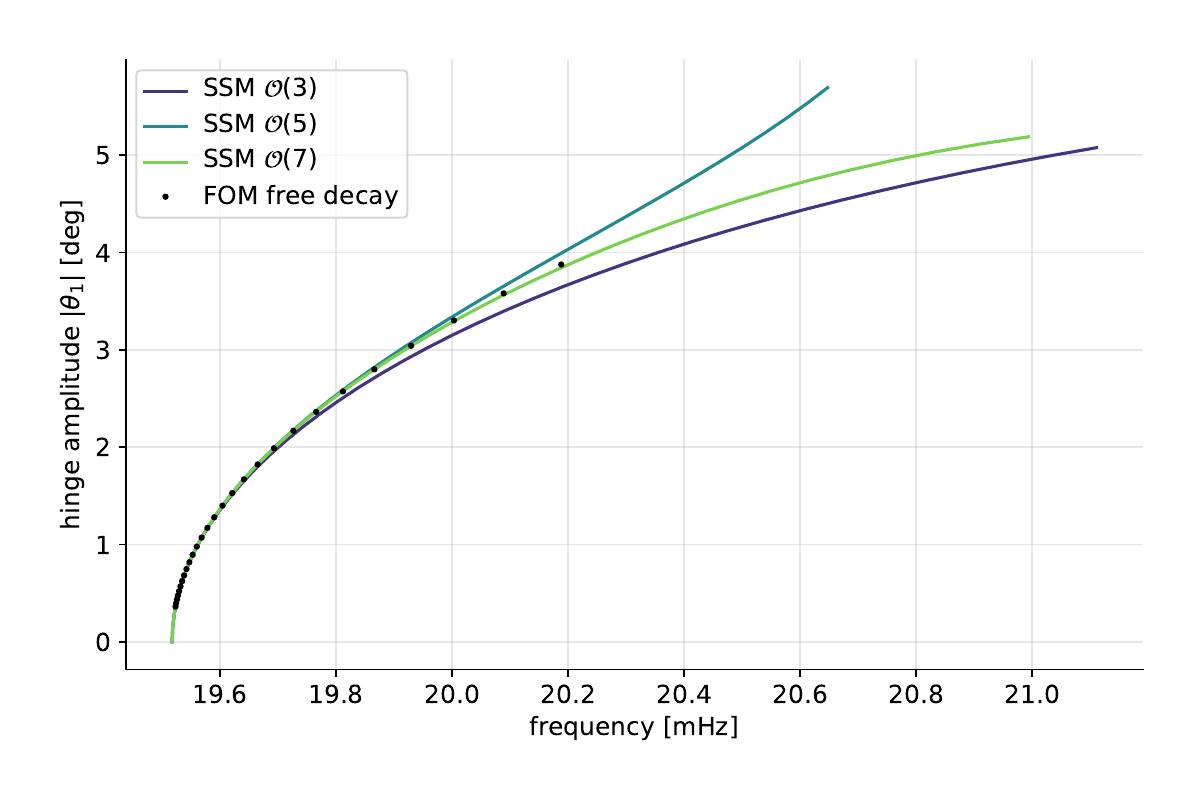}
\vspace{-0.75cm}
\caption{Backbone curve}
\end{subfigure}
\hfill
\begin{subfigure}{0.49\textwidth}
\includegraphics[width=\textwidth]{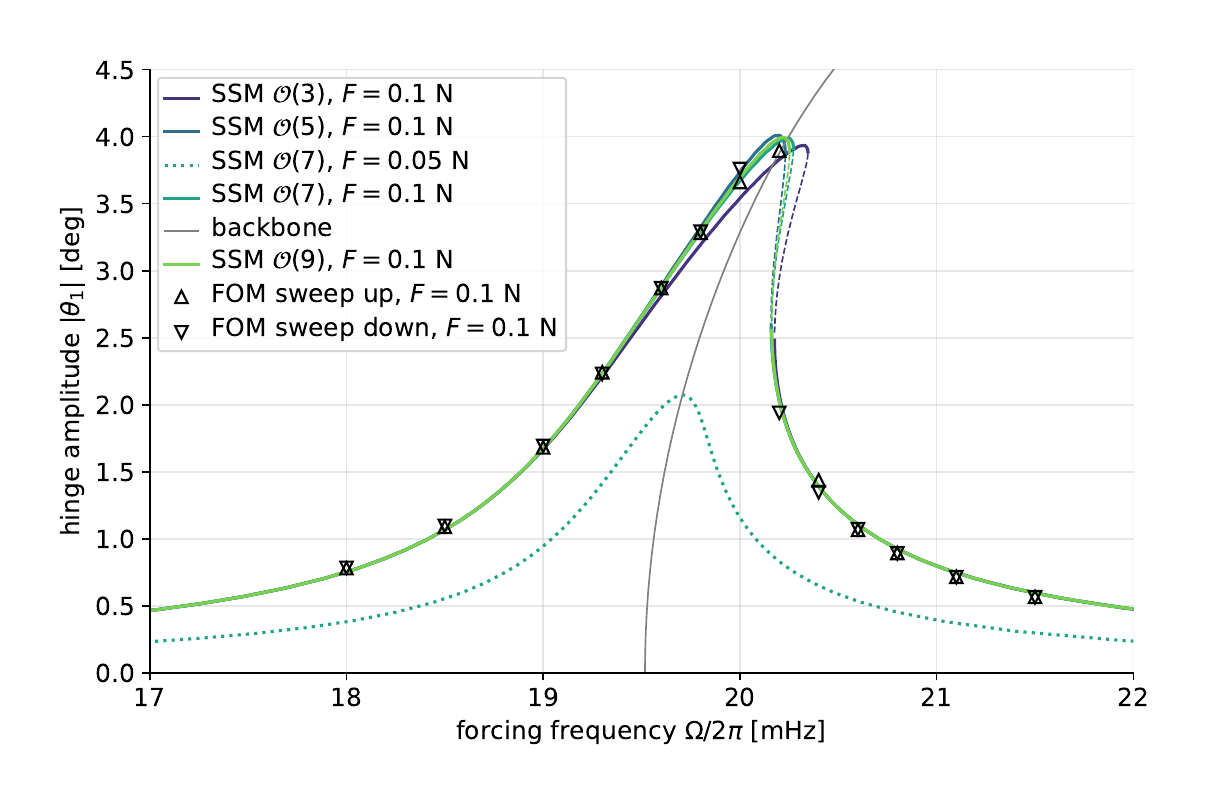}
\vspace{-0.75cm}
\caption{Frequency response function}
\end{subfigure}
\caption{Satellite: (a) backbone curves of the reduced models and of a full-order free decay; (b) frequency response curves of the reduced models (solid: stable, dashed: unstable) and full-order stepped-sine sweeps at $F=0.1$\,N ($\triangle$: increasing, $\triangledown$: decreasing frequency).}
\label{fig:satellite-response}
\end{figure}

\begin{figure}%[t]
\centering
\includegraphics[width=0.65\textwidth]{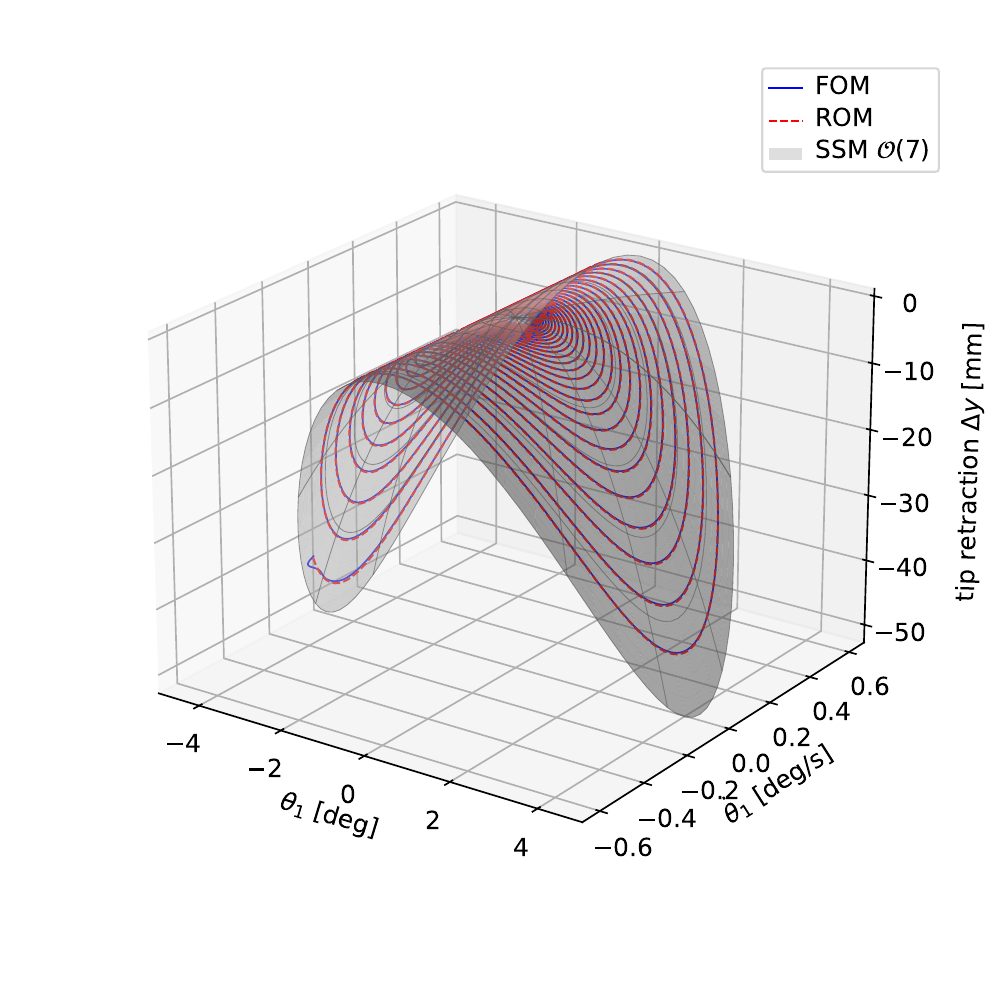}
\vspace{-1.25cm}
\caption{Satellite: SSM of order~7 in the coordinates $\theta_1$, $\dot\theta_1$, and the retraction $\Delta y$ of the $+y$ wing tip along the wing, with the free decay of the reduced and of the full-order model.}
\label{fig:satellite-ssm}
\end{figure}

The backbone curves in \fig{fig:satellite-response}a show a hardening response, with a frequency increase of $2.9\,\%$ at $3.57^\circ$ for order~7.
Orders~3, 5, and~7 coincide up to about $2.5^\circ$, and order~7 agrees with the full-order free decay from $3.88^\circ$ down to $0.36^\circ$.
At $F=0.1$\,N, the frequency response in \fig{fig:satellite-response}b folds, and orders~7 and~9 nearly coincide; order~7 has saddle-node bifurcations at $20.16$ and $20.27$\,mHz and a peak of $3.99^\circ$.
The full-order sweeps reproduce the hysteresis: the up-sweep follows the upper branch to $20.2$\,mHz and drops at $20.4$\,mHz, whereas the down-sweep jumps up at $20.0$\,mHz.
Between $19.0$ and $20.0$\,mHz, the reduced model deviates from the full model by at most $1.1\,\%$, apart from $2.4\,\%$ directly after the jump.
Away from resonance the deviation grows to $5.7\,\%$, with the reduced model below the full model on the low-frequency side and above it on the high-frequency side.
This is consistent with the quasi-static response of the non-master modes, which the leading-order forced reduced model does not contain.
\fig{fig:satellite-ssm} shows the SSM together with both free decays in the coordinates $\theta_1$, $\dot\theta_1$, and the retraction $\Delta y$ of the wing tip, which is quadratic in the flapping amplitude and reflects the geometric nonlinearity of the finite rotations.

\begin{figure}%[t]
\centering
\begin{subfigure}{0.49\textwidth}
\includegraphics[width=\textwidth]{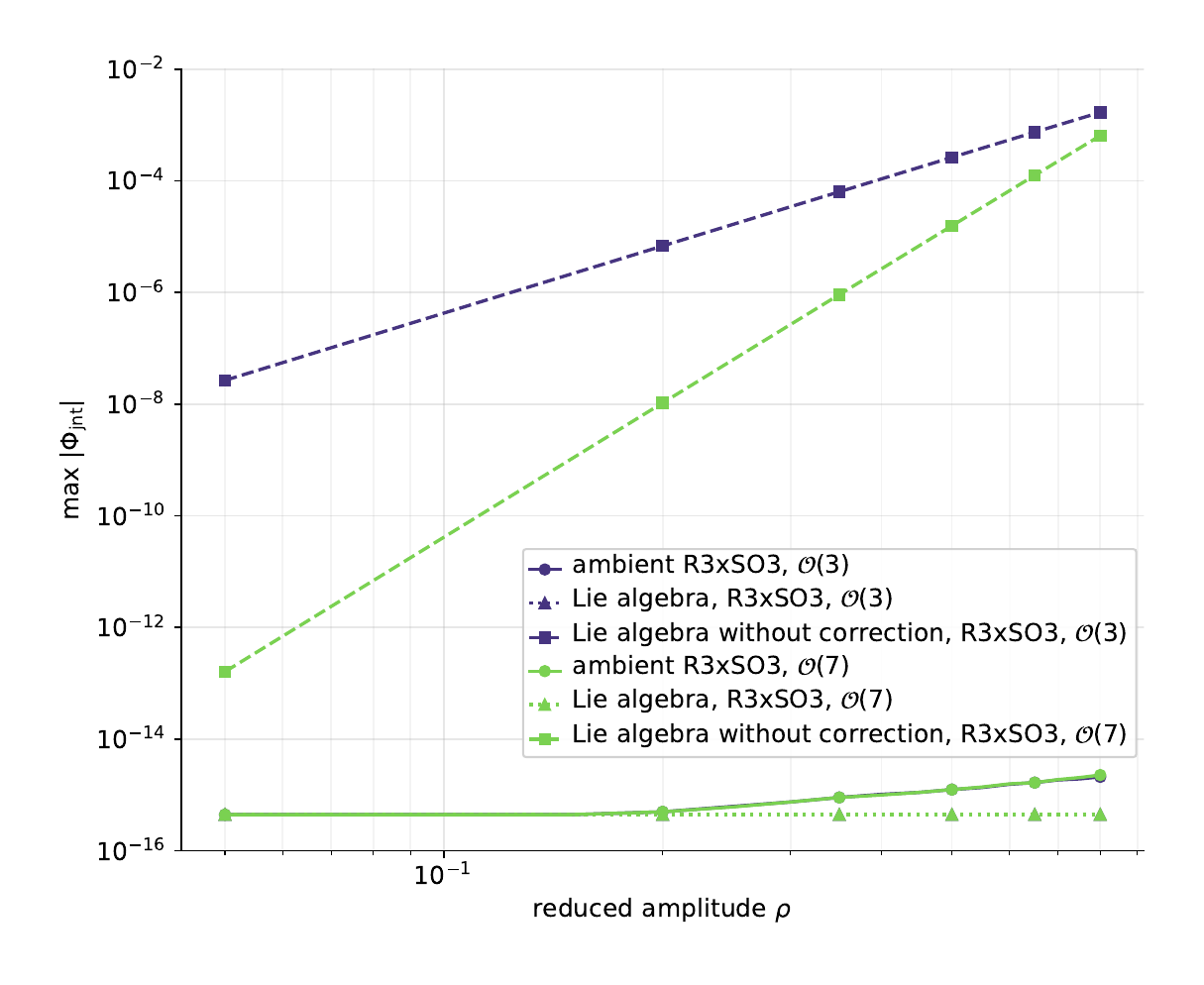}
\vspace{-0.75cm}
\caption{Joint constraints error}
\end{subfigure}
\hfill
\begin{subfigure}{0.49\textwidth}
\includegraphics[width=\textwidth]{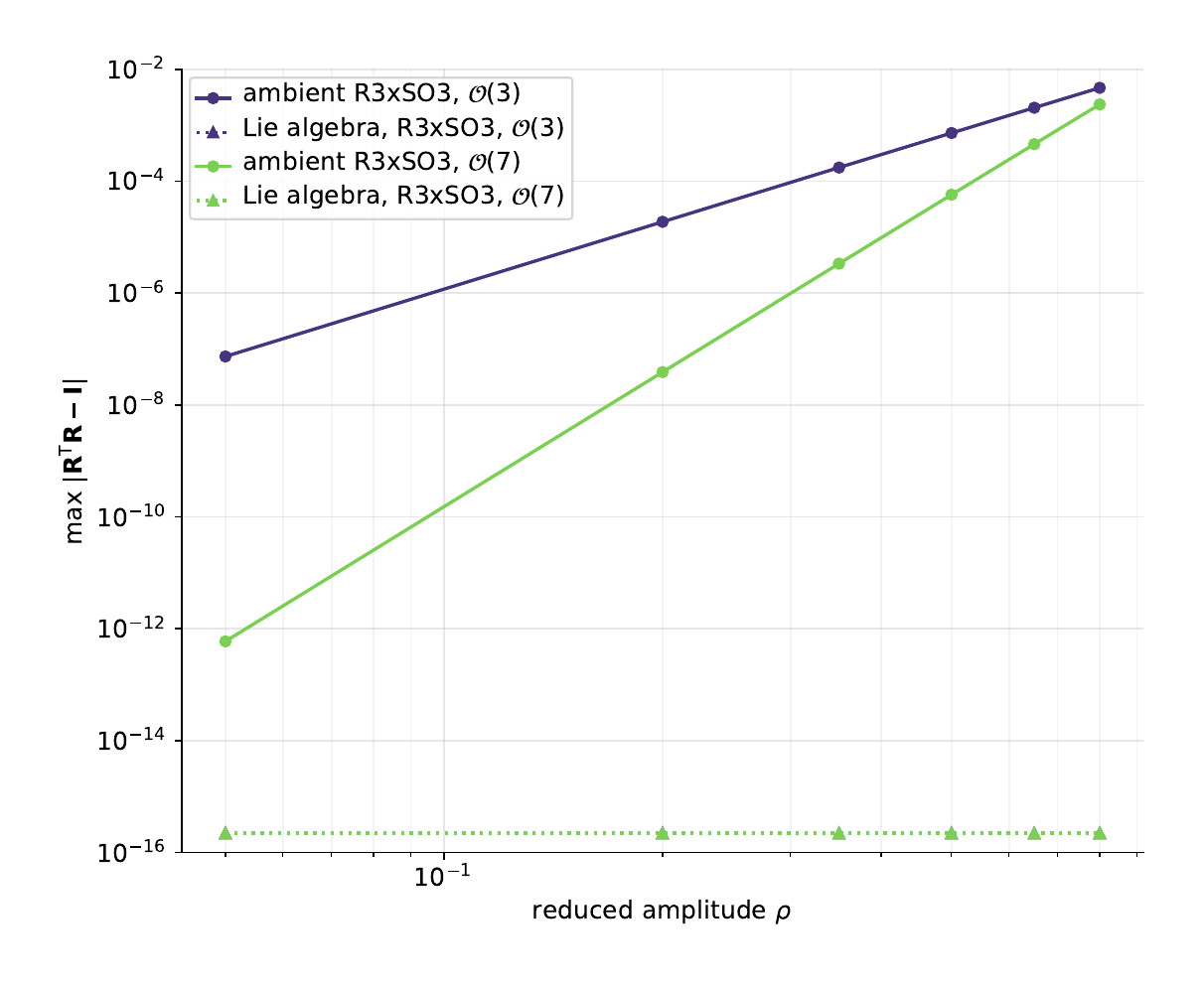}
\vspace{-0.75cm}
\caption{Configuration space constraints error}
\end{subfigure}
\caption{Satellite: maximum violation of (a) the joint constraints and (b) the orthogonality of the rotation matrices along the backbone orbits.}
\label{fig:satellite-constraints}
\end{figure}

\fig{fig:satellite-constraints} compares the constraint violations along the backbone orbits.
The ambient formulation satisfies the joint constraints to machine precision, but violates orthogonality by up to $2.3\cdot10^{-3}$ at order~7.
The rotation-preserving reconstruction \eq{eq_LieReconstruction} is orthogonal to machine precision but violates the joint constraints by up to $6.5\cdot10^{-4}$, since the bus rotates about all three axes and the joint axes along the kinematic tree are therefore not parallel.
The correction map \eq{eq_CorrectionMap} restores all joint constraints to $4.4\cdot10^{-16}$.
The SSM of order~7 is computed in about $4$\,s and a reduced frequency response curve in about $12$\,s, whereas one full-order steady state takes $4$ to $8$ minutes.

%%%%%%%%%%%%%%%%%%%%%%%%%%%%%%%%%%%%%%%%%%%%%%%%%%%%%%%%%%%%%%%%%%%%%%%%%%%%%%%%
\subsection{Ropeway system}
\label{sec_Ropeway}
\begin{figure}%[t]
\centering
\includegraphics[width=\textwidth]{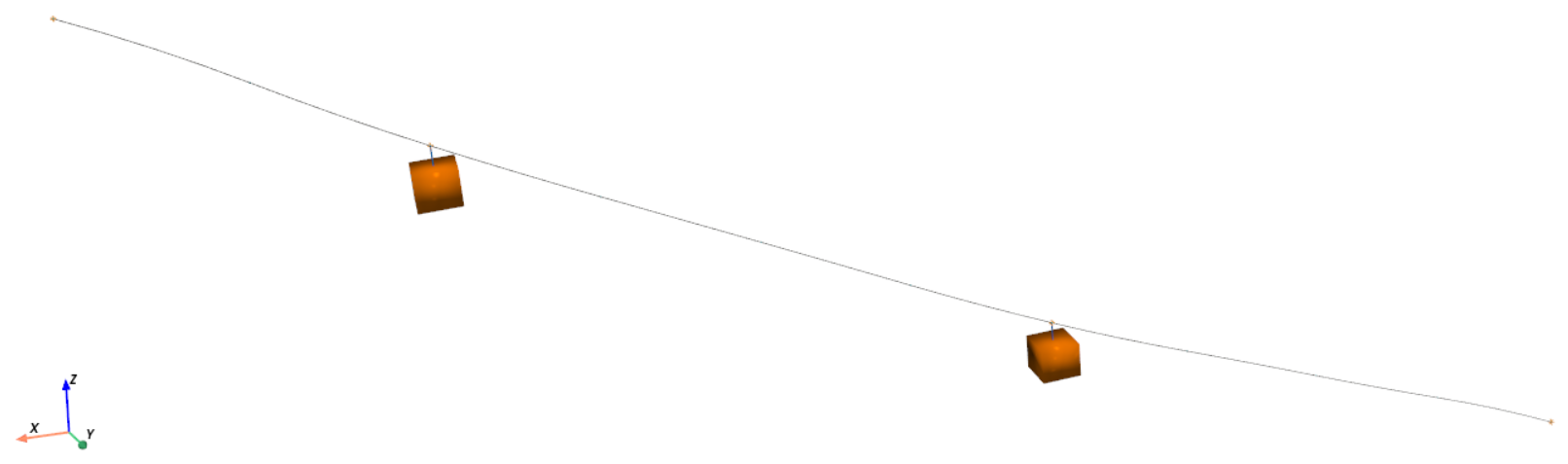}
\vspace{-0.5cm}
\caption{Section of a pre-tensioned 3S (tricable) ropeway system with two gondolas (cable cars).}
\label{fig_ropewaySystem}
\end{figure}

%%%%%%%%%%%%%%%%%%%%%%%%%%%%%%%%%%%%%%%%%%%%%%%%%%%%%%%%%%%%%%%%%%%%%%%%%%%%%%%%
%
The ropeway example reduces a section of a pre-tensioned 3S (tricable) ropeway, as shown schematically in \fig{fig_ropewaySystem}.
The track rope is a $60$\,mm locked-coil rope spanning $120$\,m between two towers with an elevation difference of $50$\,m.
It is discretized by $10$ polynomial geometrically exact $SE(3)$ beam elements.
The rope is pre-tensioned to $584$\,kN, corresponding to $15.7\,\%$ of its breaking load, by shortening the reference length of each element.
Span and tension lie in the range reported for an existing installation \cite{Becu2024}.
Two rigid $3$S cabins of mass $4600$\,kg are each suspended from the rope, at $40\,\%$ and $80\,\%$ of the span, by a rigid hanger of mass $120$\,kg and length $1.40$\,m through a revolute joint with only viscous damping and no restoring stiffness about the axis aligned with the rope.
The resulting model has $360$ states.
The linearization yields a dense cluster of lateral and longitudinal gondola-pendulum modes below $0.3$\,Hz.
Above this cluster, a well-separated vertical mode of the rope and gondolas occurs at $0.4895$\,Hz with a damping ratio of $0.55\,\%$.
Its complex-conjugate eigenvalue pair spans the master subspace ($n_r=2$), and the SSM is computed up to order~$3$ by the contour extraction of \sect{sec:implementation-aspects-contour}.
Contour extraction is used here because the Taylor-typed evaluation of \sect{sec:implementation-aspects-taylor} through the beam element, whose residual composes minimax polynomials of degree up to 19 inside the quadrature loop, required more than 2~GB of memory already at order 3 for a single element, whereas the contour sum is exact for the polynomial residual at any order and its memory does not grow with the order.
The circle radius is calibrated by a preliminary order-2 reduction such that it spans 0.05~m of gondola motion, and $N_c=48$ points are used; doubling $N_c$ leaves the coefficients unchanged. 
Order 3 suffices for the amplitudes considered, since the backbone coefficients of orders 3 and 5 agree to $10^{-10}$ relative and the predicted amplitudes to $0.4\,\%$.

\begin{figure}
\centering
\includegraphics[width=0.55\textwidth]{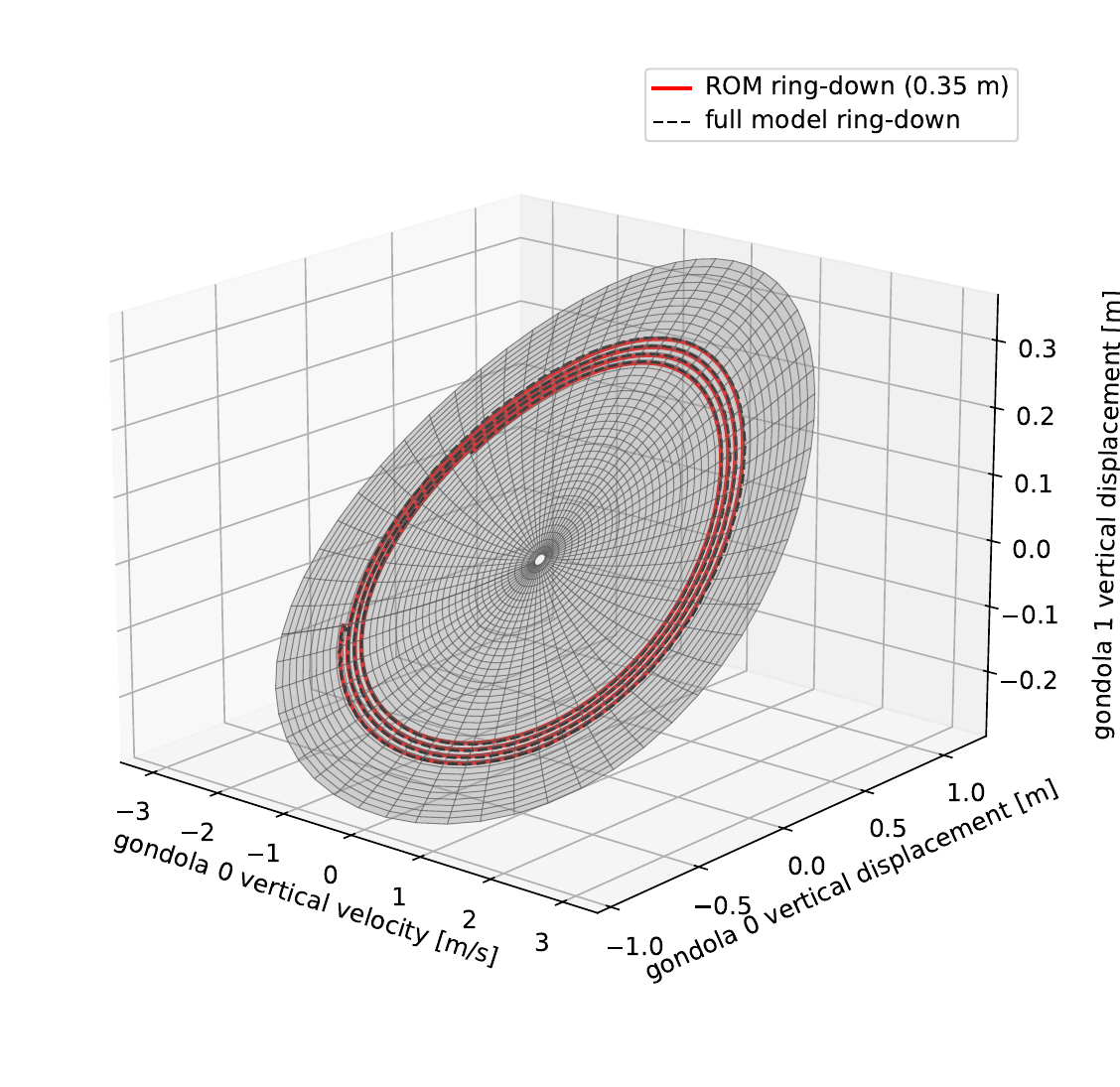}
\vspace{-0.85cm}
\caption{SSM at $\bigO(3)$ approximation for the ropeway, projected onto the vertical velocity and displacement of the first gondola and the vertical displacement of the second gondola. Full model (dashed) and reduced model (solid) start on the SSM at a $1$\,m vertical excursion of the first gondola.}
\label{fig_ropeway_SSM}
\end{figure}

\fig{fig_ropeway_SSM} shows the SSM together with a free decay of $8$\,s started at a vertical excursion of $1$\,m of the first gondola, computed with the reduced model and with the full multibody model (generalized-$\alpha$, $\rho_\infty=0.6$, step size\footnote{Halving the step size changes the free-decay trajectory by $0.065\,\%$ of the amplitude.} $5\cdot10^{-3}$\,s).
The coordinates are the vertical velocity and displacement of the first gondola together with the vertical displacement of the second gondola.
The trajectories of the reduced and the full model agree closely over the decay, with a slowly accumulating phase lag of the reduced model that reaches $5\,\%$ of the initial excursion at the end of the window.
For the excursion of $1$\,m, the frequency extracted from the decay lies $2.55\,\%$ below the linear value for the full model and $2.81\,\%$ below it for the reduced model.
For a smaller excursion of $0.5$\,m, the corresponding frequency shifts are $0.83\,\%$ and $0.84\,\%$, and the two trajectories differ by less than $0.5\,\%$ of the initial excursion.

\begin{figure}
\centering
\includegraphics[width=0.6\textwidth]{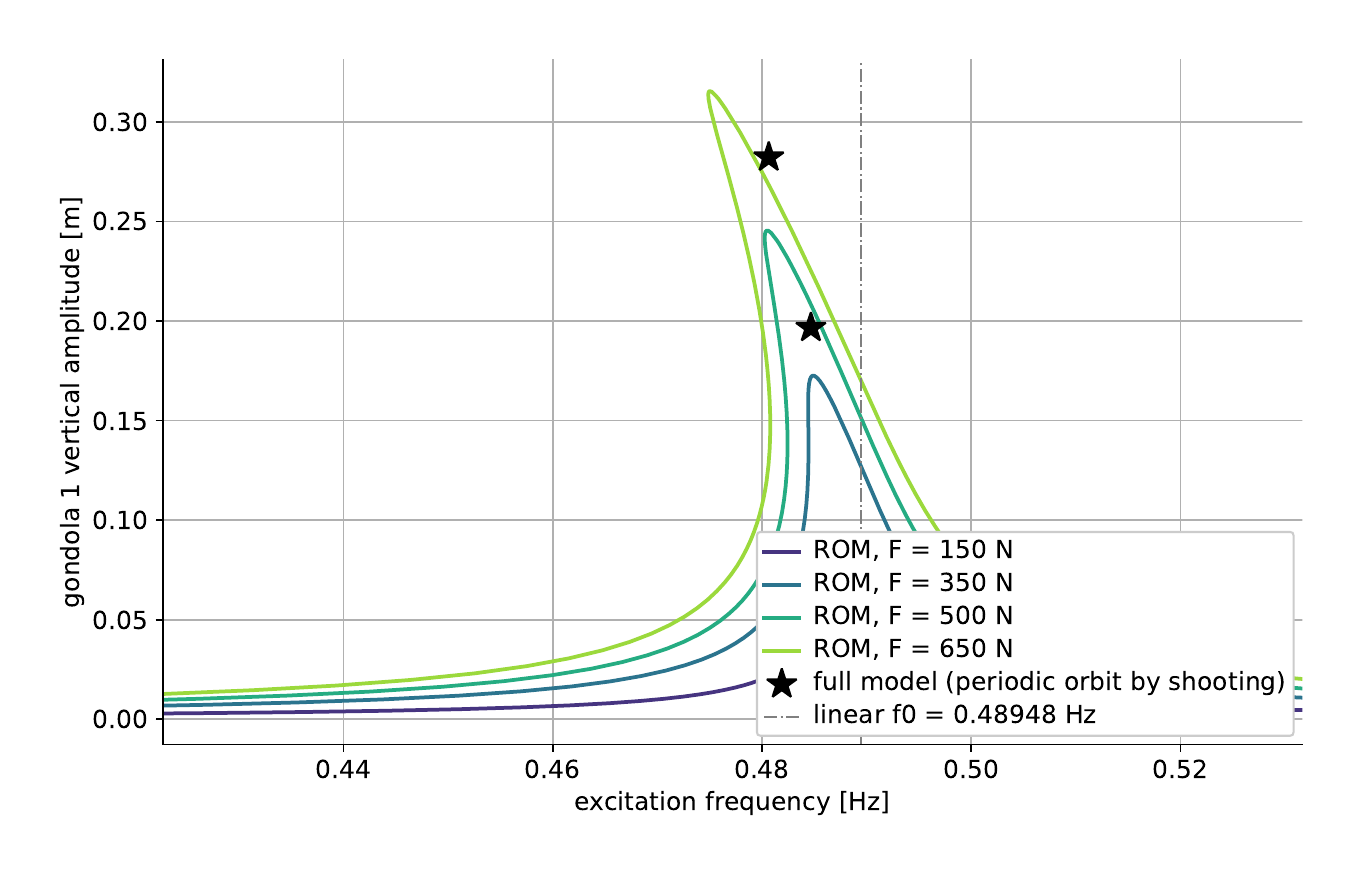}
\vspace{-0.5cm}
\caption{Frequency response curves of the ropeway: vertical amplitude at the second gondola of the reduced model at $\bigO(3)$ for four forcing amplitudes applied at the first gondola, and full-order periodic orbits at $F=500$ and $650$\,N (dots).}
\label{fig_ropeway_FRC}
\end{figure}

For harmonic forcing $F\cos(\Omega t)$ applied vertically at the first gondola, \fig{fig_ropeway_FRC} shows frequency response curves of the vertical response at the second gondola for $F=150$, $350$, $500$, and $650$\,N.
The response softens, and the peak shifts from $0.4886$\,Hz at $F=150$\,N to $0.4750$\,Hz at $F=650$\,N, where the curve steepens toward a fold.
Full-order periodic orbits at $85\,\%$ of the peak amplitude of the two largest curves, obtained by Newton iteration on the one-period return map of the master coordinate rather than by transient settling (the $0.55\,\%$ modal damping ratio makes the free transient far too slow to settle in practice), agree with the reduced model to $-5.6\,\%$ at $F=500$\,N and $+5.1\,\%$ at $F=650$\,N.
These deviations are consistent with the small frequency bias of the reduced dynamics observed in the ring-down, since near the fold, where $\mathrm{d}\Omega/\mathrm{d}\rho\to0$, a small frequency error translates into a large amplitude difference at fixed frequency.
The frequencies of the full-order orbits at fixed amplitude agree with the reduced model to within one percent.

The SSM of order~$3$ is computed in about $210$\,s and a reduced frequency response curve in about $15$\,s, whereas one full-order periodic orbit by shooting takes $10$ to $20$\,min and a full-order free decay of $8$\,s takes about $9$\,min.

%%%%%%%%%%%%%%%%%%%%%%%%%%%%%%%%%%%%%%%%%%%%%%%%%%%%%%%%%%%%%%%%%%%%%%%%%%%%%%%%
\section{Conclusions}
\label{sec_Conclusions}
%%%%%%%%%%%%%%%%%%%%%%%%%%%%%%%%%%%%%%%%%%%%%%%%%%%%%%%%%%%%%%%%%%%%%%%%%%%%%%%%
%
This work has developed a framework for the SSM-based reduction of constrained multibody systems on matrix Lie groups, covering assembled models of rigid bodies, nodal-based floating frame of reference bodies, and geometrically exact beam finite elements connected by joints that permit relative motion.
The framework rests on two ingredients.
First, a polynomial computational representation describes each group element redundantly by its components, i.e., by the entries of the rotation matrix or by the Euler parameters together with the position, and enforces group membership through configuration-space constraints.
Rotation matrix components or Euler parameters render the equations of motion of rigid and floating frame of reference bodies polynomial, and an Euler--Rodrigues representation with lifted auxiliary variables and minimax polynomials removes the logarithmic map from the $SE(3)$ beam element while retaining its screw-motion interpolation.
%The geometrically exact formulation itself remains unaltered.
Second, the parameterization method has been reformulated directly in terms of the residual of the assembled equations of motion. 
This enables the computation of SSMs for autonomous and periodically forced systems without a descriptor reformulation, decomposition of the nonlinear forces into multilinear maps, or recasting the nonlinearities into quadratic form.
 
The analysis of constraint preservation has shown that the reduced model inherits the algebraic structure of the residual.
Constraints that are linear in the computational variables, which include lower-pair joints in the rotation-matrix representation, are satisfied exactly at every truncation order, whereas nonlinear constraints, among them the configuration-space constraints, hold up to the truncation order.
The rotation-preserving parameterization in the tangent space, followed by an explicit correction along the kinematic tree, satisfies both group membership and the joint constraints to machine precision, while leaving the Taylor coefficients and the reduced dynamics unchanged.
Constraint preservation is therefore a property of the chosen residual representation and of the evaluation of the parameterization, rather than of the reduction procedure itself.
 
On the implementation side, truncated Taylor arithmetic delivers all coefficients of a given order from a single residual evaluation.
For codes that evaluate their residual only in real arithmetic, a non-intrusive coupling reconstructs the residual as an explicit polynomial from real-valued evaluations.
 
The numerical examples cover rigid and flexible pendulum systems with and
without internal resonances, a floating satellite with hinged solar-array
wings, and a ropeway system with a track rope discretized by geometrically exact beam elements and coupled to rigid-body gondolas. 
The examples demonstrate that the proposed approach yields accurate low-dimensional reduced-order models for constrained multibody systems on matrix Lie groups. 
A complete forced response curve of the reduced model is obtained 20 to 80 times faster than a single full-order steady state, so that the speed-up for an entire curve grows with the number of frequencies at which the full model would otherwise have to be evaluated.
The offline cost of the SSM computation, between less than a second and a few minutes in the examples considered, lies below that of one full-order steady state and is incurred only once per model and truncation order, since the reduced model is reused for all forcing amplitudes and frequencies.
The flexible pendulum-slider has further demonstrated that a model defined entirely in the external multibody code Exudyn is reduced through real-valued residual evaluations alone.
This coupling is subject to restrictions: the residual must be evaluable at arbitrary states in a neighborhood of the equilibrium, and the nonlinear terms that are not recovered exactly are represented by a polynomial reconstruction, whose accuracy bounds that of the reduced model.
 
Joints that close kinematic loops are not covered by the correction map, since they would require a simultaneous correction of all bodies in the loop.
Future work will address the constrained evaluation of closed kinematic loops and the application of the framework to larger flexible multibody models.

%
%% To print the credit authorship contribution details
\printcredits

%% Loading bibliography style file
%

\appendix

%%%%%%%%%%%%%%%%%%%%%%%%%%%%%%%%%%%%%%%%%%%%%%%%%%%%%%%%%%%%%%%%%%%%%%%%%%%%%%%%
\section{Linear formulation of lower-pair joints}
\label{sect_LinearLowerPairJoints}
The preservation results established in Section~4.1 apply to all residual equations that are linear in the computational variables. 
This appendix demonstrates that a broad class of lower-pair joints admits such linear constraint formulations.
To this end, two elementary linear constraints are introduced, from which several standard lower-pair joints can be constructed as particular cases.

\subsection{Elementary linear constraints}
Two elementary linear constraints are introduced in this section.
The first constraint enforces the coincidence of two attachment points, whereas the second constraint enforces the coincidence of two body-fixed directions.

\paragraph{Point-coincidence constraint.}

Let $\uv_a,\uv_b\in\RealNumbers^3$ denote constant vectors from the body-fixed reference frames of bodies $a$ and $b$ to the joint location. 
Furthermore, let $\yv_a,\yv_b\in\RealNumbers^3$ and $\RotMat_a,\RotMat_b\in SO(3)$ denote the positions and orientations of the body-fixed reference frames with respect to the global inertial frame. 
The computational configuration variables hence read
\be\label{eq_ConfigVariablesLinearJoints}
\qVec = \begin{bmatrix}
\yv_a\tp & \vec{\RotMat_{a}}\tp & \yv_b\tp & \vec{\RotMat_{b}}\tp
\end{bmatrix}\tp \eqDot
\ee
The coincidence of the joint attachment points is formulated as
\begin{equation}
\phib_p = \yv_a+\RotMat_a\uv_a - \yv_b-\RotMat_b\uv_b = \zero \eqDot
\label{eq:PointCoincidence}
\end{equation}

Since $\uv_a$ and $\uv_b$ are constant body-fixed vectors, every component of $\phib_p$ is linear in the computational configuration variables $\qVec$.
The corresponding constraint gradient matrix is
\begin{equation}
\Bm_p
=
\frac{\partial\phib_p}{\partial\qVec}
=
\begin{bmatrix}
\identity_3 &
\uv_a^T\otimes\identity_3 &
-\identity_3 &
-\uv_b^T\otimes\identity_3
\end{bmatrix} \eqComma
\label{eq:PointConstraintGradient}
\end{equation}
where $\identity_3$ denotes the $3\times 3$ identity matrix.

\paragraph{Vector-coincidence constraint.}
Let $\sv_a, \sv_b \in \RealNumbers^3$ denote constant body-fixed unit vectors. 
Coincidence of their spatial directions is formulated as
\begin{equation}
\phib_v = \RotMat_a\sv_a - \RotMat_b\sv_b = \zero \eqDot
\label{eq:VectorCoincidence}
\end{equation}
Again, every component of $\phib_v$ is linear in the
computational configuration variables
\eq{eq_ConfigVariablesLinearJoints}. The corresponding
constraint gradient matrix is
\begin{equation}
\Bm_v
=
\frac{\partial\phib_v}{\partial\qVec}
=
\begin{bmatrix}
\zero &
\sv_a\tp\otimes\identity_3 &
\zero &
-\sv_b\tp\otimes\identity_3
\end{bmatrix},
\label{eq:VectorConstraintGradient}
\end{equation}

Since $\sv_a$ and $\sv_b$ are constant unit vectors, \eq{eq:VectorCoincidence} comprises three scalar equations although only two are linearly independent.
In practical joint formulations, the redundant equation is removed by projection onto two independent directions. 
Alternatively, the redundant formulation may be retained. Although the corresponding constraint gradient matrix is rank deficient, the generalized constraint forces depend only on its column space. 
Consequently, redundant linear constraint formulations are mechanically equivalent to their projected counterparts while remaining compatible with the SSM reduction developed in this paper.

\subsection{Generic linear joint formulation}
The translational and orientational constraint equations required for a particular lower-pair joint are selected by the constant matrices $\Sm_p$ and $\Sm_v$, respectively. 
In particular, $\Sm_v$ may be chosen to project the vector-coincidence constraint onto two independent directions, thereby eliminating the redundant scalar equation while preserving linearity.
\begin{equation}
\phib_{g}
=
\begin{bmatrix}
\Sm_p\phib_p\\
\Sm_v\phib_v
\end{bmatrix}
= \zero \eqDot
\label{eq:GenericJointConstraint}
\end{equation}
The corresponding joint constraint gradient matrix is
\begin{equation}
\Bm_{\mathrm{jnt}}
=
\frac{\partial\phib_{g}}{\partial\qv}
=
\begin{bmatrix}
\Sm_p\Bm_{p}\\
\Sm_v\Bm_{v}
\end{bmatrix} \eqDot
\label{eq:GenericJointGradient}
\end{equation}

Since both elementary constraints are linear in the
computational configuration variables \eq{eq_ConfigVariablesLinearJoints} and the selection matrices are constant, every constraint generated by \eq{eq:GenericJointConstraint} is linear.
Consequently, the preservation result established in \sect{sec_LinearConstraints} applies directly to all joint formulations that can be represented by \eq{eq:GenericJointConstraint}.

\subsection{Classical lower-pair joints}
The elementary point- and vector-coincidence constraints introduced above provide the building blocks for constructing a broad class of lower-pair joints. 
By appropriate selection of translational and orientational constraint equations, many classical lower-pair joints, including spherical, revolute, cylindrical, and prismatic joints, admit formulations whose residual equations are linear with respect to the computational variables.
Universal joints may be represented by introducing an auxiliary body connected to the adjacent bodies through two revolute joints with orthogonal axes. Since both revolute joints admit linear formulations within the proposed framework, the resulting universal joint constraint equations consist exclusively of linear constraint equations.
Consequently, the preservation result established in \sect{sec_LinearConstraints} applies directly to these joint formulations.

\paragraph{Example (Revolute joint)}
A revolute joint is obtained by combining one point-coincidence constraint with a projected vector-coincidence constraint. 
This corresponds to
\be
\Sm_p = \identity_3 \eqComma
\ee
and
\be
\Sm_v=
\begin{bmatrix}
\bar{\mathbf u}_1 \tp \\
\bar{\mathbf u}_2 \tp
\end{bmatrix} \eqComma
\ee
where $\bar{\mathbf u}_1,\bar{\mathbf u}_2\in\mathbb R^3$ are two constant unit vectors satisfying

\be
\bar{\mathbf{u}}_i\tp\bar{\mathbf{s}} = 0 \eqComma
\qquad i=1,2,
\ee
with $\bar{\mathbf s}$ denoting the joint axis in the
reference configuration.
The resulting joint consists of three translational and two orientational constraint equations, yielding the five independent constraints of a revolute joint.

The present work focuses on model reduction of multibody systems formulated as differential-algebraic equations on matrix Lie groups rather than on the development of kinematic joint formulations. 
Therefore, further explicit derivations of the individual joint equations are omitted for brevity and are left to the reader.

\section{Application to Exudyn}
\label{app:exudyn}
%\label{sec:exudyn-application}
%%%%%%%%%%%%%%%%%%%%%%%%%%%%%%%%%%%%%%%%%%%%%%%%%%%%%%%%%%%%%%%%%%%%%%%%%%%%%%%%
%
Exudyn \cite{Gerstmayr2024} is an open-source, C++-based Python library for efficient simulation of flexible multibody dynamics systems.
The broad usability of Python allows an Exudyn multibody model to be coupled with external Python packages.

The multibody formulation used by Exudyn is mainly based on redundant coordinates and hence fits into the considered modeling approach of constrained multibody systems introduced in \sect{sect_LieGroupEOM}.
Although Lie group time integration \cite{Holzinger2021, Holzinger2025} is available, Exudyn does not generally employ the first-order representation of the equations of motion as presented in \sects{sect_LieGroupEOM}{sec:csp-constraint}.
Instead, it is based on a second-order formulation of the equations of motion with joint and configurations space constraints imposed on the position level using Lagrange multipliers. 
The Lie group formulation employed for rigid body motion is $\RealNumbers^3 \times S^3$.
Consequently, the form of the assembled residual of Exudyn differs from residual representation \eq{eq_CompleteResidual}.
However, the proposed residual formulation of the parameterization is not affected by the underlying formulation of the governing equations. 

The section is structured as follows. 
The first part addresses the residual structure of Exudyn. 
The second part describes the extraction and processing of the relevant residuals required to compute an SSM reduced-order model using the proposed residual formulation.

\subsection{Assembled residual of Exudyn}\label{sect_AssembledExudynResidual}
Exudyn formulates the equations of motion directly in the computational configuration variables $\qVec$ of \sect{sect_CompRep} and their time derivatives, which we denote for notational convenience by
\be
\qVecExuDot \equiv \dot{\qVec} \eqComma \qquad \qVecExuDDot \equiv \ddot{\qVec} \eqDot
\ee
They are related to the velocity $\vv$ and the acceleration $\dot\vv$ by
\begin{align}
\vv &= \kinRecMatInv(\qVec)\,\qVecExuDot \eqComma \label{eq_ComputationalVelocityExtraction} \\
\dot{\vv} &= \kinRecMatInv(\qVec)\,\qVecExuDDot + \dot{\kinRecMatInv}(\qVec,\qVecExuDot)\,\qVecExuDot \eqDot \label{eq_AccelerationExtraction}
\end{align}
Left multiplication of \eq{eq_MechanicalResidual} by $\kinRecMat(\qVec)$ (\eq{eq_CompRecMat}) then yields the mechanical residual of Exudyn on $\RealNumbers^3\times\SThree$ for holonomic constraints,
\be
\widehat{\MRes}(\qVec,\qVecExuDot,\qVecExuDDot,\lagMulb,\auxVarb) = \widehat{\Mb}(\qVec)\qVecExuDDot + \widehat{\Bm}\jnt\tp(\qVec)\lagMulb - \widehat{\gb}(\qVec, \qVecExuDot, \auxVarb, \Omega t, \epsilon) = \zero \eqComma
\label{eq_MechResExu}
\ee
where\footnote{Note, for any Lie group whose rotational part is represented by Euler parameters and whose translational part is represented in Cartesian coordinates the term $\dot{\Pm}(\qVec,\qVecExuDot)\qVecExuDot=\zero$.
Therefore, the inertia term $\Qm(\qVec)\Mm(\qVec)\dot{\Pm}(\qVec,\qVecExuDot)\dot{\qVec}$ does not show up in \eq{eq_MechResExu}.}
\begin{align}
\widehat{\Mb}(\qVec) &:= \kinRecMat(\qVec) \Mm(\qVec) \kinRecMatInv(\qVec) \in \RealNumbers^{\dimq \times \dimq} \eqComma \\
\widehat{\gb}(\qVec, \qVecExuDot, \auxVarb, \Omega t, \epsilon) &:= \kinRecMat(\qVec) \left[\gb\ine(\qVec, \kinRecMatInv(\qVec)\qVecExuDot) + \gb\int(\qVec, \kinRecMatInv(\qVec)\qVecExuDot) + \gb\aux(\qVec,\kinRecMatInv(\qVec)\qVecExuDot,\auxVarb,t) \right. \nonumber \\
& - \left.\epsilon\,\gb_\textrm{ext}(\qVec, \kinRecMatInv(\qVec)\qVecExuDot, \Omega t) \right] \in \RealNumbers^{\dimq}, \\
\widehat{\Bm}\jnt\tp(\qVec) :&= \kinRecMat(\qVec) \Bm\jnt\tp(\qVec) \in \RealNumbers^{\dimq \times \dimJnt} \eqComma 
\end{align}

The matrix $\kinRecMatInv(\qVec)$ in the above equations is computed as
\be\label{eq_CompRecMatInv}
\kinRecMatInv(\qVec) = \left(\Em\tp\Em\right)^{-1}\Em\tp (e \otimes q^{-1}) \in \RealNumbers^{\dimG \times \dimq } \eqComma
\ee
and satisfies
\begin{equation}
\Pm(\qVec)\kinRecMat(\qVec)=\Im \in \RealNumbers^{\dimG \times \dimG} \eqDot
\end{equation}

Exudyn does not distinguish joint and configuration space constraints.
Instead, both constraints are combined into the a single constraint residual
\begin{equation}
\widehat{\phib}(\qVec)
=
\begin{bmatrix}
\phib\csp(\qVec)\\
\phib\jnt(\qVec)
\end{bmatrix}
=\zero \eqDot
\end{equation}
The residual assembled form the respective residual blocks of Exudyn is finally obtained as
\begin{equation}
\widehat{\CRes}(\qVec,\qVecExuDot,\qVecExuDDot,\lagMulb)
=
\begin{bmatrix}
\widehat{\MRes}(\qVec,\qVecExuDot,\qVecExuDDot,\lagMulb, \auxVarb) \\
\dot{\qVec} - \qVecExuDot \\
\dot\auxVarb - \ab(\qVec,\kinRecMatInv(\qVec)\qVecExuDot,\auxVarb,t) \\
\widehat{\phib}(\qVec)
\end{bmatrix}
=\Null \eqDot
\label{eq_CompleteResidualSecondOrder}
\end{equation}

Note, in \sect{sect_LieGroupEOM}'s own general formulation, $\qVec$ and $\vv$ are the independent state components, related only through the kinematic residual itself. 
In contrast to the assembled residual \eq{eq_CompleteResidual}, in which $\vv$ is an independent state component, Exudyn uses $\qVecExuDot$ as the velocity-level state.
The velocity $\vv$ therefore does not appear as an argument of $\widehat{\MRes}$ and enters only through \eq{eq_ComputationalVelocityExtraction}.
For holonomic constraints, the recursion of
\sect{sec:computing-parameterizations} can therefore be carried out
directly on $(\qVec,\dot\qVec,\lagMulb)$, without ever introducing a
separate kinematic residual row.

While the residual \eq{eq_CompleteResidualSecondOrder} may appear simpler than the residual of the first-order formulation \eq{eq_CompleteResidual}, the mechanical residual becomes more nonlinear due to \eqs{eq_ComputationalVelocityExtraction}{eq_AccelerationExtraction}. 
For rigid multibody systems, for example, the polynomial order increases from two in \eq{eq_MechanicalResidual} to at least three in the Exudyn formulation \eq{eq_MechResExu}.

%%%%%%%%%%%%%%%%%%%%%%%%%%%%%%%%%%%%%%%%%%%%%%%%%%%%%%%%%%%%%%%%%%%%%%%%%%%%%%%%
\subsection{Interfacing Exudyn's residuals}
%%%%%%%%%%%%%%%%%%%%%%%%%%%%%%%%%%%%%%%%%%%%%%%%%%%%%%%%%%%%%%%%%%%%%%%%%%%%%%%%
%
The residual is obtained from Exudyn's internal assembly routine \texttt{ComputeNewtonResidual}, which its generalized-$\alpha$ integrator \cite{Arnold2007} calls in every Newton iteration to assemble the mechanical and constraint rows.
Here, the routine is called directly at arbitrary trial states, independently of the time integration.

\texttt{ComputeNewtonResidual} returns the mechanical and the constraint rows in a vector of length $\dimq+\dimJnt$, the constraint rows unscaled and the mechanical rows scaled.
Denoting the mechanical rows by
\be
\overline{\MRes}(\qVec,\qVecExuDot,\qVecExuDDot,\lagMulb,\auxVarb) :=
\texttt{ComputeNewtonResidual}(\ldots)[1{:}\dimq] \eqComma
\ee
the mechanical residual \eq{eq_MechResExu} follows as
\be
\widehat{\MRes} = \beta'\,\overline{\MRes} \eqComma
\qquad
\beta' = \frac{1}{\beta h^2} \eqComma
\qquad
\beta = \tfrac{1}{4}\bigl(\gamma+\tfrac{1}{2}\bigr)^2 \eqComma
\qquad
\gamma = \tfrac{1}{2}+\alpha_f-\alpha_m \eqComma
\ee
where $h$ denotes the step size of the generalized-$\alpha$ method and $\alpha_m=(2\rho_\infty-1)/(\rho_\infty+1)$ and $\alpha_f=\rho_\infty/(\rho_\infty+1)$ follow from the spectral radius $\rho_\infty$ \cite{Arnold2007}.
Since no time step is performed, $h$ is merely the step size set in the solver and cancels in this rescaling.

Since the constraint rows are unscaled, they can be read off directly as
\be
\widehat\phib(\qVec) = \texttt{ComputeNewtonResidual}(\ldots)[\dimq{+}1{:}\dimq{+}\dimJnt] \eqDot
\ee
Together, $\widehat\MRes$ and $\widehat\phib$ assemble directly into
the assembled residual $\widehat{\CRes}$ of
\eq{eq_CompleteResidualSecondOrder}.

The equilibrium point is obtained first, from a real-valued Newton iteration on $\hat{\mathcal F}$. 
The mass matrix $\hat\Mm$ and the constraint Jacobian $\hat\Bm\tp\jnt$ are taken from Exudyn's assembled system matrices (\appx{sect_AssembledExudynResidual}), while the dependence of $\hat\Mm$, $\hat\Bm\tp\jnt$, $\hat{\boldsymbol\Phi}$, and $\hat\gv$ on $\qVec$ is recovered by polynomial fitting.

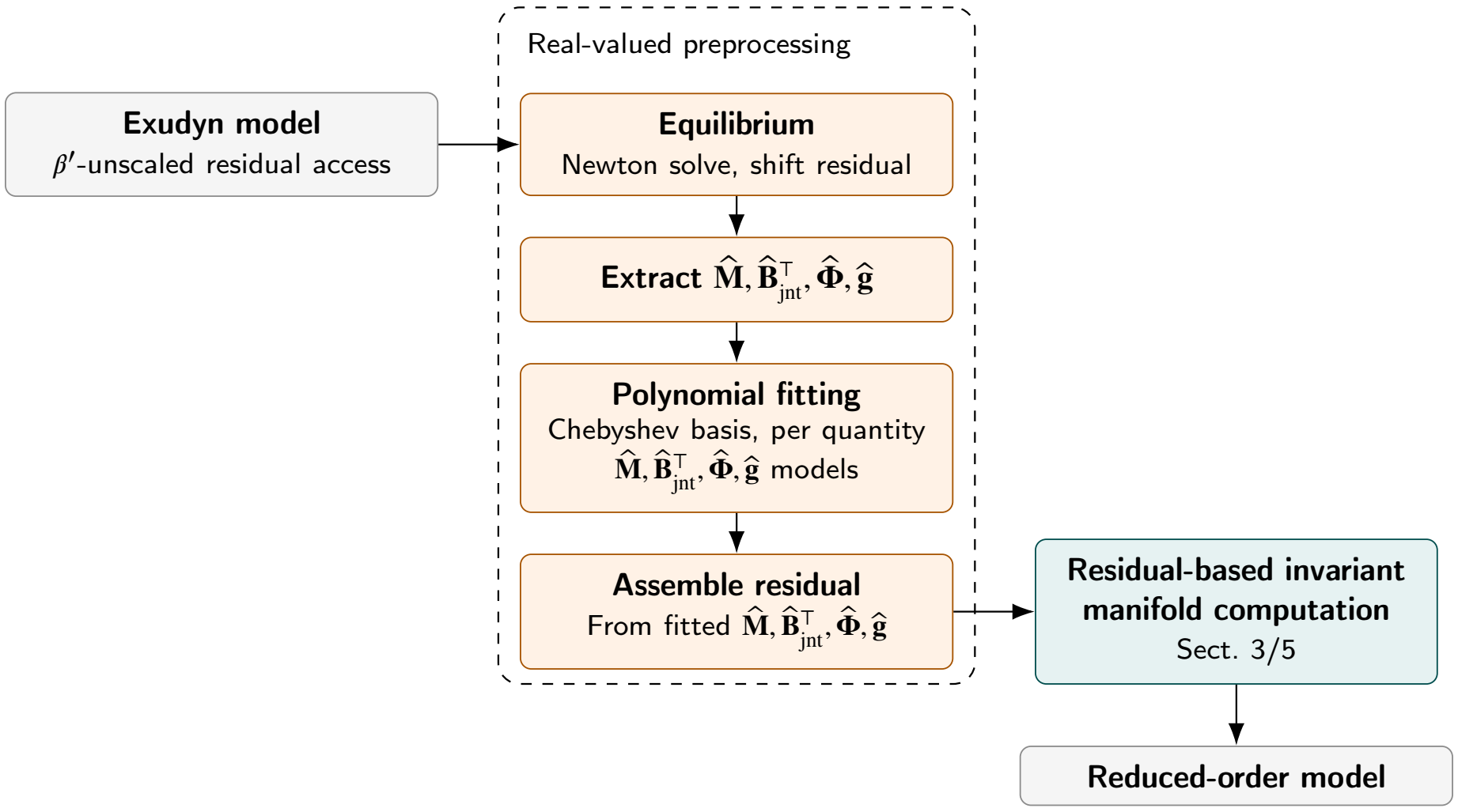
\begin{figure}
\centering
%\resizebox{\linewidth}{!}{%
\scalebox{0.85}{
% exudyn_pipeline_fig_horizontal.tex -- Exudyn coupling pipeline, horizontal
% layout for a full-width (figure*) placement in a double-column document.
%
% Import into the main document with, e.g.:
%
%   \begin{figure*}[t]
%     \centering
%     \input{exudyn_pipeline_fig_horizontal.tex}
%     \caption{...}
%     \label{fig:exudyn-pipeline}
%   \end{figure*}
%
% Requires, in the MAIN document's own preamble:
%   \usepackage{tikz}
%   \usetikzlibrary{arrows.meta, positioning, fit, backgrounds}
%
% Uses the paper's own macros directly (\Mb, \Bm, \gb, \phib, \jnt, \tp,
% \aux, \qVec) -- these must already be defined in the main document;
% this file does not define them itself.
%
% Font sizes match exudyn_pipeline_fig.tex (the vertical version) --
% \small/\footnotesize, not \scriptsize/\tiny -- since the smaller sizes
% triggered a MiKTeX on-the-fly bitmap font generation failure.

\begin{tikzpicture}[
  node distance=6mm,
  every node/.style={font=\small},
  box/.style={
    rectangle, rounded corners=3pt, draw, align=center,
    text width=3.9cm, inner sep=5pt, line width=0.4pt
  },
  graybox/.style={box, fill=black!4, draw=black!45},
  coralbox/.style={box, fill=orange!10, draw=orange!65!black},
  tealbox/.style={box, text width=3.6cm, fill=teal!10, draw=teal!60!black},
  arr/.style={-{Latex[length=2.2mm, width=1.6mm]}, line width=0.5pt},
  lbl/.style={font=\footnotesize, align=center},
]

% -- real-valued preprocessing (inside the dashed container), stacked ------
\node[coralbox] (equilibrium)
  {\textbf{Equilibrium} \\ \footnotesize Newton solve, shift residual};

\node[coralbox, below=4mm of equilibrium] (extract)
  {\textbf{Extract} $\widehat{\Mb},\widehat{\Bm}\jnt\tp,\widehat{\phib},\widehat{\gb}$ \\
   %\footnotesize Direct Exudyn API calls \\
   %\footnotesize $\widehat{\gb}\aux$ not extracted
   };

\node[coralbox, below=4mm of extract] (fitting)
  {\textbf{Polynomial fitting} \\
   \footnotesize Chebyshev basis, per quantity \\
   \footnotesize $\widehat{\Mb},\widehat{\Bm}\jnt\tp,\widehat{\phib},\widehat{\gb}$ models};

\node[coralbox, below=4mm of fitting] (assemble)
  {\textbf{Assemble residual} \\
   \footnotesize From fitted $\widehat{\Mb},\widehat{\Bm}\jnt\tp,\widehat{\phib},\widehat{\gb}$};

% -- left of container: Exudyn model, vertically aligned with Equilibrium --
\node[graybox, left=8mm of equilibrium] (exudyn)
  {\textbf{Exudyn model} \\ \footnotesize $\beta'$-unscaled residual access};

\node[coordinate, above=7mm of equilibrium.north] (topspacer) {};

\begin{pgfonlayer}{background}
\node[draw, dashed, rounded corners=6pt, line width=0.5pt,
      fit=(topspacer)(equilibrium)(extract)(fitting)(assemble),
      inner xsep=6pt, inner ysep=4pt]
      (container) {};
\end{pgfonlayer}
\node[lbl, anchor=north west, xshift=5pt, yshift=-4pt] at (container.north west)
  {Real-valued preprocessing};

% -- right of container: generic recursion, vertically aligned with Assemble
\node[tealbox, right=8mm of assemble] (recursion)
  {\textbf{Residual-based invariant manifold computation} \\ \footnotesize Sect.~3/5};

% output now BELOW recursion, not to its right -- shortens the figure
% horizontally at the cost of a bit more height on the right side
\node[graybox, below=6mm of recursion] (output)
  {\textbf{Reduced-order model}};

% -- arrows -----------------------------------------------------------------
\draw[arr] (exudyn.east) -- (equilibrium.west);
\draw[arr] (equilibrium) -- (extract);
\draw[arr] (extract) -- (fitting);
\draw[arr] (fitting) -- (assemble);
\draw[arr] (assemble.east) -- (recursion.west);
\draw[arr] (recursion.south) -- (output.north);

\end{tikzpicture}}
\caption{Procedure for computing an SSM reduced-order model from an Exudyn multibody model.}
\label{fig_exudynPipeline}
\end{figure}

The complete procedure is summarized in \fig{fig_exudynPipeline}.

%%%%%%%%%%%%%%%%%%%%%%%%%%%%%%%%%%%%%%%%%%%%%%%%%%%%%%%%%%%%%%%%%%%%%%%%%%%%%%%%
\section{Computation of the minimax polynomial coefficients}
\label{app_MinimaxCoefficients}
%%%%%%%%%%%%%%%%%%%%%%%%%%%%%%%%%%%%%%%%%%%%%%%%%%%%%%%%%%%%%%%%%%%%%%%%%%%%%%%%
%
This appendix summarizes the construction of the minimax polynomial
approximations used for the scalar functions appearing in
\sect{sect_PolyInterpolation}.  The coefficients are generated by the Remez
exchange algorithm using arbitrary-precision arithmetic and are subsequently
verified against the exact transcendental functions.

%%%%%%%%%%%%%%%%%%%%%%%%%%%%%%%%%%%%%%%%%%%%%%%%%%%%%%%%%%%%%%%%%%%%%%%%%%%%%%%%
\subsection{Target functions}
%%%%%%%%%%%%%%%%%%%%%%%%%%%%%%%%%%%%%%%%%%%%%%%%%%%%%%%%%%%%%%%%%%%%%%%%%%%%%%%%
%
The polynomial beam formulation depends on the analytic scalar functions
\be
S(a)=\frac{\sin\sqrt a}{\sqrt a},
\qquad
C(a)=\cos\sqrt a,
\qquad
F(4a)=\frac{1-\frac{\sqrt{4a}}{2}\cot\!\left(\frac{\sqrt{4a}}{2}\right)}{4a},
\ee
together with several auxiliary coefficient functions appearing in the element
matrices.

%%%%%%%%%%%%%%%%%%%%%%%%%%%%%%%%%%%%%%%%%%%%%%%%%%%%%%%%%%%%%%%%%%%%%%%%%%%%%%%%
\subsection{Minimax approximation}
%%%%%%%%%%%%%%%%%%%%%%%%%%%%%%%%%%%%%%%%%%%%%%%%%%%%%%%%%%%%%%%%%%%%%%%%%%%%%%%%
%
Each function is approximated by a polynomial
\be
P_n(a)=\sum_{i=0}^{n}c_i a^i,
\ee
chosen to minimize the maximum approximation error over the prescribed fitting
interval.  The coefficients are obtained with the Remez exchange algorithm,
which exploits Chebyshev's equioscillation theorem: among all degree-$n$
polynomials, the unique minimax polynomial is characterized by an error curve
that attains equal alternating extrema at $n+2$ points. Starting from
Chebyshev reference points, the algorithm repeatedly solves the interpolation
system enforcing alternating errors and updates the extremal points until
convergence.

%%%%%%%%%%%%%%%%%%%%%%%%%%%%%%%%%%%%%%%%%%%%%%%%%%%%%%%%%%%%%%%%%%%%%%%%%%%%%%%%
\subsection{Choice of approximation norm}
%%%%%%%%%%%%%%%%%%%%%%%%%%%%%%%%%%%%%%%%%%%%%%%%%%%%%%%%%%%%%%%%%%%%%%%%%%%%%%%%
%
The approximation norm is selected individually for each function and is listed in Tab.~\ref{tab_MinimaxOverview}.
Functions are fitted either in the absolute norm,
\be
\max_a|P_n(a)-f(a)| \eqComma
\ee
or in the relative norm,
\be
\max_a\frac{|P_n(a)-f(a)|}{|f(a)|} \eqDot
\ee
The function $C(a)=\cos\sqrt a$ vanishes inside the fitting interval and is therefore fitted in the absolute norm, which avoids the singular weighting of a relative-error fit.
The remaining functions are bounded and of order one on the fitting interval, so that both norms differ only by a constant factor and either yields coefficient tables accurate to double-precision round-off.
The relative norm is used for the functions of the strain tangent, \tab{tab_PB_functions}, and for $\rho$, whose relative error enters the element forces directly.

%%%%%%%%%%%%%%%%%%%%%%%%%%%%%%%%%%%%%%%%%%%%%%%%%%%%%%%%%%%%%%%%%%%%%%%%%%%%%%%%
\subsection{Numerical construction}
%%%%%%%%%%%%%%%%%%%%%%%%%%%%%%%%%%%%%%%%%%%%%%%%%%%%%%%%%%%%%%%%%%%%%%%%%%%%%%%%
%
The coefficients are generated with arbitrary-precision arithmetic using the
following procedure.

\begin{enumerate}
\item Select the fitting interval for the target function.
\item Evaluate the exact function with extended precision.
\item Run the Remez exchange algorithm to obtain the degree-$n$ minimax
      polynomial.
\item Verify that the resulting error satisfies the expected equioscillation
      property.
\item Round the coefficients to double precision and verify that the rounded
      tables still remain below double-precision round-off over the entire
      fitting interval.
\end{enumerate}

The implementation performs the fitting in arbitrary-precision arithmetic with 60 decimal digits, and with 260 digits for the functions with a removable singularity at $a=0$, so that the cancellation in their closed forms near $a=0$ is absorbed by the working precision. 
Only the final coefficients are rounded to double precision.

%%%%%%%%%%%%%%%%%%%%%%%%%%%%%%%%%%%%%%%%%%%%%%%%%%%%%%%%%%%%%%%%%%%%%%%%%%%%%%%%
\subsection{Polynomial evaluation}
%%%%%%%%%%%%%%%%%%%%%%%%%%%%%%%%%%%%%%%%%%%%%%%%%%%%%%%%%%%%%%%%%%%%%%%%%%%%%%%%

The resulting coefficients are evaluated with Horner's scheme,
\be
P_n(a)
=
c_0+a(c_1+a(c_2+\cdots+a c_n)) \eqComma
\ee
requiring only multiply-add operations and preserving the polynomial structure
of the beam formulation.

%%%%%%%%%%%%%%%%%%%%%%%%%%%%%%%%%%%%%%%%%%%%%%%%%%%%%%%%%%%%%%%%%%%%%%%%%%%%%%%%
\subsection{Approximation degrees}
%%%%%%%%%%%%%%%%%%%%%%%%%%%%%%%%%%%%%%%%%%%%%%%%%%%%%%%%%%%%%%%%%%%%%%%%%%%%%%%%
%
\tab{tab_MinimaxOverview} summarizes the polynomial degrees used in this work.
The polynomial for $\hat a(q_0)=\arccos^2(q_0)$ is used only to initialize the lifted coordinate $a$ from the reference configuration and does not enter the residual.
\begin{table}[ht]
\centering
\caption{Polynomial degrees of the minimax approximations.}
\label{tab_MinimaxOverview}
\begin{tabular}{lcc}
\toprule
Function & Norm & Degree\\
\midrule
$C(a)$ & absolute & 8\\
$S(a)$ & absolute & 8\\
$F(4a)$ & absolute & 13\\
$\rho(a)=1/S(a)$ & relative & 14\\
$\hat a(q_0)=\arccos^2(q_0)$ & absolute & 19\\
$f_1(a)$ & relative & 9\\
$f_2(a)$ & relative & 9\\
$f_3(a)$ & relative & 9\\
$f_4(a)$ & relative & 9\\
$k_2(a)$ & relative & 15\\
\bottomrule
\end{tabular}
\end{table}

%%%%%%%%%%%%%%%%%%%%%%%%%%%%%%%%%%%%%%%%%%%%%%%%%%%%%%%%%%%%%%%%%%%%%%%%%%%%%%%%
\subsection{Coefficient generation and verification}
%%%%%%%%%%%%%%%%%%%%%%%%%%%%%%%%%%%%%%%%%%%%%%%%%%%%%%%%%%%%%%%%%%%%%%%%%%%%%%%%
%
The complete implementation used to generate and verify every coefficient table
is reproduced in the accompanying supplementary material. The script implements
the Remez exchange algorithm, automatically searches for admissible polynomial
degrees when requested, and verifies that the rounded double-precision
coefficients attain errors below machine precision over the prescribed fitting
interval.

%%%%%%%%%%%%%%%%%%%%%%%%%%%%%%%%%%%%%%%%%%%%%%%%%%%%%%%%%%%%%%%%%%%%%%%%%%%%%%%%
\section{Implementation of the polynomial geometrically exact beam element}
\label{app_PolyBeam}
This appendix collects the equations required to implement the two-node polynomial beam element of \sect{sect_PolyInterpolation}.
The element is equivalent to the $SE(3)$ element of \cite{Sonneville2014}.
Throughout, $\widetilde{\mathbf{a}}\mathbf{b}=\mathbf{a}\times\mathbf{b}$, $[\mathbf{a},\mathbf{b}]:=\widetilde{\mathbf{a}}\widetilde{\mathbf{b}}+\widetilde{\mathbf{b}}\widetilde{\mathbf{a}}$, and vectors in $\RealNumbers^6$ are ordered as translational part $(\cdot)_U$ followed by rotational part $(\cdot)_\Omega$.
The adjoint operator of $\mathbf{y}=[\mathbf{y}_U\tp \, \mathbf{y}_\Omega\tp]\tp\in\RealNumbers^6$ is
\begin{equation}
\operatorname{ad}(\mathbf{y})=\begin{bmatrix}\vecToSkew{\mathbf{y}}_\Omega & \vecToSkew{\mathbf{y}}_U\\ \zero & \vecToSkew{\mathbf{y}}_\Omega\end{bmatrix} .
\label{eq_PB_ad}
\end{equation}
All scalar functions, $S$, $C$, $F$, and $\rho=1/S$ of \eq{eq_SFdef} and $f_1,\ldots,f_4$, $k_2$ of \tab{tab_PB_functions}, are evaluated from the minimax polynomials of \appx{app_MinimaxCoefficients}.
Their derivatives, denoted by a prime and taken with respect to the argument, follow by differentiating these polynomials coefficientwise.

%-------------------------------------------------------------------------------
\subsection{Element state and kinematic equations}
\label{app_PB_state}
An element connects the nodes $A$ and $B$ with the nodal unknowns and body-frame velocities of \eq{eq_nodalUnknowns}, whose kinematic reconstruction and unit constraint are given in \sect{sect_BeamNodalUnknowns}.
The nodal velocities are collected in $\vv^e=[\vv_A\tp\,\,\vv_B\tp]\tp\in\RealNumbers^{12}$.
The relative quantities $q_0$ and $\pV$ follow from \eq{eq_relativeErpFromNodes}, the scalars $\rho$ and $\mu$ from \eq{eq_bridgeScalars}, and the relative velocity from \eq{eq_relVelocity}, written compactly as $\Delta\vv=\vv_B-\Ad^{-1}\vv_A$ with the $6\times6$ matrix $\Ad^{-1}$ read off from \eq{eq_relVelocity}.
The lifted scalars $\mathbf{y}_e=(a,r,r_1,\ldots,r_{n_g})$ are additional configuration coordinates without inertia, one station reciprocal $r_k$ per quadrature point $u_k=s_k/L$, $k=1,\ldots,n_g$.
They are pinned by \eq{eq_liftedIdentities}, which contributes $2+n_g$ rows to the joint constraints and, as a generalized force, to the rows of the lifted coordinates themselves, see \sect{sect_PolyInterpolationMotion}.
The gradient of \eq{eq_liftedIdentities} with respect to the nodal rotations follows from the relative scalar of \eq{eq_relativeErpFromNodes} and the nodal kinematics \eq{eq_QmatER}.

%-------------------------------------------------------------------------------
\subsection{Chart tangent and conversion Jacobian}
\label{app_PB_chart}
The configuration $\psi(\pV)=(\RotMat(q_0,\pO),\Tm_\Omega\tp\pU)$ has the body-frame variation $\delta\mathbf{h}=\Tm(\pV)\,\delta\pV$ with
\begin{equation}
\Tm(\pV)=\begin{bmatrix}\Tm_\Omega & \Tm_+\\ \zero & \Tm_\Omega\end{bmatrix},\qquad
\Tm^{-1}(\pV)=\begin{bmatrix}\Tm_\Omega^{-1} & \Tm_-\\ \zero & \Tm_\Omega^{-1}\end{bmatrix},
\label{eq_PB_chartT}
\end{equation}
where $\Tm_\Omega$ is given by \eq{eq_Rtheta}, $\Tm_\Omega^{-1}=q_0\identity+\tfrac12\vecToSkew{\pO}$ is the transpose of the matrix in the third line of \eq{eq_relativeErpFromNodes}, and
\begin{equation}
\begin{aligned}
\Tm_-&=-\tfrac14\bigl(\pO\tp\pU\bigr)\,r\,\identity+\tfrac12\vecToSkew{\pU} ,\\
\Tm_+&=\tfrac14\bigl(\pO\tp\pU\bigr)\,r^3\bigl(\identity+\tfrac14\vecToSkew{\pO}^{\,2}\bigr)-\tfrac12\vecToSkew{\pU}+\tfrac14 r\,[\pO,\pU] .
\end{aligned}
\label{eq_PB_chartTpm}
\end{equation}
The Jacobian of the bridge relation \eq{eq_bridge}, $\mathbf{G}(\pV):=\partial\dv/\partial\pV$, has coinciding diagonal blocks and a division-free inverse,
\begin{equation}
\mathbf{G}=\begin{bmatrix}\Pbridge_U & \Pbridge_{U\Omega}\\ \zero & \Pbridge_U\end{bmatrix},\qquad
\mathbf{G}^{-1}=\begin{bmatrix}\Pbridge_U^{-1} & -\Pbridge_U^{-1}\Pbridge_{U\Omega}\Pbridge_U^{-1}\\ \zero & \Pbridge_U^{-1}\end{bmatrix},
\label{eq_PB_G}
\end{equation}
where $\Pbridge_U$ is given by \eq{eq_bridge}, its inverse $\Pbridge_U^{-1}$ in \sect{sect_PolyInterpolationMotion}, and
\begin{equation}
\begin{aligned}
\Pbridge_{U\Omega}&=\pU\,\boldsymbol{\rho}_{,\Omega}\tp+\bigl(\dv_\Omega\tp\pU\bigr)\,\dv_\Omega\,\boldsymbol{\mu}_{,\Omega}\tp+\mu\,\dv_\Omega\,\dv_U\tp+\mu\bigl(\dv_\Omega\tp\pU\bigr)\Pbridge_U,\\
\boldsymbol{\rho}_{,\Omega}&=-\tfrac12S'(a)\,\rho^2r\,\dv_\Omega,\qquad
\boldsymbol{\mu}_{,\Omega}=2F'(4a)\,r^2\,\dv_\Omega+\tfrac14F(4a)\,r^3\pO .
\end{aligned}
\label{eq_PB_XUOmega}
\end{equation}

%-------------------------------------------------------------------------------
\subsection{Strain tangent}
\label{app_PB_strainTangent}
The strain tangent $\boldsymbol{\mathcal{T}}(\dv)$ is defined by $\Delta\vv=\boldsymbol{\mathcal{T}}(\dv)\,\dot{\dv}$ and coincides with the tangent operator of $\exp_{SE(3)}$ in \cite{Sonneville2014}.
For any chart, $\Tm(\pV)\,\delta\pV$ and $\boldsymbol{\mathcal{T}}(\dv)\,\delta\dv$ are the same body-frame variation, so that
\begin{equation}
\boldsymbol{\mathcal{T}}(\dv)=\Tm(\pV)\,\mathbf{G}^{-1}(\pV),\qquad
\boldsymbol{\mathcal{T}}^{-1}(\dv)=\mathbf{G}(\pV)\,\Tm^{-1}(\pV) .
\label{eq_PB_chartIdentity}
\end{equation}
With $\dv=[\dv_U\tp \,\, \dv_\Omega\tp]\tp$, $a=\frac14\dv_\Omega\tp\dv_\Omega$ of \eq{eq_bridgeScalars} evaluated from $\dv$ rather than taken from the lifted coordinate, and $\varsigma:=\dv_\Omega\tp\dv_U$, both are evaluated directly as
\begin{equation}
\begin{aligned}
\boldsymbol{\mathcal{T}}&=\begin{bmatrix}\boldsymbol{\mathcal{T}}_\Omega & \boldsymbol{\mathcal{T}}_U\\ \zero & \boldsymbol{\mathcal{T}}_\Omega\end{bmatrix}, &
\boldsymbol{\mathcal{T}}_\Omega&=\identity+f_1\vecToSkew{\dv_\Omega}+f_2\vecToSkew{\dv_\Omega}^{\,2}, &
\boldsymbol{\mathcal{T}}_U&=f_1\vecToSkew{\dv_U}+f_2[\dv_U,\dv_\Omega]+f_3\varsigma\,\vecToSkew{\dv_\Omega}+f_4\varsigma\,\vecToSkew{\dv_\Omega}^{\,2},\\
\boldsymbol{\mathcal{T}}^{-1}&=\begin{bmatrix}\boldsymbol{\mathcal{T}}_\Omega^{-1} & \boldsymbol{\mathcal{T}}_U^{-}\\ \zero & \boldsymbol{\mathcal{T}}_\Omega^{-1}\end{bmatrix}, &
\boldsymbol{\mathcal{T}}_\Omega^{-1}&=\identity+\tfrac12\vecToSkew{\dv_\Omega}+k_1\vecToSkew{\dv_\Omega}^{\,2}, &
\boldsymbol{\mathcal{T}}_U^{-}&=\tfrac12\vecToSkew{\dv_U}+k_1[\dv_U,\dv_\Omega]+k_2\varsigma\,\vecToSkew{\dv_\Omega}^{\,2},
\end{aligned}
\label{eq_PB_strainTangent}
\end{equation}
with the scalar functions of \tab{tab_PB_functions}, all evaluated at $a$.
The closed forms in \tab{tab_PB_functions} are removable singularities at $a=0$ and serve only to generate the minimax polynomials; they are never evaluated.
$f_1,\ldots,f_4$ are entire, and $k_1$, $k_2$ are analytic for $\phi<2\pi$.

\begin{table}[t]
\centering
\caption{Scalar functions of the strain tangent \eq{eq_PB_strainTangent}, with $\phi=2\sqrt{a}$, $z=\phi^2=4a$, $\alpha=\sin\phi/\phi$, and $\beta=2(1-\cos\phi)/\phi^2$ \cite[Eqs.~(A.6), (A.7), (A.12), (A.14)]{Sonneville2014}.}
\label{tab_PB_functions}
\small
\begin{tabular}{llll}
\toprule
function & closed form & function & closed form\\
\midrule
$f_1(a)$ & $-\beta/2$ & $f_4(a)$ & $\bigl(\beta/2-3f_2\bigr)/z$\\
$f_2(a)$ & $(1-\alpha)/z$ & $k_1(a)$ & $(\beta-\alpha)/(\beta z)=F(4a)$\\
$f_3(a)$ & $(\beta-\alpha)/z$ & $k_2(a)$ & $(1+\alpha-2\beta)/(\beta z^2)$\\
\bottomrule
\end{tabular}
\end{table}

%-------------------------------------------------------------------------------
\subsection{Interpolation operators}
\label{app_PB_interpolation}
The strain rate and the body-frame velocity of the cross-section at $u=s/L$ follow from the nodal velocities as
\begin{equation}
\dot{\dv}=\boldsymbol{\Pi}\,\vv^e,\qquad
\vv(u)=\mathbf{Z}(u)\,\vv^e,\qquad
\mathbf{Z}(u)=\bigl[\,\identity_6-\Tm^*(u)\quad\Tm^*(u)\,\bigr],
\label{eq_PB_PiZ}
\end{equation}
which are the operators $\mathbf{P}$ and $\mathbf{Q}(s)$ of \cite{Sonneville2014}.
In the Euler--Rodrigues chart, with $\pV(s)=[\pU\tp(s) \,\, \pO\tp(s)]\tp$ the motion vector at $s=uL$ of \eq{eq_scaling}, they read
\begin{equation}
\boldsymbol{\Pi}=\mathbf{G}(\pV)\bigl[-\Tm^{-1}(\pV)\,\Ad^{-1}\quad\Tm^{-1}(\pV)\bigr],\qquad
\Tm^*(u)=u\,\Tm\bigl(\pV(s)\bigr)\,\mathbf{G}^{-1}\bigl(\pV(s)\bigr)\,\mathbf{G}(\pV)\,\Tm^{-1}(\pV) .
\label{eq_PB_PiZchart}
\end{equation}
Applying \eq{eq_PB_chartIdentity} at $\pV$ and at $\pV(s)$, whose exponential coordinate is $u\dv$, collapses \eq{eq_PB_PiZchart} to
\begin{equation}
\boldsymbol{\Pi}=\bigl[-\boldsymbol{\mathcal{T}}^{-1}(\dv)\,\Ad^{-1}\quad\boldsymbol{\mathcal{T}}^{-1}(\dv)\bigr],\qquad
\Tm^*(u)=u\,\boldsymbol{\mathcal{T}}(u\dv)\,\boldsymbol{\mathcal{T}}^{-1}(\dv) ,
\label{eq_PB_PiZcollapsed}
\end{equation}
so that neither $\Tm$, $\mathbf{G}$, nor $\pV(s)$ is formed in the element arrays.
$\boldsymbol{\mathcal{T}}(u\dv)$ is \eq{eq_PB_strainTangent} evaluated at $u\dv$, i.e., with the scalar functions at $u^2a$.
Eqs.~\eqref{eq_PB_chartT}--\eqref{eq_PB_XUOmega} and \eq{eq_PB_PiZchart} are therefore not needed for the implementation, but they relate the polynomial element to the chart and allow a direct verification of \eq{eq_PB_PiZcollapsed}.

%-------------------------------------------------------------------------------
\subsection{Strain, internal forces, and inertia}
\label{app_PB_forces}
With the constant strain $\boldsymbol{\varepsilon}=(\dv-\dv_0)/L$ of \sect{sect_PolyInterpolationMotion} and its rate $\dot{\boldsymbol{\varepsilon}}=\boldsymbol{\Pi}\,\vv^e/L$, the internal forces are
\begin{equation}
\gb\int^{e}=\boldsymbol{\Pi}\tp\bigl(\mathbf{K}_c\,\boldsymbol{\varepsilon}+\mathbf{C}_c\,\dot{\boldsymbol{\varepsilon}}\bigr),
\label{eq_PB_gint}
\end{equation}
where $\mathbf{K}_c=\operatorname{diag}(EA,GA_2,GA_3,GJ,EI_2,EI_3)$ is the cross-section stiffness and $\mathbf{C}_c$ an optional diagonal matrix of strain-rate (Kelvin--Voigt) damping coefficients.
Both vanish under rigid-body motion.
With the cross-section inertia $\mathbf{M}_c=\operatorname{diag}(\rho_mA\,\identity,\,\mathbf{J}_c)$ per unit length, the mass matrix, the quadratic velocity forces, and the external forces are
\begin{equation}
\Mb^{e}=\intop_0^L\mathbf{Z}\tp\mathbf{M}_c\mathbf{Z}\,ds,\qquad
\gb\ine^{e}=\intop_0^L\mathbf{Z}\tp\Bigl(\mathbf{M}_c\dot{\mathbf{Z}}-\operatorname{ad}\bigl(\mathbf{Z}\vv^e\bigr)^{\!\top}\mathbf{M}_c\mathbf{Z}\Bigr)ds\;\vv^e,\qquad
\gb\ext^{e}=\intop_0^L\mathbf{Z}\tp\bar{\mathbf{f}}\,ds,
\label{eq_PB_inertia}
\end{equation}
where $\bar{\mathbf{f}}$ is the distributed load in the body frame of the cross-section, e.g., $\bar{\mathbf{f}}=[(\rho_mA\,\RotMat\tp(s)\,\mathbf{g})\tp \,\, \zero\tp]\tp$ for gravity $\gv$ with $\RotMat(s)$ from \eq{eq_field}.
Loads constant in time, such as gravity, are assembled together with the internal forces and enter the equilibrium, whereas harmonic loads form $\epsilon\,\gb\ext$ in \eq{eq_DynEOM}.
The rate $\dot{\mathbf{Z}}=[-\dot{\Tm}^*\;\;\dot{\Tm}^*]$ follows from differentiating \eq{eq_PB_PiZcollapsed} with $\Tm_i:=\boldsymbol{\mathcal{T}}^{-1}(\dv)$ and $\dot{\dv}=\Tm_i\Delta\vv$,
\begin{equation}
\dot{\Tm}^*(u)=u\begin{bmatrix}
    \dot{\boldsymbol{\mathcal{T}}}(u\dv, \,u\dot{\dv})\Tm_i-\boldsymbol{\mathcal{T}}(u\dv) \Tm_i & \dot{\boldsymbol{\mathcal{T}}}(\dv;\dot{\dv})\,\Tm_i
    \end{bmatrix} \eqComma
\label{eq_PB_Tstardot}
\end{equation}
where the factor $\Tm_i\dot{\boldsymbol{\mathcal{T}}}(\dv,\dot{\dv})\Tm_i$ does not depend on $u$ and is formed once per element.
The directional derivative $\dot{\boldsymbol{\mathcal{T}}}(\dv,\dot{\dv})$ of \eq{eq_PB_strainTangent} follows from the product rule with $\dot a=\frac12\dv_\Omega\tp\dot{\dv}_\Omega$ and $\dot\varsigma=\dot{\dv}_\Omega\tp\dv_U+\dv_\Omega\tp\dot{\dv}_U$,
\begin{equation}
\begin{aligned}
\dot{\boldsymbol{\mathcal{T}}}_\Omega&=f_1'\dot a\,\vecToSkew{\dv_\Omega}+f_1\dot{\widetilde{\dv}}_\Omega+f_2'\dot a\,\vecToSkew{\dv_\Omega}^{\,2}+f_2\bigl(\dot{\widetilde{\dv}}_\Omega\vecToSkew{\dv_\Omega}+\vecToSkew{\dv_\Omega}\dot{\widetilde{\dv}}_\Omega\bigr),\\
\dot{\boldsymbol{\mathcal{T}}}_U&=f_1'\dot a\,\vecToSkew{\dv_U}+f_1\dot{\widetilde{\dv}}_U+f_2'\dot a\,[\dv_U,\dv_\Omega]+f_2\bigl([\dot{\dv}_U,\dv_\Omega]+[\dv_U,\dot{\dv}_\Omega]\bigr)\\
&\quad+\bigl(f_3'\dot a\,\varsigma+f_3\dot\varsigma\bigr)\vecToSkew{\dv_\Omega}+f_3\varsigma\,\dot{\widetilde{\dv}}_\Omega
+\bigl(f_4'\dot a\,\varsigma+f_4\dot\varsigma\bigr)\vecToSkew{\dv_\Omega}^{\,2}+f_4\varsigma\bigl(\dot{\widetilde{\dv}}_\Omega\vecToSkew{\dv_\Omega}+\vecToSkew{\dv_\Omega}\dot{\widetilde{\dv}}_\Omega\bigr),
\end{aligned}
\label{eq_PB_calTdot}
\end{equation}
which is polynomial in $a$ and linear in $\dot{\dv}$; the lower-left block of $\dot{\boldsymbol{\mathcal{T}}}$ vanishes and its lower-right block equals $\dot{\boldsymbol{\mathcal{T}}}_\Omega$.

%-------------------------------------------------------------------------------
\subsection{Element arrays and assembly}
\label{app_PB_arrays}
The integrals \eq{eq_PB_inertia} are evaluated by Gauss--Legendre quadrature with abscissae $u_k\in(0,1)$ and weights $w_k$, the same stations at which the reciprocals $r_k$ are lifted,
\begin{equation}
\Mb^{e}=L\sum_{k}w_k\,\mathbf{Z}_k\tp\mathbf{M}_c\mathbf{Z}_k,\qquad
\gb\ine^{e}=L\sum_{k}w_k\,\mathbf{Z}_k\tp\Bigl(\mathbf{M}_c\dot{\mathbf{Z}}_k-\operatorname{ad}\bigl(\mathbf{Z}_k\vv^e\bigr)^{\!\top}\mathbf{M}_c\mathbf{Z}_k\Bigr)\vv^e,\qquad
\gb\ext^{e}=L\sum_{k}w_k\,\mathbf{Z}_k\tp\bar{\mathbf{f}}(u_k),
\label{eq_PB_quadrature}
\end{equation}
with $\mathbf{Z}_k=\mathbf{Z}(u_k)$ and $\dot{\mathbf{Z}}_k=\dot{\mathbf{Z}}(u_k)$.
In this work, $n_g=3$.
The integrands are polynomials in $u$ of high degree, since $\boldsymbol{\mathcal{T}}(u\dv)$ contains the minimax polynomials evaluated at $u^2a$.
At a relative rotation of $0.35$\,rad, the mass matrix obtained with three, four, and six points differs from the exactly integrated one by $5\cdot10^{-7}$, $3\cdot10^{-10}$, and $10^{-14}$, respectively, whereas $\boldsymbol{\varepsilon}$ and $\gb\int^{e}$ are constant over the element and require no quadrature.
The element arrays are scattered into $\Mb(\qVec)$, $\gb\ine$, $\gb\int$, and $\gb\ext$ of the mechanical residual \eq{eq_MechanicalResidual}, the nodal kinematic reconstruction of \sect{sect_BeamNodalUnknowns} and the additive reconstruction $\dot{\mathbf{y}}_e=\vv_{\mathbf{y}}$ of the lifted coordinates form the corresponding rows of $\kinRecMat(\qVec)$, and $\mathbf{c}_e$ enters $\phib\jnt$ as described in \sect{sect_PolyInterpolationMotion}, all within the assembled residual \eq{eq_CompleteResidual}.
Algorithm~\ref{alg_PB_element} summarizes the evaluation of one element.

\begin{algorithm}[htb]
\caption{Evaluation of one polynomial geometrically exact beam element.}
\label{alg_PB_element}
\KwIn{Nodal unknowns $\qVec_A,\qVec_B$ and velocities $\vv^e$; lifted scalars $a$, $r$, $r_1,\ldots,r_{n_g}$; element data $L$, $\dv_0$, $\mathbf{K}_c$, $\mathbf{C}_c$, $\mathbf{M}_c$, $\bar{\mathbf{f}}$, $(u_k,w_k)_{k=1}^{n_g}$}
\KwOut{$\Mb^{e}$, $\gb\ine^{e}$, $\gb\int^{e}$, $\gb\ext^{e}$; $\mathbf{c}_e$}
$q_0$ and $\pV$ from \eq{eq_relativeErpFromNodes}; $\rho$, $\mu$ from \eq{eq_bridgeScalars}\;
$\dv$ from \eq{eq_bridge}; $\boldsymbol{\varepsilon}=(\dv-\dv_0)/L$\;
$\RotMat$, $\tv$ from \eq{eq_RotErp}, \eq{eq_Rtheta}; $\Delta\vv$ and $\Ad^{-1}$ from \eq{eq_relVelocity}\;
$\Tm_i=\boldsymbol{\mathcal{T}}^{-1}(\dv)$ from \eq{eq_PB_strainTangent}; $\boldsymbol{\Pi}$ from \eq{eq_PB_PiZcollapsed}; $\dot{\dv}=\Tm_i\Delta\vv$\;
$\gb\int^{e}=\boldsymbol{\Pi}\tp\bigl(\mathbf{K}_c\boldsymbol{\varepsilon}+\mathbf{C}_c\dot{\dv}/L\bigr)$ from \eq{eq_PB_gint}\;
$\Tm_i\dot{\boldsymbol{\mathcal{T}}}(\dv;\dot{\dv})\Tm_i$ from \eq{eq_PB_calTdot}, once per element\;
Identities $\mathbf{c}_e$ from \eq{eq_liftedIdentities}\;
\For{$k=1,\ldots,n_g$}{
    $\boldsymbol{\mathcal{T}}(u_k\dv)$ and $\dot{\boldsymbol{\mathcal{T}}}(u_k\dv;u_k\dot{\dv})$ from \eq{eq_PB_strainTangent} and \eq{eq_PB_calTdot}, scalar functions at $u_k^2a$\;
    $\Tm^*$, $\mathbf{Z}_k$ from \eq{eq_PB_PiZcollapsed}, \eq{eq_PB_PiZ}; $\dot{\Tm}^*$, $\dot{\mathbf{Z}}_k$ from \eq{eq_PB_Tstardot}\;
    \If{$\bar{\mathbf{f}}$ depends on the configuration}{
        $\pV(s_k)$, $q_0(s_k)$ from \eq{eq_scaling}, \eq{eq_scalingScalars}; $\RotMat(s_k)$, $\rv(s_k)$ from \eq{eq_field} with $r_k$\;
    }
    Add the contributions of station $k$ to $\Mb^{e}$, $\gb\ine^{e}$, $\gb\ext^{e}$ according to \eq{eq_PB_quadrature}\;
}
\end{algorithm}

\bibliographystyle{elsarticle-num}

% Loading bibliography database
\bibliography{references}

% Biography
%\bio{}
% Here goes the biography details.
%\endbio

%\bio{pic1}
% Here goes the biography details.
%\endbio

\end{document}